\documentclass[12pt]{amsart}
\usepackage{ShaharStyle}

\title{Lusztig Families of Distinguished Representations for Symmetric Spaces}
\author{Shahar Dagan}
\thanks{Email: \texttt{shahar.dagan@weizmann.ac.il}} 
\begin{document}

\begin{abstract}
    Let $G$ be a connected reductive group over an algebraic closure
    of a finite field of odd characteristic. Let $F$ be a Frobenius
    endomorphism and $\theta$ an involution commuting with $F$.
    Set $H\coloneqq G^\theta$, and write $H^\circ$ for its identity component.
    We prove that every $H^F$-distinguished irreducible character $\chi$
    of $G^F$ belongs to the same Lusztig family as its
    $\theta$-twisted dual $\chi^\vee\circ\theta$.
    If $Z(G)$ is smooth and connected and every simple factor of
    $G_{\opr{der}}$ is of classical type, then every
    $(H^\circ)^F$-distinguished $\chi\in\Irr(G^F)$ satisfies
    $\chi=\chi^\vee\circ\theta$.
    We give an example in symplectic groups showing that this equality can fail
    when the center is disconnected.

    For groups with smooth connected center, we obtain a necessary condition
    for distinction in terms of Weyl-group representations.
    We use the convolution action of the monodromic Hecke category
    on $H^\circ$-equivariant monodromic sheaves on $G/U$,
    where $U$ is a maximal unipotent subgroup.
    The weighted Grothendieck groups give a symmetric Hecke
    module, which we specialize and filter using an operator formula
    from \cite{dagan-monodromic-hecke}.
    The associated graded is a direct sum of representations induced
    from characters of cross-action stabilizers.
    The proofs combine the stabilizer criterion with root-system
    calculations and a reduction to simple and Galois cases.
\end{abstract}
\maketitle
% !TEX root = ../main.tex
\section{Introduction}
\label{sec:introduction}
\setcounter{tocdepth}{2}
\tableofcontents

% !TEX root = ../main.tex
\subsection{Statement of main results}
\label{subsec:Statement-of-main-results}
Our first main result concerns the Lusztig families of distinguished characters.
\begin{introtheorem}[cf. Theorem~\ref{thm:lp-conjecture}]
    \label{intthm:lp-conjecture}
    Let $G$ be a connected reductive group over an algebraic closure
    of a finite field of odd characteristic.
    Let $F$ be a Frobenius endomorphism of $G$, and let $\theta$ be an
    involution of $G$ commuting with $F$.
    Set $H\coloneqq G^\theta$, the subgroup of $\theta$-fixed points.
    Let $\chi$ be an $H^F$-distinguished irreducible character of $G^F$.
    Then $\chi$ and its $\theta$-twisted dual $\chi^\vee\circ\theta$ lie in the same Lusztig family.
\end{introtheorem}

For classical groups with connected center, we obtain equality.
\begin{introtheorem}[cf. Theorem~\ref{thm:classical-twisted-self-duality}]
    \label{intthm:classical-twisted-self-duality}
    Let $G,F,\theta,H$ be as in Theorem~\ref{intthm:lp-conjecture}.
    Suppose that $Z(G)$ is smooth and connected and every simple
    factor of $G_{\opr{der}}$ is of classical type.
    Write $H^\circ$ for the identity component of $H$.
    Every $(H^\circ)^F$-distinguished $\chi\in\Irr(G^F)$ satisfies
    \[
        \chi=\chi^\vee\circ\theta.
    \]
\end{introtheorem}

Theorem~\ref{thm:distinguished-family-uniform-invariance}
gives an intermediate result in every simple connected-center type:
all uniform almost characters in a family containing a distinguished
irreducible character are fixed by twisted duality.
In classical types, scalar products with uniform almost characters
separate irreducible characters, giving equality in simple
connected-center type.
The Galois conclusions in Lemma~\ref{lem:lp-for-galois-type}
and Proposition~\ref{prop:classical-galois-twisted-self-duality}
are consequences of Lusztig's results
\cite[Theorem~1.4(a)--(c)]{lusztig2000g}.
The reduction to simple and Galois cases then gives both main theorems.
For groups with connected center and exceptional simple derived
factors, we prove family invariance but do not establish character equality
in general.
In \S\ref{sec:failure-of-equality-in-disconnected-center}, we give
an example where the center is disconnected and equality fails.

% !TEX root = ../main.tex
\subsection{Background and motivation}
\label{subsec:background-and-motivation}

Let $G$ be a connected reductive group over a non-Archimedean local
field $k$, let $\theta$ be a $k$-rational involution, and put
$H\coloneqq G^\theta$.
An irreducible smooth complex representation $\pi$ of $G(k)$ is
$H(k)$-distinguished if
\[
    \Hom_{H(k)}(\pi,\bC)\ne0.
\]
The Lapid--Prasad conjecture predicts that the $L$-packet containing an
$H(k)$-distinguished irreducible smooth representation is preserved
by $\theta$-twisted duality $\pi\mapsto\pi^\vee\circ\theta$
\cite[Conjecture~1.2]{Gla18}.

Over $p$-adic fields of odd residual characteristic,
Hakim--Murnaghan reduce distinction
to conditions involving cuspidal representations of finite reductive
quotients and quadratic character twists
\cite[Theorem~5.26 and Remark~5.27]{hakim2008distinguished}.
If $G$ splits over a tamely ramified extension and $p$ does not
divide the order of its Weyl group, Fintzen's theorem
\cite[Theorem~8.1]{fintzen2021types} and Kaletha's verification of
Hypothesis~C \cite[Lemma~3.5.2]{kaletha2019regular} make
Hakim--Murnaghan's results applicable to all supercuspidal representations.

For tamely ramified groups over $p$-adic fields with sufficiently large
residual characteristic, Chong Zhang proves twisted self-duality for
distinguished regular depth-zero and epipelagic supercuspidal
representations \cite[Corollaries~3.17 and~3.20]{zhang2020distinguished}.
He also proves invariance of their $L$-packets under
$\pi\mapsto\pi^\vee\circ\theta$
\cite[Corollary~4.18]{zhang2020distinguished}.

For the finite-field problem, let $G$ be a connected reductive group
over $\overline{\bF}_q$, with $q$ odd, and let $F$ be the Frobenius
endomorphism of an $\bF_q$-structure on $G$.
Let $\theta$ be an involution commuting with $F$, and put
$H\coloneqq G^\theta$.
An irreducible character of $G^F$ is $H^F$-distinguished when its
restriction to $H^F$ contains the trivial character.
The question is how distinction constrains the relation between a
character $\chi$ and its twisted dual $\chi^\vee\circ\theta$.

Lusztig computes $K^F$-averages of Deligne--Lusztig virtual characters for
$F$-stable closed subgroups $H^\circ\subseteq K\subseteq H$
\cite[Theorem~3.3]{lusztig2013symmetric}.
He also gives a necessary condition on the torus parameters of
Deligne--Lusztig cohomology containing a $K^F$-distinguished constituent
\cite[Theorem~3.5]{lusztig2013symmetric}.
For the Galois pair $(G^{F^2},G^F)$, with involution induced by $F$,
he determines the multiplicities
of $G^F$-invariants in irreducible representations of $G^{F^2}$ when
$Z(G)$ is connected and $G/Z(G)$ is simple
\cite[Theorem~1.4(a)--(d)]{lusztig2000g}.
Under these hypotheses, the Galois case of our family-invariance theorem follows.

In \cite[Theorem~A]{kapon2025distinguished}, Kapon proved twisted self-duality for every
distinguished irreducible
representation of $\GL_n(\bF_{q^m})$, with $q$ odd and $m\in\{1,2\}$,
relative to a $\bF_q$-rational involution of
$\opr{Res}_{\bF_{q^m}/\bF_q}\GL_n$.

The relative Langlands program provides a broader context for
distinction and its expected symmetries.
For Galois pairs, Prasad proposes descriptions of distinction and
its multiplicities in terms of Langlands parameters, with a
canonical quadratic character twist
\cite[\S13]{prasad2015arelative}.
Yang studies a finite-field analogue through Shintani base change
\cite[Theorem~1.1]{yang2020distinguished}.
On the spectral side, Sakellaridis--Venkatesh formulate conjectures
relating the $L^2$-spectra of spherical varieties to their dual groups
\cite[\S16.2]{sakellaridis2017periodsharmonicanalysisspherical}.
Ben-Zvi--Sakellaridis--Venkatesh propose a duality for a class of
Hamiltonian varieties for Langlands dual groups, relating
automorphic periods to $L$-functions
\cite[\S\S1, 3.5]{ben2024relative}.
Leslie develops a theory of endoscopy for symmetric varieties
using their dual groups and relates this theory to relative
Langlands duality \cite[\S\S7--9]{leslie2024symmetric}.
A finite-field instance of relative Langlands duality appears in
Wang's work relating theta correspondence and Gan--Gross--Prasad
multiplicities through Lusztig's Jordan decomposition
\cite[\S\S1.2, 7]{wang2024lusztig}.

% !TEX root = ../main.tex
\subsection{Further research}
\label{subsec:further-research}
\subsubsection{The Lapid--Prasad conjecture}
\label{subsubsec:lapid-prasad-conjecture}
Hakim--Murnaghan's results
\cite[Theorem~5.26 and Remark~5.29]{hakim2008distinguished}
provide a route from our finite-field results to the original
Lapid--Prasad conjecture for tame supercuspidal representations.
Pursuing this route requires extending our results to the quadratic
distinction characters in their reduction and lifting the resulting
duality relations to the full depth-zero inducing representation.

Using family invariance to establish the packet formulation requires
comparing Lusztig families with $L$-packets compatibly with Yu's
construction and twisted duality.
Aubert--Ciubotaru--Romano provide a model for such a comparison in
the unipotent setting: for $\SL_n$ and $\PGL_n$, they relate
combinations of characters in unipotent $L$-packets to Lusztig's
nonabelian Fourier transform on finite reductive quotients
\cite[Proposition~13.9 and Theorem~14.2]{aubertCiubotaruRomano2025}.
Extending this comparison to the tame supercuspidal data relevant
here remains a further problem.

The equality for connected center gives a more packet-like conclusion
in classical types.
Assume every simple factor of $G_{\opr{der}}$ is of classical type,
and let $\chi\in\Irr(G^F)$ be $H^F$-distinguished.
Put $G_{\opr{ad}}\coloneqq G/Z(G)$, and let $G_{\opr{ad}}^F$
act on $\Irr(G^F)$ by conjugation.
Theorem~\ref{intthm:classical-twisted-self-duality} and the
regular-embedding reduction in
\S\ref{subsubsec:regular-embeddings-and-smooth-coverings}
imply that $G_{\opr{ad}}^F\cdot\chi$ contains
$\chi^\vee\circ\theta$, even when $Z(G)$ is disconnected.
The sum over this orbit is fixed by twisted duality:
\[
    \left(\sum_{\psi\in G_{\opr{ad}}^F\cdot\chi}\psi\right)^\vee
    \circ\theta
    =\sum_{\psi\in G_{\opr{ad}}^F\cdot\chi}\psi.
\]

Beyond supercuspidals, Offen's necessary conditions for distinction
under parabolic induction
\cite[Theorem~4.2 and Corollary~5.2]{offen2017parabolic}
suggest studying how twisted duality of inducing data constrains
the $L$-packets of distinguished irreducible quotients.
This requires accounting for the character twists in those conditions
and controlling the passage from induced representations to their
irreducible quotients.
Related induction arguments appear in our work on multiplicity bounds
\cite[\S\S3--5]{dagan2025boundmultiplicitiessymmetricpairs}.

\subsubsection{Twisted self-duality in exceptional types}
\label{subsubsec:twisted-self-duality-in-exceptional-types}
Theorem~\ref{thm:distinguished-family-uniform-invariance} does not
give twisted self-duality in exceptional types, since uniform
almost characters need not separate irreducible characters.
A possible approach is to refine Proposition~\ref{prop:occurrence}
using an asymptotic analogue of the symmetric Hecke module.
Lusztig's construction for twisted involutions provides a model
\cite[Theorem~0.6]{lusztig2012asymptotic}.

A relative $a$-function could provide a compatible truncation of
the symmetric Hecke module category.
Retaining the categorical structure may yield an occurrence
criterion for individual family parameters, analogous to the
square-root modules in Lusztig's Galois distinction formula
\cite[\S1.3 and Theorem~1.4(c)]{lusztig2000g}.
The remaining goal would be to show that every parameter detected
by this criterion is fixed by twisted duality.
Such a refinement could extend
Theorem~\ref{intthm:classical-twisted-self-duality}
to exceptional types with connected center.

% !TEX root = ../main.tex
\subsection{Overview of the method}
\label{subsec:overview-of-the-method}
Family invariance for $H^F$-distinguished characters transfers along
products, smooth coverings, regular embeddings, and
$\mathsf{FI}$-rational equivalences.
Character equality for $(H^\circ)^F$-distinguished characters transfers
along products, smooth coverings, and the cyclic equivalences of
Lemma~\ref{lem:f-cyclic-rationally-equivalent}.
Proposition~\ref{prop:connected-center-reduction} reduces family
invariance for arbitrary $G$, and character equality for connected
$Z(G)$, to simple connected-center triples and Galois triples arising
from groups of simple connected-center type.

For connected center,
Proposition~\ref{prop:distinguished-uniform-almost-characters}
associates to each distinguished irreducible character $\chi$ a pair
$(\cL,E)$ indexing its Lusztig family.
Here $\cL$ is a Kummer local system on a fixed torus $T$,
and $E$ is a simple module of the relative Weyl group $W_\cL\coloneqq \{w\in W \,|\, w\cL \cong \cL\}$.
The associated uniform almost character $\rho^\cL_{\widetilde E}$ is
$(H^\circ)^F$-distinguished.
Proposition~\ref{prop:occurrence} then implies that $E$ occurs in
the specialized symmetric Hecke module $M^{H^\circ}_{\cL}$.

The module $M^{H^\circ}_{\cL}$ is constructed by convolution of sheaves.
For $B=TU$ and $\fo\coloneqq W\cdot\cL$, the mixed monodromic
Hecke category $\cH_\fo$ of Lusztig--Yun
\cite[\S\S2.7, 3.1]{lusztig2020endoscopy} acts by convolution on
the category $\cM^{H^\circ}_\fo$ of $H^\circ$-equivariant complexes
on $G/U$ with torus monodromy in $\fo$.
Passing to weighted Grothendieck groups gives a Hecke-algebra action,
as in \cite[\S3.13]{lusztig2020endoscopy} and
Henderson's construction \cite[Proposition~3.4]{henderson2003character}.
Specializing the weight variable $u^{1/2}=-1$ and taking the $\cL$-summand gives
$M^{H^\circ}_{\cL}$ with its $W_\cL$-action.
We recall these constructions in
\S\ref{subsec:hecke-category-and-symmetric-module}.

The standard basis $\{m_\xi\}$ of $M^{H^\circ}_{\cL}$ is indexed by
the set $\Gamma^{H^\circ}_{\cL}$ of $H^\circ$-orbits
in $G/B$ equipped with irreducible equivariant local systems on
their inverse images in $G/U$, with torus monodromy $\cL$.
Lemma~\ref{lem:combi-parameter-rel-standard-obj} writes each parameter
as $\xi=(f,\mu)$: the element $f\in G$ represents the orbit
$H^\circ fB/B$, and the component-group character $\mu$ determines
the local system.

\cite{dagan-monodromic-hecke} computes the Hecke operators $\mathbf T_s$
associated to simple reflections of $W$ on the standard basis $\{m_\xi\}$,
including reflections that change monodromy.
The two results recalled in Lemmas~\ref{lem:cross-action-well-defined}
and~\ref{lem:specialized-rank-one-leading-term} give the cross action
of $W$ on the parameter set $\Gamma^{H^\circ}_{\fo}$.
They show that, after specialization, each simple-reflection operator acts by a signed
cross-action term plus terms of larger anti-invariant torus rank
(Definition~\ref{def:theta-degree}).
Lemma~\ref{lem:relative-module-filtration-stability} therefore gives a stable
filtration whose associated graded has basis vectors $\overline m_\xi$.
The group $W_\cL$ permutes the lines $\Span\{\overline m_\xi\}$.

Let $W_\xi$ be the stabilizer of $\xi$ under cross action, and let
$\eta_\xi$ be its character on $\Span\{\overline m_\xi\}$.
Each orbit of basis lines spans a module induced from such a character.
By semisimplicity and Frobenius reciprocity,
Corollary~\ref{cor:relative-module-graded-constituents} says that
$E$ occurs in $M^{H^\circ}_{\cL}$ if and only if there exists $\xi$ such that $E$ is $(W_\xi,\eta_\xi)$-distinguished.
Theorem~\ref{thm:family-distinction-criterion} packages these implications as a
necessary condition for distinction.

For the parameter $\xi=(f,\mu)$ supplied by this criterion, the
conjugate involution
$\theta_f\coloneqq\opr{Int}_{f^{-1}}\circ\theta\circ\opr{Int}_f$
preserves $T$ and induces an automorphism $\delta_f$ of $W_\cL$.
Proposition~\ref{prop:induced-weyl-group-involution}
identifies the family of $\chi^\vee\circ\theta$ with the family
indexed by $E\circ\delta_f$.
By Proposition~\ref{prop:non-fixed-module-criteria}, failure of
$E\cong E\circ\delta_f$ can arise only from constituents that do not
match under $\delta_f$ on exchanged irreducible factors or from degenerate
representations on preserved factors of type $D_{2m}$ where $\delta_f$ acts outerly.
For Galois triples, $\theta_f$ exchanges all ambient root-system
components. Lemma~\ref{lem:diagonal-subgroup-in-parameter-stabilizer}
therefore gives $E^{\Delta_{\delta_f}^{\mathrm{sw}}}\ne0$, and
Proposition~\ref{prop:non-fixed-module-criteria} gives
$E\cong E\circ\delta_f$.
Lemma~\ref{lem:lp-for-galois-type} gives another proof of Galois
family invariance, which follows from Lusztig's
\cite[Theorem~1.4(a),(b)]{lusztig2000g}.
Proposition~\ref{prop:classical-galois-twisted-self-duality}
deduces classical character equality from Lusztig's
parameter formula \cite[Theorem~1.4(c)]{lusztig2000g}
and compatibility of Jordan decomposition with twisted duality.

For simple connected-center triples,
Lemma~\ref{lem:involution-away-from-a-d-e6-are-inner}
shows that away from types $A$, $D$, and $E_6$, $\delta_f$
is inner.
In types $A$, $D$, and $E_6$,
Lemmas~\ref{lem:diagonal-subgroup-in-parameter-stabilizer}
and~\ref{lem:even-sign-subgroup-in-parameter-stabilizer}
exclude the obstructions for $E\cong E\circ\delta_f$
for each of these simple root-system types.
Furthermore,
Lemma~\ref{lem:exchanged-components-have-type-a}
implies that exchanged components have
type $A$, whose families are singletons.
Degenerate $D_{2m}$ representations also form singleton families.
Lemma~\ref{lem:fixed-label-implies-fixed-family} therefore extends
fixation of $E$ to its entire Weyl-group family.
Lemma~\ref{lem:unif-alm-char-theta-inv-functoriality} gives exact
fixation of every uniform almost character in the associated
character family, proving
Theorem~\ref{thm:distinguished-family-uniform-invariance}.

Family invariance follows and proves
Theorem~\ref{intthm:lp-conjecture} after the reductions.
For classical groups with connected center,
Lemma~\ref{lem:classical-uniform-separation} shows that uniform
almost characters separate irreducible characters.
The fixed uniform coefficients therefore imply character equality
in simple connected-center classical type.
Together with the classical Galois cases, the equality-preserving
reduction proves
Theorem~\ref{intthm:classical-twisted-self-duality}.

% !TEX root = ../main.tex
\subsection{Acknowledgments}
\label{subsec:acknowledgments}
I thank Guy Shtotland, Guy Kapon, and Shachar Carmeli for helpful conversations.
I thank my advisor, Avraham Aizenbud, for introducing
me to the problem, for his guidance throughout this work, and for
teaching me mathematics.

I was partially supported by ISF grant 1781/23 and BSF grant 2022193.

% !TEX root = ../main.tex
\section{Preliminaries}
\label{sec:preliminaries}
% !TEX root = ../main.tex
\subsection{Conventions and distinction}
\label{subsec:conventions-and-distinction}

Let $G$ be a connected reductive group over $k\coloneqq\overline{\bF}_q$,
where $q$ is odd. Fix a Frobenius endomorphism $F:G\to G$ and an
involution $\theta$ commuting with $F$. Put $H\coloneqq G^\theta$,
and write $H^\circ$ for its identity component; both are $F$-stable.
Fix a prime $\ell\ne\opr{char}(k)$ and a compatible identification
of roots of unity of order prime to $\opr{char}(k)$ in $k$ and
$\overline{\bQ}_\ell$.

Write $\opr{CF}(G^F)$ for the space of $\overline{\bQ}_\ell$-valued
class functions on $G^F$, and $\Irr(G^F)$ for its irreducible characters.
For $a\in\opr{CF}(G^F)$, put $a^\vee(g)\coloneqq a(g^{-1})$.
If $V$ has character $\chi_V$, then $\chi_V^\vee$ is the character
of the contragredient representation $V^\vee$.
A reductive group has \emph{connected center} if its scheme-theoretic
center is smooth and connected.

\begin{definition}
    \label{def:qf-distinction}
    Let $Q\subseteq G$ be an $F$-stable closed subgroup. Define the
    averaging functional
    \[
        \opr{Av}_{Q^F}:\opr{CF}(G^F)\longrightarrow\overline{\bQ}_\ell,
        \qquad
        \opr{Av}_{Q^F}(f)\coloneqq\frac{1}{|Q^F|}\sum_{h\in Q^F}f(h).
    \]
    A class function $f\in\opr{CF}(G^F)$ is \emph{$Q^F$-distinguished}
    if $\opr{Av}_{Q^F}(f)\ne0$.
    We use the same notation for the linear extension to class functions
    with values in any $\overline{\bQ}_\ell$-algebra.
\end{definition}

For a finite-dimensional $\overline{\bQ}_\ell$-representation $V$
of $G^F$ with character $\chi_V$,
\[
    \opr{Av}_{Q^F}(\chi_V)=\dim V^{Q^F}.
\]
Since $V^{H^F}\subseteq V^{(H^\circ)^F}$, $H^F$-distinction
implies $(H^\circ)^F$-distinction for genuine characters.
This implication need not hold for virtual characters.

For a finite group $K$, a subgroup $L\subseteq K$, and a character
$\eta:L\to\overline{\bQ}_\ell^\times$, a finite-dimensional
$\overline{\bQ}_\ell$-representation $E$ of $K$ is
\emph{$(L,\eta)$-distinguished} if $\Hom_L(\eta,E)\ne0$.
We call $E$ \emph{$L$-distinguished} if it is $(L,\opr{triv})$-distinguished
where $\opr{triv}$ is the trivial character.

% !TEX root = ../main.tex
\subsection{Sheaves and weighted Grothendieck groups}
\label{subsec:sheaves}
For a separated scheme $X$ of finite type over $k$,
write $\cD(X)\coloneqq D_c^b(X,\overline{\bQ}_\ell)$
for the bounded constructible derived category of $\ell$-adic sheaves on $X$.
For an algebraic group $Q$ acting on $X$,
we write $\cD^Q(X)$ for the equivariant
constructible derived category as in
\cite[\S1.10.2]{lusztig2020endoscopy}.

\begin{definition}[cf. {\cite[\S\S2.1--2.2]{lusztig2020endoscopy}}]
    \label{def:kummer-local-systems}
    A \emph{Kummer local system} is a rank-one
    $\overline{\bQ}_{\ell}$-local system $\cL$ on a torus $T$
    equipped with a multiplicative isomorphism
    $m_T^*\cL\simeq\cL\boxtimes\cL$ and a trivialization of $\cL_1$,
    satisfying the associativity and unit axioms.
    Here $m_T:T\times T\to T$ is multiplication and $1\in T$ is the identity.
    We require $[n]^*\cL$ to be trivial for some $n>0$ prime to
    $\opr{char}(k)$, where $[n](t)\coloneqq t^n$.
\end{definition}
Choose $n$ as in Definition~\ref{def:kummer-local-systems} and set
$T[n]\coloneqq[n]^{-1}(1)$. There is a unique character
$\chi_{\cL}:T[n]\to\overline{\bQ}_\ell^\times$ such that
\[
    \cL\simeq([n]_*\overline{\bQ}_\ell)_{\chi_\cL}.
\]
Under the fixed identification of roots of unity, write
$\chi_{\cL}=\lambda|_{T[n]}$ with $\lambda\in X^*(T)$, determined
modulo $nX^*(T)$.
The local system $\cL$ corresponds to the rational character
class $\lambda/n+X^*(T)$ in $(X^*(T)\otimes_\bZ\bQ)/X^*(T)$.

\begin{definition}[cf. {\cite[\S2.5]{lusztig2020endoscopy}}]
    \label{def:monodromic}
    Let $a:T\times X\to X$ be a torus action and $\cL$ a Kummer
    local system. Choose $n,\chi_{\cL}$ as above, and put
    $a_n(t,x)\coloneqq a(t^n,x)$. Define
    \[
        \cD_{\cL}(X)\coloneqq
        \{K\in\cD^T(X,a_n)\mid
        z\text{ acts as }\chi_{\cL}(z)\opr{id}_K
        \text{ for all }z\in T[n]\}.
    \]
    Its objects are
    \emph{$(T,\cL)$-monodromic complexes}, with twisted equivariance
    $a^*K\simeq\cL\boxtimes K$.
\end{definition}

\begin{definition}[cf. {\cite[\S5.1]{BBD82}}]
    \label{def:weil-structure}
    Let $F_X:X\to X$ define a model $X_0$ of $X$ over $\bF_{q^r}$,
    where $r\in\bZ_{>0}$.
    A \emph{Weil structure} on $K\in\cD(X)$ is an
    isomorphism $\phi:F_X^*K\xrightarrow{\sim}K$.

    Put $\cD^{\mathrm{mix}}(X,F_X)\coloneqq
        D_m^b(X_0,\overline{\bQ}_\ell)$, the full subcategory of
    $D_c^b(X_0,\overline{\bQ}_\ell)$ whose cohomology sheaves admit
    finite filtrations with pointwise pure quotients of integral weights.
    A constructible sheaf $V$ on $X_0$ is \emph{pointwise pure of weight}
    $j\in\bZ$ if, at every closed point $x$ of degree $d$ over $\bF_{q^r}$,
    each local geometric Frobenius eigenvalue $\alpha$ on $V_{\bar x}$
    is algebraic over $\bQ$ and satisfies
    \[
        |\sigma(\alpha)|=q^{rdj/2}
        \qquad\text{for every embedding }\sigma:\bQ(\alpha)\hookrightarrow\bC.
    \]
\end{definition}

For $V\in \cD^{\mathrm{mix}}(\Spec{\bF_{q^r}})$ concentrated
in degree $0$,
write $\mathrm{Wt}_\bullet V$ for its weight filtration,
with $\opr{Gr}^{\mathrm{Wt}}_jV$ pure of weight $j$.
Put $\cA\coloneqq\bZ[u^{1/2},u^{-1/2}]$ and
$R_{\mathrm{Fr}}(q^r)\coloneqq
    K_0(\cD^{\mathrm{mix}}(\Spec{\bF_{q^r}}))$, with multiplication
induced by tensor product. Define the ring homomorphism
\[
    \opr{wt}:R_{\mathrm{Fr}}(q^r)\longrightarrow\cA,\qquad
    [V]\longmapsto\sum_{i,j}(-1)^i(-u^{1/2})^j
    \dim\opr{Gr}^{\mathrm{Wt}}_jH^i(V).
\]

\begin{definition}[cf. {\cite[\S\S1.10.1--1.10.2]{lusztig2020endoscopy}}]
    \label{def:tate-twists-and-normalization}

    Let $\mathcal C^{\mathrm{mix}}$ be a triangulated category with an
    action of $\cD^{\mathrm{mix}}(\Spec{\bF_{q^r}})$ exact in each variable.
    Its \emph{weighted Grothendieck group} and the \emph{weighted class}
    of $V\in\mathcal C^{\mathrm{mix}}$ are
    \[
        K_0^{\mathrm{wt}}(\mathcal C^{\mathrm{mix}})
        \coloneqq\cA\otimes_{R_{\mathrm{Fr}}(q^r),\opr{wt}}
        K_0(\mathcal C^{\mathrm{mix}}),\qquad
        \chi_u(V)\coloneqq1\otimes[V].
    \]
    % The mixed equivariant and monodromic categories used below carry this scalar action.

\end{definition}
For the equivariant monodromic categories
$\mathcal C^{\mathrm{mix}}=\cD^{\mathrm{mix},Q}_\cL(X,F_X)$ used below,
the commuting $Q$- and $T$-actions have finitely many $Q\times T$-orbits.
Assume that $F_X$ preserves every orbit and $F_X^*$ fixes every isomorphism
class of irreducible $Q$-equivariant $(T,\cL)$-monodromic local systems
on those orbits.
The weight-polynomial $\chi_u$ description of
\cite[pp.~27--28]{henderson2003character} identifies the weighted Grothendieck group
with the free $\cA$-module on simple perverse objects,
equipped with pure weight-zero Weil structures.
Base change preserves these classes and their weights. Hence base change
to $\bF_{q^{rm}}$, with $m\geq1$, induces an $\cA$-linear
isomorphism
\[
    K_0^{\mathrm{wt}}\!\left(\cD^{\mathrm{mix},Q}_\cL(X,F_X)\right)
    \xrightarrow{\sim}
    K_0^{\mathrm{wt}}\!\left(\cD^{\mathrm{mix},Q}_\cL(X,F_X^m)\right).
\]
For an irreducible equivariant monodromic local system $\mathscr E$
on one orbit, two compatible pure weight-zero Weil structures
$\phi,\phi'$ differ by a scalar of weight zero. Thus
$\chi_u(\mathscr E,\phi)=\chi_u(\mathscr E,\phi')$ in the weighted
Grothendieck group of that orbit.

For a finite group $L$ acting trivially on $\Spec{\bF_{q^r}}$,
the natural identification
\[
    K_0^{\mathrm{wt}}\bigl(\cD^{\mathrm{mix},L}(\Spec{\bF_{q^r}})\bigr)
    \simeq\cA\otimes_\bZ K_0(\opr{Rep}_{\overline{\bQ}_\ell}(L))
\]
expresses the weighted class of a mixed complex $V$ with an $L$-action
commuting with Frobenius as the \emph{representation-valued weight polynomial}
\[
    \chi_u(V)=\sum_{i,j}(-1)^i(-u^{1/2})^j
    [\opr{Gr}^{\mathrm{Wt}}_jH^i(V)].
\]
For $L=1$, this is the scalar weight polynomial $\opr{wt}([V])$.

Fix a square root $q^{1/2}\in\overline{\bQ}_\ell$ and let $\overline{\bQ}_\ell(1/2)$ be
the one-dimensional sheaf on $\Spec{\bF_{q^r}}$ with Frobenius
eigenvalue $(q^{1/2})^{-r}$. For $n\in\bZ$, put
\[
    V(n/2)\coloneqq V\otimes\overline{\bQ}_\ell(1/2)^{\otimes n},
    \qquad V\langle n\rangle\coloneqq V[n](n/2).
\]
Here $[n]$ is the cohomological shift and $V(1)$ is the usual Tate twist.

An exact functor $F:\mathcal C^{\mathrm{mix}}\to(\mathcal C')^{\mathrm{mix}}$
compatible with the scalar action induces an $\cA$-linear map $[F]$
on weighted Grothendieck groups, with $[F](\chi_u(V))=\chi_u(FV)$.

Triangle additivity and the weight of the half Tate twist give
\begin{equation}
    \label{eq:weighted-shift-conventions}
    \begin{split}
        \chi_u(V[1])              & =-\chi_u(V),\qquad
        \chi_u(V(1/2))=-u^{-1/2}\chi_u(V),              \\
        \chi_u(V\langle n\rangle) & =u^{-n/2}\chi_u(V).
    \end{split}
\end{equation}
Evaluation at $u^{1/2}=-1$ sends a Frobenius vector space to its
dimension, and a complex on a point to its Euler characteristic.

For $V$ with an action of a finite group $L$ commuting with Frobenius,
taking characters identifies the representation-valued $\chi_u(V)$ with
a polynomial-valued class function. Taking dimensions recovers the scalar
weight polynomial. Exactness of $L$-invariants gives
$\opr{Av}_L(\chi_u(V))=\chi_u(V^L)$.

% !TEX root = ../main.tex
\subsection{Lusztig theory}
\label{subsec:lusztig-theory}
In this subsection we recall definitions and results from
Lusztig theory on the classification of irreducible characters of
a connected reductive group $G$ over a finite field $\bF_q$, following
\cite[Chapters~2--4]{lusztig1984characters}.
For this subsection, we assume that the center $Z(G)$ is connected.
For $Z(G)$ arbitrary, we use the
definition of Lusztig families in \cite[\S5]{digneMichel1990},
which agrees with Definition~\ref{def:lusztig-families-of-irreducible-characters}
when $Z(G)$ is connected.
In \S\ref{sec:reduction-to-simple-and-galois-cases},
we reduce the main theorem to the connected-center case.

Fix an $F$-stable Borel subgroup $B$ and an $F$-stable,
$\theta$-stable maximal torus $T\subseteq B$.
Such a choice exists by
\cite[Lemma~3.1.1]{helminck1994symmetric}.
Put $U\coloneqq R_u(B)$ and $W\coloneqq N_G(T)/T$.
Let $\Phi^+$, $S$, and $\ell_W$ be the positive roots, simple
reflections, and length determined by $B$.
For $s\in S$, write $\alpha_s$ for the positive simple root with
$s=s_{\alpha_s}$.

Let $R^G_T(\psi)$ be the Deligne--Lusztig virtual character,
as defined in \cite[\S\S1.17--1.20]{deligneLusztig1976}.
This virtual character has the following functoriality.

\begin{lemma}[cf. {\cite[\S1.20 and Theorem~4.2]{deligneLusztig1976}}]
    \label{lem:dl-functoriality}
    Let $\varphi:G\to G$ be an algebraic automorphism commuting
    with $F$. Let $T'\subseteq G$ be an $F$-stable maximal torus
    and $\psi:(T')^F\to\overline{\bQ}_\ell^\times$ a character. Then
    \[
        R_{T'}^G(\psi)\circ\varphi
        =R_{\varphi^{-1}(T')}^G(\psi\circ\varphi),
        \qquad
        \bigl(R_{T'}^G(\psi)\bigr)^\vee=R_{T'}^G(\psi^{-1}).
    \]
\end{lemma}

For every character $\psi:T^F\to\overline{\bQ}_\ell^\times$,
there is a unique, up to isomorphism, $F$-stable Kummer
local system $\cL_\psi$ on $T$ whose normalized Frobenius
isomorphism $\phi_\psi:F^*\cL_\psi\xrightarrow{\sim}\cL_\psi$ satisfies
\[
    \psi(t)=\Tr((\phi_\psi)_t,(\cL_\psi)_t),
    \qquad t\in T^F.
\]
Here the normalization is $(\phi_\psi)_1=\mathrm{id}$;
see \cite[\S2.2]{lusztig2020endoscopy}.

\begin{definition}[cf. {\cite[\S2.4]{lusztig2020endoscopy}}]
    \label{def:relative-weyl-group}
    For a Kummer local system $\cL$ on $T$, define
    \[
        W'_\cL\coloneqq\{w\in W\mid w\cL\cong\cL\},
        \qquad
        \Phi_\cL\coloneqq\{\alpha\in\Phi(G,T)\mid
        (\alpha^\vee)^*\cL\cong\overline{\bQ}_\ell\},
    \]
    where $w\cL\coloneqq(w^{-1})^*\cL$.
    The set $\Phi_\cL$ is a root subsystem.
    Set
    \[
        \Phi^+_\cL\coloneqq\Phi_\cL\cap\Phi^+,\qquad
        W^\circ_\cL\coloneqq\langle s_\alpha\mid\alpha\in\Phi_\cL\rangle.
    \]
    Write $S_\cL$ for the simple reflections of $W^\circ_\cL$
    determined by $\Phi^+_\cL$.
\end{definition}
In general, $S_\cL$ is not a subset of $S$.

\begin{lemma}[cf. {\cite[Theorem~5.13]{deligneLusztig1976}}]
    \label{lem:connected-center-character-weyl-group}
    If $G$ has connected center, then
    \[
        W'_\cL=W^\circ_\cL=W(\Phi_\cL),
    \]
    for every Kummer local system $\cL$ on $T$.
\end{lemma}

Under our connected-center assumption, write
$W_\cL\coloneqq W'_\cL=W^\circ_\cL$ from now on.

Suppose that $W\cdot\cL$ is $F$-stable. Following
\cite[\S2.1]{lusztig1984characters}, define
\[
    Z_{\cL,F}\coloneqq\{z\in W\mid z\cL\cong F^*\cL\}.
\]
Since $G$ has connected center, $Z_{\cL,F}=vW_\cL$, where $v$
is its unique element of minimal $\ell_W$-length
\cite[\S2.19]{lusztig1984characters}.

\begin{definition}
    \label{lem:kummer-to-rational-character}
    Let $z\in W$. Choose a representative $\dot z\in N_G(T)$
    and $g_z\in G$ with $g_z^{-1}F(g_z)=\dot z$. Put
    \[
        T_z\coloneqq g_zTg_z^{-1},\qquad
        \cL_z\coloneqq(\opr{Int}_{g_z^{-1}})^*(F^*\cL).
    \]
    The torus $T_z$ is $F$-stable, and $\cL_z$ is $F$-stable
    if and only if $z\in Z_{\cL,F}$. In this case, let
    $\phi_z:F^*\cL_z\xrightarrow{\sim}\cL_z$ be normalized by
    $(\phi_z)_1=\mathrm{id}$, and define
    \[
        \psi_z(t)\coloneqq\Tr((\phi_z)_t,(\cL_z)_t),
        \qquad t\in T_z^F.
    \]
    By \cite[\S2.2]{lusztig2020endoscopy}, $\psi_z$ is a character
    of $T_z^F$.
\end{definition}
The pair $(T_z,\psi_z)$ is independent of $\dot z$ and $g_z$
up to $G^F$-conjugacy by Lang's theorem.
Put $F_z\coloneqq\opr{Int}_{\dot z}\circ F$ on $T$.
The map $\opr{Int}_{g_z}:T\to T_z$ intertwines $F_z$ and $F$.
The Frobenius-stability criterion follows from
\[
    F_z^*(F^*\cL)\cong F^*\cL
    \quad\Longleftrightarrow\quad
    z\cL\cong F^*\cL.
\]

Normalized trace characters commute with pullback along
$F$-equivariant homomorphisms and with duality.
In particular, $\theta^*\cL_z^{-1}$ on $\theta(T_z)$ has
character $\psi_z^{-1}\circ\theta$.

\begin{remark}[cf. {\cite[\S2.2]{lusztig2020endoscopy}}]
    \label{rem:uniform-kummer-exponent}
    Throughout, put $p\coloneqq\opr{char}(k)$ and
    $N\coloneqq |G^F|_{p'}$, the largest divisor of $|G^F|$
    prime to $p$.
    If $\cL$ is a Kummer local system on $T$ with $F$-stable
    $W$-orbit, then
    \[
        [N]^*\cL\cong\overline{\bQ}_\ell.
    \]
\end{remark}

Write $F_W$ for the automorphism of $W$ induced by $F$. Set
\begin{equation}
    \label{eq:transporter-automorphism}
    \gamma_v:W_\cL\longrightarrow W_\cL,
    \qquad \gamma_v(w)\coloneqq F_W(vwv^{-1}).
\end{equation}
The automorphism $\gamma_v$ preserves $S_\cL$ by
\cite[(2.15.1)--(2.15.2)]{lusztig1984characters}.
Following \cite[\S\S3.1--3.2 and 3.6]{lusztig1984characters},
form $W_\cL\rtimes_{\gamma_v}\bZ$, where $(1,1)$ acts on
$W_\cL$ by $\gamma_v$. Let $\Irr(W_\cL)^{\gamma_v}$
consist of $\gamma_v$-stable $E\in\Irr(W_\cL)$. For such $E$,
choose an extension $\widetilde E$ such that
$\Gamma_{\widetilde E}\coloneqq\widetilde E(1,1)$ has finite order. Thus
\[
    \Gamma_{\widetilde E}E(w)\Gamma_{\widetilde E}^{-1}=E(\gamma_v(w)).
\]

\begin{definition}[cf. {\cite[(3.7.1)]{lusztig1984characters}}]
    \label{def:uniform-almost-character}
    The \emph{uniform almost character} associated with
    $\widetilde E$ is
    \[
        \rho^\cL_{\widetilde E}\coloneqq
        \frac{1}{|W_\cL|}\sum_{w\in W_\cL}
        \opr{Tr}(\Gamma_{\widetilde E}E(w))R_{T_{vw}}^G(\psi_{vw}).
    \]
\end{definition}
By Schur's lemma, for two extensions $\widetilde E$ and
$\widetilde E'$, $\Gamma_{\widetilde E}=\zeta \cdot \Gamma_{\widetilde E'}$ for some $\zeta\in \overline{\bQ}_{\ell}$; see also \cite[\S0.4]{digneMichel1990}.
As both $\Gamma_{\widetilde E}$ and $\Gamma_{\widetilde E'}$ have finite order, $\zeta$ is a root of unity.
Therefore, for a fixed $E\in \Irr(W_\cL)^{\gamma_v}$,
$\rho^\cL_{\widetilde E}$ is well defined up to a root of unity.

Use the standard inner product in $\opr{CF}(G^F)$.
For $E,E'\in\Irr(W_\cL)^{\gamma_v}$, consider the relation
\begin{equation}
    \label{eq:family-defining-relation}
    E\sim E' \iff
    \exists\chi\in\Irr(G^F):\quad
    \langle\chi,\rho^\cL_{\widetilde E}\rangle_{G^F}\ne0
    \quad\text{and}\quad
    \langle\chi,\rho^\cL_{\widetilde E'}\rangle_{G^F}\ne0.
\end{equation}

\begin{definition}
    \label{def:weyl-group-families}
    The equivalence classes for the equivalence relation
    generated by \eqref{eq:family-defining-relation}
    are called \emph{Lusztig families of $\gamma_v$-stable representations}.
    For $E\in\Irr(W_\cL)^{\gamma_v}$, write
    $\opr{Fam}^{W_\cL}_{\gamma_v}(E)$ for the class containing $E$.
\end{definition}

\begin{definition}
    \label{def:lusztig-families-of-irreducible-characters}
    For $E\in\Irr(W_\cL)^{\gamma_v}$, define
    the \emph{Lusztig family} of irreducible $G^F$-characters
    associated to $E$ by
    \[
        \opr{Fam}^{G^F}_\cL(E)\coloneqq
        \left\{\chi\in\Irr(G^F)\;\middle|\;
        \begin{array}{l}
            \langle\chi,\rho^\cL_{\widetilde E'}\rangle_{G^F}\ne0
            \\ \text{for some }\gamma_v\text{-stable }
            E'\in\opr{Fam}^{W_\cL}_{\gamma_v}(E)
        \end{array}\right\}.
    \]
\end{definition}

Our exposition uses Kummer local systems on $T$ to describe
Lusztig series and their families. In the dual-group formulation
(see e.g. \cite[\S\S1--2]{digneMichel1990}), Lusztig series are often
indexed by semisimple conjugacy classes in $(G^*)^F$.
Here $G^*$ is the dual of $G$, and $F$ also denotes the dual Frobenius.
Under our connected-center assumption, duality for tori identifies
each $F$-stable orbit $W\cdot\cL$ with such a class.

If $s\in(G^*)^F$ represents the class associated to $W\cdot\cL$,
then the Lusztig series indexed by $s$ consists of the irreducible
characters occurring in $R_{T_z}^G(\psi_z)$ for some
$z\in Z_{\cL,F}$. Thus every family
$\opr{Fam}^{G^F}_\cL(E)$ belongs to the series indexed by $s$.
This identifies our local-system labels with the semisimple
labels used in the dual-group formulation.

Definitions~\ref{def:weyl-group-families}
and~\ref{def:lusztig-families-of-irreducible-characters}
also recover Lusztig's family decomposition within each series.
For $E\in\Irr(W_\cL)^{\gamma_v}$, let $\cF_W$ be its family
in the sense of \cite[\S4.2]{lusztig1984characters}. Then
\[
    \opr{Fam}^{W_\cL}_{\gamma_v}(E)
    =\cF_W\cap\Irr(W_\cL)^{\gamma_v},
\]
and $\opr{Fam}^{G^F}_\cL(E)$ is the character family attached to
$\cF_W$ by \cite[Theorem~4.23]{lusztig1984characters}.
The comparison follows from Lusztig's multiplicity formula:
coefficients vanish between distinct families, while the uniform
almost character attached to the special representation of $\cF_W$
has nonzero coefficient at every character in the attached family
\cite[(4.1.4), (4.14.2), and \S4.26]{lusztig1984characters}.

\begin{lemma}[cf. {\cite[Proposition~6.3]{digneMichel1990}}]
    \label{lem:classical-unipotent-uniform-separation}
    Let $G$ be a connected reductive group with Frobenius $F$.
    Assume that every simple
    factor of $G_{\opr{der}}$ has classical type.
    Let $\chi,\chi'\in\Irr(G^F)$ be unipotent.
    If $\langle \chi, \rho^\cL_{\widetilde E} \rangle = \langle \chi', \rho^\cL_{\widetilde E} \rangle$
    for every uniform almost character $\rho^\cL_{\widetilde E}$ with
    $\cL\cong \underline{\overline{\bQ}_\ell}$, then $\chi=\chi'$.
\end{lemma}
\begin{proof}
    Pullback from $G/Z(G)$ identifies unipotent characters and
    preserves their uniform coefficients, so we may assume that
    $G$ is adjoint.
    By Lusztig's multiplicity formula
    \cite[Theorem~4.23]{lusztig1984characters}, the uniform
    coefficients are the pairings considered in
    \cite[Proposition~6.3]{digneMichel1990}.

    The reductions in \cite[\S\S6.1--6.2]{digneMichel1990}
    account for products and Frobenius permutations of components.
    For classical components, the family groups are of the form $C_2^{r}$
    where $C_2$ is the group with $2$ elements and $r\in \bZ_{\geq0}$.
    Thus none of the exceptional cases excluded in
    \cite[Proposition~6.3]{digneMichel1990} occurs.
    The uniform coefficients therefore determine the unipotent
    character, which proves $\chi=\chi'$.
\end{proof}

\begin{lemma}
    \label{lem:classical-uniform-separation}
    Assume that $Z(G)$ is connected and every simple
    factor of $G_{\opr{der}}$ has classical type.
    If $\chi,\chi'\in\Irr(G^F)$ satisfy
    \begin{equation}
        \label{eq:uniform-coefficients}
        \langle\chi,\rho^\cL_{\widetilde E}\rangle_{G^F}
        =\langle\chi',\rho^\cL_{\widetilde E}\rangle_{G^F}
    \end{equation}
    for every uniform almost character, then $\chi=\chi'$.
\end{lemma}
\begin{proof}
    The vanishing of uniform coefficients between distinct families
    and the nonzero special coefficients imply that $\chi$ and
    $\chi'$ belong to the same Lusztig family, hence the same
    Lusztig series.
    Let $s\in(G^*)^F$ represent the semisimple conjugacy class
    associated to their common Lusztig series.
    Since $Z(G)$ is connected, $(G^*)_{\opr{der}}$ is simply
    connected and $C_{G^*}(s)$ is connected.
    Its root system has only classical components.

    Jordan decomposition identifies the coefficients in
    \eqref{eq:uniform-coefficients} with the uniform coefficients
    of the unipotent characters associated to $\chi$ and $\chi'$,
    up to a common sign
    \cite[Theorem~7.1(i)]{digneMichel1990}.
    Lemma~\ref{lem:classical-unipotent-uniform-separation}
    gives equality of the Jordan correspondents.
    Injectivity of Jordan decomposition gives $\chi=\chi'$.
\end{proof}

% !TEX root = ../main.tex
\section{Reduction to simple and Galois cases}
\label{sec:reduction-to-simple-and-galois-cases}
Throughout this section we use the following definitions.
\begin{definition}\label{def:frobenius-involution-triple}
    Let $\mathsf{FI}$ denote the category of \emph{Frobenius--involution triples}. Its objects are triples $(G,F,\theta)$
    in which $G$ is a connected reductive group, $F:G\to G$ is a
    Frobenius endomorphism, and $\theta:G\to G$ is an involution commuting with $F$.
    A morphism from $(G,F,\theta)$ to $(G',F',\theta')$ is a homomorphism of algebraic groups $f:G\to G'$ such that $f\circ F=F'\circ f$ and $f\circ\theta=\theta'\circ f$.
\end{definition}

\begin{definition}
    \label{def:simple-connected-center-type}
    A connected reductive group $G$ is of \emph{simple connected-center type}
    if $Z(G)$ is smooth and connected and $G_{\opr{der}}$ is simple
    and simply connected.
    An object $(G,F,\theta)\in\mathsf{FI}$ is of simple connected-center type if $G$ is.
\end{definition}

\begin{definition}
    \label{def:galois-type-triple}
    A triple $(G,F,\theta)\in\mathsf{FI}$ is of \emph{Galois type}
    if there is a connected reductive group $G'$ with Frobenius $F'$
    such that
    \[
        (G,F,\theta)\cong
        (G'\times G',(F'\times F')\circ\opr{swap},\opr{swap}),
        \qquad \opr{swap}(x,y)\coloneqq(y,x).
    \]
    We say that $(G,F,\theta)$ \emph{arises from a group of simple connected-center type}
    if $G'$ can be chosen to be of simple connected-center type.
\end{definition}
Proposition~\ref{prop:connected-center-reduction}(1) reduces
Theorem~\ref{intthm:lp-conjecture} to triples of simple
connected-center type and Galois triples arising from groups of
simple connected-center type. We first establish transfer laws and
settle the torus and group cases, then prove the reduction.
Proposition~\ref{prop:connected-center-reduction}(2) gives the
reduction for character equality for groups with connected center.
% !TEX root = ../main.tex
\subsection{Transfer laws}
\label{subsec:transfer-laws}
We establish transfer laws for family invariance under $H^F$-distinction
and character equality under $(H^\circ)^F$-distinction.
Following \cite{taylor2019rootdata}, we construct
regular embeddings and smooth coverings compatibly with involutions.

\subsubsection{Products}
\label{subsubsec:products}
\begin{lemma}
    \label{lem:transfer-lp-along-products}
    \begin{enumerate}
        \item  Let $(G,F,\theta),(G',F',\theta')\in\mathsf{FI}$.
              Set $H\coloneqq G^\theta$ and $H'\coloneqq G'^{\theta'}$.
              If $G$ and $G'$ have smooth connected centers and both triples
              satisfy Theorem~\ref{intthm:lp-conjecture}, then
              $(G\times G',F\times F',\theta\times\theta')$ satisfies
              Theorem~\ref{intthm:lp-conjecture}.
        \item  If $\alpha^\vee\circ\theta=\alpha$ and
              $(\alpha')^\vee\circ\theta'=\alpha'$ for every
              $(H^\circ)^F$-distinguished $\alpha\in\Irr(G^F)$ and every
              $(H'^\circ)^{F'}$-distinguished $\alpha'\in\Irr(G'^{F'})$, then
              $\chi^\vee\circ(\theta\times\theta')=\chi$ for every
              $(H^\circ)^F\times(H'^\circ)^{F'}$-distinguished
              $\chi\in\Irr(G^F\times G'^{F'})$.
    \end{enumerate}

\end{lemma}
\begin{proof}
    Every irreducible character
    of $G^F\times G'^{F'}$ has the form $\chi=\alpha\boxtimes\alpha'$, and
    \[
        \opr{Av}_{H^F\times H'^{F'}}(\chi)
        =\opr{Av}_{H^F}(\alpha)\cdot\opr{Av}_{H'^{F'}}(\alpha').
    \]
    If $\chi$ is $H^F\times H'^{F'}$-distinguished, then $\alpha$ and
    $\alpha'$ are $H^F$- and $H'^{F'}$-distinguished, respectively.
    Moreover,
    \[
        \chi^\vee\circ(\theta\times\theta')
        =(\alpha^\vee\circ\theta)\boxtimes
        ((\alpha')^\vee\circ\theta').
    \]
    In (1),
    $\alpha$ and $\alpha^\vee\circ\theta$ belong to the same family,
    as do $\alpha'$ and $(\alpha')^\vee\circ\theta'$.
    Families of the product are products of families by
    Definition~\ref{def:lusztig-families-of-irreducible-characters},
    which proves the family assertion.

    Since
    $((G\times G')^{\theta\times\theta'})^\circ=H^\circ\times H'^\circ$,
    the same averaging formula with $H,H'$ replaced by $H^\circ,H'^\circ$
    proves (2).
\end{proof}

\subsubsection{Rational equivalences}
\label{subsubsec:rational-equivalences}
\begin{definition}
    \label{def:rational-equivalence}
    Pairs $(G,F)$ and $(G',F')$ consisting of a connected reductive
    group and a Frobenius endomorphism are \emph{rationally equivalent}
    if there is an isomorphism
    \[
        \phi:G^F\xrightarrow{\sim}G'^{F'}
    \]
    whose induced bijection on irreducible characters preserves and
    reflects Lusztig families. Triples
    $(G,F,\theta),(G',F',\theta')\in\mathsf{FI}$ are
    \emph{$\mathsf{FI}$-rationally equivalent} if $\phi$ can be chosen
    to satisfy $\phi\theta=\theta'\phi$ on $G^F$.
\end{definition}
Such an equivariant isomorphism $\phi$ identifies $(G^\theta)^F$ with
$(G'^{\theta'})^{F'}$. Pullback along $\phi$ commutes with twisted
duality, so Theorem~\ref{intthm:lp-conjecture} is invariant under
$\mathsf{FI}$-rational equivalence. The equivalences in
Lemma~\ref{lem:f-cyclic-rationally-equivalent} also identify $((G^\theta)^\circ)^F$ with $((G'^{\theta'})^\circ)^{F'}$.

\begin{definition}[cf. {\cite[1.23]{taylor2019rootdata}}]
    \label{def:f-cyclic}
    Fix a direct-product decomposition
    $G=G_1\times\cdots\times G_n$ into connected reductive subgroups. The decomposition is \emph{$F$-cyclic} if
    $F(G_i)=G_{i+1}$, with indices taken modulo $n$. If $G$ admits such a decomposition we call $(G,F)$ \emph{cyclic}.

    A triple $(G,F,\theta)\in\mathsf{FI}$ is
    \emph{$\mathsf{FI}$-cyclic} if $(G,F)$ is cyclic and $\theta$ permutes the factors $G_i$. We say that $(G,F,\theta)$ \emph{arises from a group of simple connected-center type} if $G_1$ can be chosen to be of simple connected-center type.
\end{definition}
\begin{lemma}
    \label{lem:cyclic-projection-preserves-families}
    Let $(G=G_1\times\cdots\times G_n,F)$ be cyclic with connected
    center. Projection onto $G_1$ gives a rational equivalence
    \[
        p_1:G^F\xrightarrow{\sim}G_1^{F^n}.
    \]
\end{lemma}
\begin{proof}
    By \cite[\S8.2, Proposition~8.3, and Corollary~8.8]{taylor2019rootdata},
    $p_1$ is an isomorphism and $p_1^*$ is compatible with Deligne--Lusztig induction and Lusztig series.
    The cyclic-factor reduction in \cite[\S6.1]{digneMichel1990}
    identifies Weyl-group families and preserves the scalar-product
    formulas relating uniform almost characters to irreducible characters.
    Applied to Definition~\ref{def:uniform-almost-character}, this
    reduction identifies the uniform almost characters under $p_1^*$,
    up to nonzero scalars, with their Weyl-group family labels.
    Since $p_1^*$ is an isometry, it preserves and reflects Lusztig
    families by
    Definition~\ref{def:lusztig-families-of-irreducible-characters}.
\end{proof}

$F$-cyclic triples have the following classification, up to rational equivalence.
\begin{lemma}
    \label{lem:f-cyclic-rationally-equivalent}
    Let $(G,F,\theta)\in\mathsf{FI}$ be $\mathsf{FI}$-cyclic with connected
    center.
    Then $(G,F,\theta)$ is $\mathsf{FI}$-rationally equivalent to a triple $(G',F',\theta')\in\mathsf{FI}$ of one of the following forms:
    \begin{enumerate}
        \item $(G_1,F^n|_{G_1},\theta|_{G_1})$ if $\theta$ preserves $G_1$;
        \item $(G_1\times G_1,(F_0\times F_0)\circ\opr{swap},\opr{swap})$
              otherwise. In this case we have $n=2m$, $\theta(G_i)=G_{i+m}$,
              $F_0\coloneqq (\theta\circ F^m)|_{G_1}$, and
              $\opr{swap}(x,y)\coloneqq (y,x)$.
    \end{enumerate}
    In both cases, the equivalence can be chosen to identify
    $((G^\theta)^\circ)^F$ with $((G'^{\theta'})^\circ)^{F'}$.
\end{lemma}
Case (2) is an example of a triple of Galois type (Definition~\ref{def:galois-type-triple}).
\begin{proof}
    The permutation induced by $\theta$ commutes with the $n$-cycle
    induced by $F$, so it is $i\mapsto i+m$ with $2m=0$ modulo $n$.
    If $m=0$, the rational equivalence $p_1$ of
    Lemma~\ref{lem:cyclic-projection-preserves-families} intertwines
    $\theta$ with $\theta|_{G_1}$, giving (1).
    Here $G^\theta=\prod_iG_i^{\theta|_{G_i}}$.
    The inverse $p_1^{-1}:x\mapsto(x,F(x),\ldots,F^{n-1}(x))$ therefore
    identifies $((G_1^{\theta|_{G_1}})^\circ)^{F^n}$ with
    $((G^\theta)^\circ)^F$.

    Otherwise $n=2m$. On $G^F$ we have
    $p_1\theta=F_0p_1$, and $F_0^2=F^n|_{G_1}$.
    Projection onto the first coordinate gives a rational equivalence
    \[
        q_1:(G_1\times G_1)^{(F_0\times F_0)\circ\opr{swap}}
        \xrightarrow{\sim}G_1^{F^n},\qquad (x,F_0(x))\longmapsto x,
    \]
    by Lemma~\ref{lem:cyclic-projection-preserves-families}. Since $q_1\circ\opr{swap}=F_0q_1$,
    the isomorphism $q_1^{-1}p_1$ intertwines $\theta$ with
    $\opr{swap}$, giving (2).
    The subgroup $G^\theta$ is the product of the graphs of
    $\theta|_{G_i}:G_i\to G_{i+m}$ for $1\le i\le m$, so it is
    connected. The target fixed subgroup is the connected diagonal
    $\Delta G_1$. Thus $q_1^{-1}p_1$ identifies the rational points
    of the fixed identity components.
\end{proof}

\subsubsection{Regular embeddings and smooth coverings}
\label{subsubsec:regular-embeddings-and-smooth-coverings}
We recall how regular embeddings and smooth coverings act on Lusztig
families and deduce transfer laws for the main theorem.
\begin{definition}[cf. {\cite[1.13]{taylor2019rootdata}}]
    \label{def:regular-embedding}
    A closed embedding $\iota:G\hookrightarrow G_r$ of
    connected reductive groups is a \emph{regular embedding} if
    $Z(G_r)$ is smooth and connected and
    $\iota(G_{\opr{der}})=(G_r)_{\opr{der}}$.
    A morphism in $\mathsf{FI}$ is a regular embedding if
    its underlying morphism of algebraic groups is.
\end{definition}
\begin{lemma}[{\cite[\S14.7]{taylor2016generalized}}]
    \label{lem:families-under-regular-embeddings}
    Let $\iota:(G,F)\hookrightarrow(G',F')$ be a regular embedding intertwining $F$ and $F'$. For any Lusztig family
    $\cF_{G'}\subseteq\Irr(G'^{F'})$
    all irreducible constituents of the restrictions $\iota^*\chi'$
    for $\chi'\in\cF_{G'}$,
    belong to a single Lusztig family of $G^F$.
\end{lemma}

\begin{corollary}
    \label{cor:transfer-lp-along-regular-embedding}
    Let $\iota:(G,F,\theta)\hookrightarrow
        (G_r,F_r,\theta_r)$ be a regular
    embedding in $\mathsf{FI}$. If $(G_r,F_r,\theta_r)$ satisfies
    Theorem~\ref{intthm:lp-conjecture}, then so does $(G,F,\theta)$.
\end{corollary}
\begin{proof}
    Set $H_r \coloneqq G_r^{\theta_r}$ and $H\coloneqq G^\theta = H_r \cap G$.
    Let $\chi\in \Irr(G^F)$ be $H^F$-distinguished.
    By Mackey theory and Frobenius reciprocity,
    there exists an $H_r^{F_r}$-distinguished $\chi_r\in \Irr(G_r^{F_r})$
    such that
    \[
        \bigl\langle\chi,\Res_{G^F}^{G_r^{F_r}}\chi_r\bigr\rangle_{G^F}>0.
    \]
    Restriction commutes with taking duals and twisting by involutions, so
    \[
        \bigl\langle\chi^\vee\circ\theta,
        \Res_{G^F}^{G_r^{F_r}}(\chi_r^\vee\circ\theta_r)
        \bigr\rangle_{G^F}>0.
    \]
    By assumption and Lemma~\ref{lem:families-under-regular-embeddings},
    $\chi$ and $\chi^\vee\circ\theta$ are in the same family.
\end{proof}
\begin{remark}
    \label{rem:fixed-subgroups-under-regular-embedding}
    In the proof of Corollary~\ref{cor:transfer-lp-along-regular-embedding},
    we used the fact that $H=H_r\cap G$. It is generally
    not true that $H^\circ = H_r^\circ \cap G$.
    The reduction to connected center therefore uses $H^F$-distinction.
    From \S\ref{sec:character-distinction-to-family-distinction} onward,
    we use $(H^\circ)^F$-distinction.

    Among the reduction steps used here, regular embeddings transfer family invariance but need not transfer character equality.
\end{remark}

\begin{definition}[{\cite[1.20]{taylor2019rootdata}}]
    \label{def:smooth-covering}
    A surjective homomorphism $\pi:\widetilde G\to G$ of connected
    reductive groups is a \emph{smooth covering} if $\ker\pi$ is a
    central torus, $\widetilde G_{\opr{der}}$ is simply connected,
    and if $Z(G)$ is connected and/or smooth, then $Z(\widetilde G)$ is connected and/or smooth.
    A morphism in $\mathsf{FI}$ is a smooth covering if its underlying
    morphism of algebraic groups is one.
\end{definition}
\begin{lemma}[{\cite[Corollary~9.3]{digneMichel1990}}]
    \label{lem:families-under-smooth-coverings}
    Let $\pi:(\widetilde G,\widetilde F)\twoheadrightarrow(G,F)$
    be a smooth covering intertwining $\widetilde F$ and $F$, with
    $Z(G)$ smooth and connected. For $\chi,\chi'\in\Irr(G^F)$,
    the inflations $\chi\circ\pi,\chi'\circ\pi \in \Irr(\widetilde G^{\widetilde F})$ belong to the
    same Lusztig family if and only if $\chi$ and $\chi'$ do.
\end{lemma}

\begin{corollary}
    \label{cor:transfer-lp-along-smooth-covering}
    Let $\pi:(\widetilde G,\widetilde F,\widetilde\theta)
        \twoheadrightarrow(G,F,\theta)$ be a smooth covering in
    $\mathsf{FI}$. Set $\widetilde H\coloneqq\widetilde G^{\widetilde\theta}$
    and $H\coloneqq G^\theta$.
    \begin{enumerate}
        \item  If $Z(G)$ is smooth and connected and
              $(\widetilde G,\widetilde F,\widetilde\theta)$ satisfies
              Theorem~\ref{intthm:lp-conjecture}, then so does $(G,F,\theta)$.
        \item  If $\widetilde\chi^\vee\circ\widetilde\theta=\widetilde\chi$
              for every $(\widetilde H^\circ)^{\widetilde F}$-distinguished
              $\widetilde\chi\in\Irr(\widetilde G^{\widetilde F})$, then
              $\chi^\vee\circ\theta=\chi$ for every $(H^\circ)^F$-distinguished
              $\chi\in\Irr(G^F)$.
    \end{enumerate}

\end{corollary}
\begin{proof}
    Since $\ker\pi$ is a connected torus, Lang's theorem makes
    $\widetilde G^{\widetilde F}\to G^F$ surjective.
    Inflation is therefore injective on characters and preserves
    irreducibility.
    For $\chi\in\Irr(G^F)$, set $\widetilde\chi\coloneqq\chi\circ\pi$.
    Equivariance gives
    \[
        \widetilde\chi^\vee\circ\widetilde\theta
        =(\chi^\vee\circ\theta)\circ\pi.
    \]
    If $\chi$ is $H^F$-distinguished, then $\widetilde\chi$ is
    $\widetilde H^{\widetilde F}$-distinguished because
    $\pi(\widetilde H^{\widetilde F})\subseteq H^F$.
    Under the family hypotheses, $\widetilde\chi$ and
    $\widetilde\chi^\vee\circ\widetilde\theta$ belong to the same family.
    Lemma~\ref{lem:families-under-smooth-coverings} gives (1).

    The connected image $\pi(\widetilde H^\circ)$ lies in $H^\circ$.
    Thus, if $\chi$ is $(H^\circ)^F$-distinguished, then
    $\widetilde\chi$ is
    $(\widetilde H^\circ)^{\widetilde F}$-distinguished.
    Under the equality hypothesis,
    $\widetilde\chi^\vee\circ\widetilde\theta=\widetilde\chi$,
    and injectivity of inflation gives $\chi^\vee\circ\theta=\chi$.
\end{proof}

We construct regular embeddings and smooth coverings for triples in $\mathsf{FI}$.
The following construction generalizes the regular embedding and smooth covering
\[
    \SL_n
    \hookrightarrow
    \GL_n
    \twoheadrightarrow
    \PGL_n,
\]
with compatible involutions.
\begin{lemma}
    \label{lem:equivariant-central-pushout}
    Let $(G,F,\theta)\in\mathsf{FI}$, with $F$ defining an
    $\mathbb F_q$-structure, and let $A\subseteq Z(G)$ be a
    $\theta$-stable closed subgroup scheme defined over $\mathbb F_q$
    such that $Z(G/A)$ is smooth and connected. There exists $(G_A,F_A,\theta_A)\in \mathsf{FI}$ such that the quotient map factors as
    \[
        (G,F,\theta)\overset{\iota}{\hookrightarrow}
        (G_A,F_A,\theta_A)\overset{\pi}{\twoheadrightarrow}
        (G/A,F_{G/A},\theta_{G/A}),
    \]
    where $\iota$ is a regular embedding and $\ker\pi$ is a central
    torus. Moreover, $\pi$ is a smooth covering if
    $G_{\opr{der}}$ is simply connected.
\end{lemma}
\begin{proof}
    Choose an $F$-stable maximal torus $T\subseteq G$. Since
    $A\subseteq Z(G)\subseteq T$, the closed embedding
    \[
        j:A\hookrightarrow S\coloneqq T\times T,\qquad a\longmapsto(a,\theta(a)),
    \]
    intertwines $F$ with $F_S\coloneqq F|_T\times F|_T$ and $\theta$ with
    $\theta_S\coloneqq\opr{swap}$. Form the central fppf quotient
    \[
        G_A\coloneqq(G\times S)/\Delta A,\qquad
        \Delta A\coloneqq\{(a,j(a)^{-1}):a\in A\}.
    \]
    The group $G_A$ is connected reductive by
    \cite[Corollary~3.3.5]{Con14}. The maps $F\times F_S$ and
    $\theta\times\theta_S$ descend to the quotient Frobenius $F_A$
    and a commuting involution $\theta_A$.

    The maps $\iota(g)=[g,1]$ and $\pi([g,s])=gA$ are equivariant,
    and $\iota$ is a closed embedding since $j$ is. Identifying
    $S$ with its image $s\mapsto[1,s]$, we have
    \[
        \ker\pi=S,\qquad (G_A)_{\opr{der}}=\iota(G_{\opr{der}}).
    \]
    Indeed, $[a,s]=[1,j(a)s]$ for $a\in A$, and $G_A=\iota(G)S$.
    The center formula for central quotients
    (\cite[Corollary~3.3.5]{Con14}) gives an exact sequence
    \[
        1\longrightarrow S\longrightarrow Z(G_A)
        \longrightarrow Z(G/A)\longrightarrow1.
    \]
    The map $Z(G_A)\to Z(G/A)$ is an $S$-torsor.
    Since $Z(G/A)$ is smooth and connected, $Z(G_A)$ is smooth and
    connected. The center and derived-group identities prove that
    $\iota$ is regular and, if $G_{\opr{der}}$ is simply connected,
    that $\pi$ is a smooth covering.
\end{proof}

\begin{lemma}
    \label{lem:equivariant-embedding-and-covering}
    Let $(G,F,\theta)\in\mathsf{FI}$.
    \begin{enumerate}
        \item There exists a regular embedding
              $(G,F,\theta)\hookrightarrow
                  (G_r,F_r,\theta_r)$ in
              $\mathsf{FI}$.
        \item If $G$ has connected center, there is a smooth
              covering $(\widetilde G,\widetilde F,\widetilde\theta)
                  \twoheadrightarrow(G,F,\theta)$ in $\mathsf{FI}$,
              with $Z(\widetilde G)$ a torus.
    \end{enumerate}
\end{lemma}
\begin{proof}
    For (1), apply Lemma~\ref{lem:equivariant-central-pushout} to
    $A=Z(G)$. Since $G/A$ is adjoint, the center sequence in that
    proof gives $Z(G_A)\cong S$.

    For (2), let $S\coloneqq Z(G)$ and let
    $p:\widetilde{G_{\opr{der}}}\twoheadrightarrow G_{\opr{der}}$
    be the simply connected covering. $F|_{G_{\opr{der}}}$ and
    $\theta|_{G_{\opr{der}}}$ lift uniquely to this covering
    (\cite[Exercise~6.5.2(i),(iii)]{Con14}); uniqueness makes the lifts
    commute and the latter an involution. Together with the
    restrictions to $S$, they make multiplication
    \[
        m:G_0=\widetilde{G_{\opr{der}}}\times S\twoheadrightarrow G,
        \qquad(d,s)\longmapsto p(d)s,
    \]
    a central isogeny in $\mathsf{FI}$. Its kernel is defined over
    the field determining $F$ and stable under the lifted involution.
    Apply Lemma~\ref{lem:equivariant-central-pushout} to $G_0$ with
    these maps and $A\coloneqq\ker m\subseteq Z(G_0)$.
    Since $G_0/A\cong G$ and $(G_0)_{\opr{der}}$ is simply connected,
    the lemma gives a smooth covering of $(G,F,\theta)$.
\end{proof}

\begin{definition}
    \label{def:group-type-triple}
    A triple $(G,F,\theta)\in\mathsf{FI}$ is of \emph{group type}
    if there is a connected reductive group $G'$ with Frobenius $F'$
    such that
    \[
        (G,F,\theta)\cong
        (G'\times G',F'\times F',\opr{swap}).
    \]
    In this case, $(G^\theta)^F$ is identified with the diagonal copy of $G'^{F'}$
    in $G'^{F'}\times G'^{F'}$. We say that $(G,F,\theta)$ \emph{arises from a group of simple connected-center type} if $G'$ can be chosen to be of simple connected-center type.
\end{definition}

\begin{lemma}
    \label{lem:equivariant-factors-separation}
    Let $(G,F,\theta)\in\mathsf{FI}$ be such that $Z\coloneqq Z(G)$ is a torus and $G_{\opr{der}}$ is simply connected and nontrivial. There is a smooth
    covering $\mu:(\widehat G,\widehat F,\widehat\theta)
        \twoheadrightarrow(G,F,\theta)$ in $\mathsf{FI}$ such that
    $Z(\widehat G)$ is a torus, $\mu$ induces an isomorphism on the derived groups, and
    \[
        (\widehat G,\widehat F,\widehat\theta)
        \cong\prod_{i=1}^n(G_i,F_i,\theta_i),
    \]
    where each factor $(G_i,F_i,\theta_i)$ is either
    $\mathsf{FI}$-cyclic arising from a group of simple connected-center
    type, or of group type.
\end{lemma}
\begin{proof}
    Write $G_{\opr{der}}=\prod_{a\in I}D_a$ as
    the product of its simple simply connected factors
    (\cite[Proposition~6.4.4]{Con14}). Set
    \[
        \widehat G_a=ZD_a,\qquad \widehat G=\prod_{a\in I}\widehat G_a.
    \]
    Each $\widehat G_a$ is connected reductive with derived group $D_a$ and center $Z$. Therefore,
    each $\widehat G_a$ is of simple connected-center type.

    The multiplication map $\mu:\widehat G\to G$ is surjective,
    \[
        \ker\mu=\{(z_a)\in Z^I:\prod_a z_a=1\}
        \cong Z^{|I|-1},
    \]
    is a central torus, $Z(\widehat G)=Z^I$ is smooth and connected, and $\mu$ identifies $\widehat G_{\opr{der}}$ with $G_{\opr{der}}=\prod_a D_a$. Therefore, $\mu$ is a smooth covering.

    The maps $F$ and $\theta$ preserve $Z$ and induce commuting
    permutations $\sigma$ and $\tau$ of $I$, with $\tau^2=1$.
    Their restrictions $\widehat G_a\to\widehat G_{\sigma(a)}$ and
    $\widehat G_a\to\widehat G_{\tau(a)}$ define commuting maps
    $\widehat F$ and $\widehat\theta$ on $\widehat G$, making $\mu$
    equivariant. Here $\widehat\theta^2=1$, and $\widehat F$ is
    Frobenius: it is the restriction to $\widehat G\subseteq G^I$
    of the product Frobenius followed by the permutation $\sigma$,
    which defines a finite permutation twist of the rational structure.
    The kernel and derived-group identities show that $\mu$ is a
    smooth covering in $\mathsf{FI}$.

    For each $\sigma$-orbit $\omega\subseteq I$, put
    $Q_\omega\coloneqq\prod_{a\in\omega}\widehat G_a$ and
    $F_\omega\coloneqq\widehat F|_{Q_\omega}$. The involution $\tau$
    acts on the set of these orbits. A $\tau$-fixed orbit $\omega$
    gives the $\mathsf{FI}$-cyclic factor
    $(Q_\omega,F_\omega,\widehat\theta|_{Q_\omega})$, arising from
    any $\widehat G_a$ with $a\in\omega$.
    For a non-fixed orbit $\omega$, the isomorphism
    \[
        Q_\omega\times Q_\omega\xrightarrow{\sim}
        Q_\omega\times Q_{\tau(\omega)},\qquad (x,y)\longmapsto(x,\hat{\theta}(y)),
    \]
    intertwines $F_\omega\times F_\omega$ with $\widehat F$
    and $\opr{swap}$ with $\widehat\theta$, giving a group-type factor.
    Grouping the factors by these $\tau$-orbits proves the decomposition.
\end{proof}

% !TEX root = ../main.tex
\subsection{Tori and triples of group type}
\label{subsec:tori-and-group-type}
We first treat tori and triples of group type.
\begin{lemma}
    \label{lem:lp-torus}
    Every triple $(G,F,\theta)\in\mathsf{FI}$ with $G$ a torus
    satisfies Theorem~\ref{intthm:lp-conjecture}.
    In fact, every $((G^\theta)^\circ)^F$-distinguished irreducible
    character $\chi$ satisfies $\chi^\vee\circ\theta=\chi$.
\end{lemma}
\begin{proof}
    Let $\chi$ be an $((G^\theta)^\circ)^F$-distinguished irreducible
    character of $G^F$. Then $\chi$ is one-dimensional and trivial
    on $((G^\theta)^\circ)^F$.
    The morphism $G\to G^\theta$, $g\mapsto g\theta(g)$, has connected
    image containing $1$, so its image lies in $(G^\theta)^\circ$.
    Thus, for every $g\in G^F$, we have
    \[
        1=\chi(g\theta(g))=\chi(g)\chi(\theta(g)) \implies
        \chi=\chi^\vee\circ\theta.
    \]
\end{proof}
\begin{lemma}
    \label{lem:lp-group-case}
    Every triple $(G,F,\theta)\in\mathsf{FI}$ of group type
    satisfies Theorem~\ref{intthm:lp-conjecture}.
    In fact, every $((G^\theta)^\circ)^F$-distinguished irreducible
    character $\chi$ satisfies $\chi^\vee\circ\theta=\chi$.
\end{lemma}
\begin{proof}
    Write $(G,F,\theta)\cong (G'\times G',F'\times F',\opr{swap})$.
    The fixed subgroup $G^\theta\cong G'$ is connected.
    An irreducible character
    $\alpha\boxtimes\beta$ of $G'^{F'}\times G'^{F'}$ is $\Delta G'^{F'}$-distinguished precisely
    when $\beta=\alpha^\vee$. Thus
    \[
        (\alpha\boxtimes\beta)^\vee\circ\opr{swap}
        =\alpha\boxtimes\beta,
    \]
    as needed.
\end{proof}

% !TEX root = ../main.tex
\subsection{The reduction}
\label{subsec:the-reduction}
We now reduce to simple connected-center triples and Galois triples
arising from groups of simple connected-center type.

\begin{proposition}
    \label{prop:connected-center-reduction}
    \begin{enumerate}
        \item Suppose Theorem~\ref{intthm:lp-conjecture} holds for
              all triples of simple connected-center type and all
              Galois-type triples arising from a group of simple
              connected-center type.
              Then Theorem~\ref{intthm:lp-conjecture} holds for all
              triples in $\mathsf{FI}$.
        \item Assume that $\chi^\vee\circ\theta=\chi$ for every
              $((G^\theta)^\circ)^F$-distinguished irreducible character
              of every triple $(G,F,\theta)$
              of simple connected-center type and every
              Galois-type triple arising from a group of simple
              connected-center type.
              Then $\chi^\vee\circ\theta=\chi$ for every
              $((G^\theta)^\circ)^F$-distinguished irreducible character
              of every triple $(G,F,\theta)$ with connected center.
    \end{enumerate}
    In (2), one may restrict all groups to those whose simple derived
    factors have Dynkin types in any fixed set.
\end{proposition}
\begin{proof}
    Lemma~\ref{lem:lp-torus} settles $(1)$ and $(2)$ when $G$ is a torus. Thus we assume $G_{\opr{der}}\ne1$.

    \emph{(1).} Let $(G,F,\theta)\in\mathsf{FI}$.
    Lemma~\ref{lem:equivariant-embedding-and-covering}(1)
    and Corollary~\ref{cor:transfer-lp-along-regular-embedding}
    reduce to the case where $Z(G)$ is smooth and connected.
    Apply Lemma~\ref{lem:equivariant-embedding-and-covering}(2)
    and Lemma~\ref{lem:equivariant-factors-separation} successively to obtain
    \[
        (\widehat G,\widehat F,\widehat\theta)
        \cong\prod_{i=1}^n(G_i,F_i,\theta_i)
        \overset{\mu}{\twoheadrightarrow}
        (G^{\opr{cov}}, F^{\opr{cov}}, \theta^{\opr{cov}})
        \overset{\pi}{\twoheadrightarrow}
        (G,F,\theta).
    \]
    Corollary~\ref{cor:transfer-lp-along-smooth-covering} and
    Lemma~\ref{lem:transfer-lp-along-products} reduce
    Theorem~\ref{intthm:lp-conjecture} to the factors
    $(G_i,F_i,\theta_i)$.

    By Lemma~\ref{lem:equivariant-factors-separation}, each factor
    is of group type or is $\mathsf{FI}$-cyclic arising from a group
    of simple connected-center type.
    Lemma~\ref{lem:lp-group-case} settles the group-type factors.
    By Lemma~\ref{lem:f-cyclic-rationally-equivalent}, each cyclic
    factor is $\mathsf{FI}$-rationally equivalent to a triple of simple
    connected-center type or a Galois-type triple arising from a group
    of simple connected-center type.
    Both types of triples satisfy
    Theorem~\ref{intthm:lp-conjecture} by the hypothesis of (1).
    Invariance under $\mathsf{FI}$-rational equivalence proves (1).

    \emph{(2).} Let $(G,F,\theta)$ have connected center.
    Apply Lemma~\ref{lem:equivariant-embedding-and-covering}(2)
    and Lemma~\ref{lem:equivariant-factors-separation} directly to
    $(G,F,\theta)$ to obtain the smooth coverings $\pi,\mu$ and
    product decomposition used in (1).
    By the equality clauses of
    Corollary~\ref{cor:transfer-lp-along-smooth-covering} and
    Lemma~\ref{lem:transfer-lp-along-products}, it suffices to prove
    $\chi^\vee\circ\theta_i=\chi$ for every
    $((G_i^{\theta_i})^\circ)^{F_i}$-distinguished
    $\chi\in\Irr(G_i^{F_i})$ and every $i$.
    Lemma~\ref{lem:lp-group-case} settles the group-type factors.
    For each cyclic factor, the equivalence in
    Lemma~\ref{lem:f-cyclic-rationally-equivalent}
    preserves distinction under fixed identity components.
    The hypothesis of (2) gives equality for the resulting simple
    connected-center or Galois-type triple.
    Pullback along the equivalence commutes with twisted duality,
    giving equality for the cyclic factor and proving (2).

    Smooth coverings preserve the Dynkin types of simple derived
    factors, and the product and cyclic reductions only select factors.
    Thus the argument for (2) also applies with the Dynkin types
    restricted to any fixed set.
\end{proof}

% !TEX root = ../main.tex
\section{Character distinction to family distinction}
\label{sec:character-distinction-to-family-distinction}
Assume that $Z(G)$ is connected, and retain the notation of
\S\ref{subsec:lusztig-theory}.
Proposition~\ref{prop:distinguished-uniform-almost-characters}
associates to each $(H^\circ)^F$-distinguished irreducible character
$\chi$ a pair $(\cL,E)$ such that
$\chi\in\operatorname{Fam}^{G^F}_\cL(E)$ and
$\rho^\cL_{\widetilde E}$ is $(H^\circ)^F$-distinguished.

In \S\ref{subsec:hecke-category-and-symmetric-module}, we
construct the monodromic Hecke category $\cH_\fo$ and its symmetric
module category $\cM_\fo^{H^\circ}$, following
\cite[\S3.13]{lusztig2020endoscopy} and
\cite[Chapter~3]{henderson2003character}.
Here $\fo\coloneqq W\cdot\cL$.
Their weighted Grothendieck groups have the standard bases described in
\S\ref{subsubsec:the-generic-module-and-frobenius}.
Specialization at $u^{1/2}=-1$ gives a
$\overline{\bQ}_\ell[W\ltimes\fo]$-module $M_\fo^{H^\circ}$,
whose $\cL$-summand is the $W_\cL$-module $M_\cL^{H^\circ}$.
Proposition~\ref{prop:occurrence} shows that $E$ occurs in
$M_\cL^{H^\circ}$.

Section~\ref{subsec:cross-action-and-specialized-operator-formula}
recalls the cross action and the specialized leading term
from \cite{dagan-monodromic-hecke}.
Using that input, we construct a stable filtration with monomial
associated graded in \S\ref{subsec:the-associated-graded-module}.
Corollary~\ref{cor:relative-module-graded-constituents} expresses
occurrence of $E$ as $(W_\xi,\eta_\xi)$-distinction for some
$\xi\in\Gamma^{H^\circ}_{\cL}$.
Here $W_\xi$ stabilizes $\xi$ under cross action, and $\eta_\xi$
is its character on the corresponding associated-graded basis line.
Theorem~\ref{thm:family-distinction-criterion} gives the resulting necessary
condition for a family to contain an $(H^\circ)^F$-distinguished
irreducible character.

% !TEX root = ../main.tex
\subsection{Distinguished uniform almost characters}
\label{subsec:distinguished-uniform-almost-characters}
The proof uses Lusztig's positivity argument
\cite[\S2.14 and proof of Theorem~8.2]{lusztig2000g}.

\begin{proposition}[cf. {\cite[\S2.14 and proof of Theorem~8.2]{lusztig2000g}}]
    \label{prop:distinguished-uniform-almost-characters}
    Let $\chi$ be a $(H^\circ)^F$-distinguished
    irreducible character of $G^F$. There exist a Kummer local system
    $\cL$ on $T$ with $F$-stable orbit $W\cdot\cL$ and a
    $\gamma_v$-stable $E\in\Irr(W_\cL)$ such that
    \[
        \chi\in\opr{Fam}^{G^F}_\cL(E),\qquad
        \opr{Av}_{(H^\circ)^F}(\rho^\cL_{\widetilde E})\ne0.
    \]
\end{proposition}

\begin{proof}
    Following
    Definition~\ref{def:lusztig-families-of-irreducible-characters},
    choose a Kummer local system $\cL$ whose $W$-orbit is $F$-stable
    and a $\gamma_v$-stable $E_0\in\Irr(W_\cL)$ such that $\chi\in \opr{Fam}^{G^F}_\cL(E_0)$.
    By \cite[6.17(ii), (6.17.2), and (6.18.2)]{lusztig1984characters},
    as applied in \cite[\S\S2.14--2.15 and proof of Theorem~8.2]{lusztig2000g},
    there is a genuine character $R$ of $G^F$ such that
    \[
        R
        = \sum_{E\in\opr{Fam}^{W_\cL}_{\gamma_v}(E_0)} c_E\rho^\cL_{\widetilde E}
        = \sum_{\chi' \in\opr{Fam}^{G^F}_\cL(E_0)} n_{\chi'} \chi'
        \qquad
        c_E \in \overline{\bQ}_\ell,\,
        n_{\chi'} \in \bZ_{\geq 0},\, n_\chi >0.
    \]
    All parameters in the first sum are $\gamma_v$-stable by
    Definition~\ref{def:weyl-group-families}. Changing the finite-order
    extensions only rescales the coefficients $c_E$ by roots of unity.
    Since $\opr{Av}_{(H^\circ)^F}(\chi)>0$ and every
    $\opr{Av}_{(H^\circ)^F}(\chi')$ is a nonnegative integer,
    we have $\opr{Av}_{(H^\circ)^F}(R)>0$.
    The first expansion therefore gives
    $\opr{Av}_{(H^\circ)^F}(\rho^\cL_{\widetilde E})\ne0$
    for some $E\in\opr{Fam}^{W_\cL}_{\gamma_v}(E_0)$.
    Definition~\ref{def:lusztig-families-of-irreducible-characters}
    then gives $\chi\in\opr{Fam}^{G^F}_\cL(E)$.
\end{proof}

% !TEX root = ../main.tex
\subsection{The Hecke category and the symmetric module}
\label{subsec:hecke-category-and-symmetric-module}
Let $B=TU$ be as in \S\ref{subsec:lusztig-theory}, and fix the
$F$-stable orbit $\fo\coloneqq W\cdot\cL$. Set
\[
    \cB\coloneqq G/B,\qquad \widetilde{\cB}\coloneqq G/U,\qquad
    \pi:\widetilde{\cB}\longrightarrow\cB.
\]

\subsubsection{Monodromic categories and convolution}
\label{subsubsec:monodromic-categories-and-convolution}
As in \cite[pp.~27--28]{henderson2003character}, choose an auxiliary
field $\bF_{q^r}$ over which $T$ splits and all elements of $\fo$
are defined. Enlarge $r$ so that every orbit and every irreducible
equivariant local system used below is fixed by $F^{r*}$.
Choose pure weight-zero structures on these local systems.
Such a common choice exists because the orbit and local-system
parameter sets are finite. We use the mixed categories of
Definition~\ref{def:weil-structure} over this auxiliary field.

\begin{definition}[cf. {\cite[\S2.7]{lusztig2020endoscopy}}]
    \label{def:hecke-category}
    For $\cL_1,\cL_2\in\fo$, the \emph{mixed Hecke category with monodromies $(\cL_1,\cL_2)$} is
    \[
        \cH_{\cL_1,\cL_2}
        \coloneqq\cD^{G,\opr{mix}}_{\cL_1\boxtimes\cL_2^{-1}}(\widetilde{\cB}^2),
    \]
    where the action is
    \[
        (g,t_1,t_2)\cdot(g_1U,g_2U)
        =(gg_1t_1^{-1}U,gg_2t_2^{-1}U).
    \]
    Write $\cH_{\cL'}\coloneqq\cH_{\cL',\cL'}$. The \emph{mixed Hecke category for the orbit $\fo$} is
    \[
        \cH_{\fo}
        \coloneqq\bigoplus_{\cL_1,\cL_2\in\fo}\cH_{\cL_1,\cL_2}.
    \]
\end{definition}

\begin{definition}
    \label{def:relative-module-category}
    For a Kummer local system $\cL'$ on $T$, the
    \emph{mixed symmetric Hecke module category with monodromy $\cL'$} and the
    \emph{mixed symmetric Hecke module category for the orbit $\fo$} are
    \[
        \cM^{H^\circ}_{\cL'}\coloneqq\cD^{H^\circ,\opr{mix}}_{\cL'}(\widetilde{\cB}),\qquad
        \cM^{H^\circ}_{\fo}
        \coloneqq\bigoplus_{\cL'\in\fo}\cM^{H^\circ}_{\cL'},
    \]
    respectively. The action defining $\cM^{H^\circ}_{\cL'}$ is
    $(h,t)\cdot gU=hgt^{-1}U$, with $(h,t)\in H^\circ\times T$.
\end{definition}

Let $p_{ij}:\widetilde{\cB}^3\to\widetilde{\cB}^2$ be the projections.
Factor $p_{13}$ as
\[
    \widetilde{\cB}^3
    \xrightarrow{\ q=\opr{id}\times\pi\times\opr{id}\ }
    \widetilde{\cB}\times\cB\times\widetilde{\cB}
    \xrightarrow{\ \overline p_{13}\ }\widetilde{\cB}^2.
\]
For $K\in\cH_{\cL_1,\cL_2}$ and $K'\in\cH_{\cL_2,\cL_3}$,
the monodromies of $p_{12}^*K\otimes p_{23}^*K'$ cancel on the middle copy of $\widetilde{\cB}$.
The resulting middle $T$-equivariance allows descent along $q$.
Write $q_\flat$ for the inverse of the pullback $q^*$. The projection $\overline p_{13}$ is proper,
with fiber $\cB$. Recall the normalized shift
$K\langle n\rangle\coloneqq K[n](n/2)$ of
Definition~\ref{def:tate-twists-and-normalization}, and put
\[
    d_0\coloneqq\dim\widetilde{\cB}+\dim T.
\]

\begin{definition}
    \label{def:hecke-convolution}
    For $K\in\cH_{\cL_1,\cL_2}$ and $K'\in\cH_{\cL_3,\cL_4}$, define their \emph{convolution} by
    \[
        K\star K'
        \coloneqq
        \begin{cases}
            \overline p_{13!}\,q_\flat
            (p_{12}^*K\otimes p_{23}^*K')\langle-d_0\rangle
            \in\cH_{\cL_1,\cL_4} & \text{if } \cL_2 \cong \cL_3 \\
            0                    & \text{otherwise.}
        \end{cases}
    \]
\end{definition}

For the action of $\cH_\fo$ on $\cM^{H^\circ}_\fo$, let $p_2:\widetilde{\cB}^2\to\widetilde{\cB}$
be the second projection and consider
\[
    \widetilde{\cB}^2
    \xrightarrow{\ q_2=\opr{id}\times\pi\ }
    \widetilde{\cB}\times\cB
    \xrightarrow{\ \overline p_1\ }\widetilde{\cB}.
\]
If $A\in\cM^{H^\circ}_{\cL_2}$ and $K\in\cH_{\cL_1,\cL_2}$, then $K\otimes p_2^*A$ descends along $q_2$. Write $q_{2\flat}$ for the inverse of the pullback $q_2^*$.
The projection $\overline p_1$ is proper.

\begin{definition}
    \label{def:relative-convolution}
    For $K\in\cH_{\cL_1,\cL_2}$ and $A\in\cM^{H^\circ}_{\cL_3}$, define their \emph{convolution} by
    \[
        K\star A
        \coloneqq
        \begin{cases}
            \overline p_{1!}\,q_{2\flat}
            (K\otimes p_2^*A)\langle-d_0\rangle
            \in\cM^{H^\circ}_{\cL_1} & \text{if } \cL_2 \cong \cL_3 \\
            0                        & \text{otherwise.}
        \end{cases}
    \]
\end{definition}

The following lemma is standard.
\begin{lemma}[cf. {\cite[\S3.1]{lusztig2020endoscopy}} and
            {\cite[pp.~26--28]{henderson2003character}}]
    \label{lem:convolution-associativity}
    Convolution commutes with forgetting the mixed structure.
    It preserves distinguished triangles and shifts in each variable.
    There are compatible associativity isomorphisms
    \[
        \begin{gathered}
            (K_{12}\star K_{23})\star K_{34}        \simeq K_{12}\star(K_{23}\star K_{34}), \\
            (K_{12}\star K_{23})\star A_3          \simeq K_{12}\star(K_{23}\star A_3)    \\
            (K_{12}\otimes V)\star A_2  \simeq (K_{12}\star A_2)\otimes V
            \simeq K_{12}\star(A_2\otimes V),                                \\
            (K_{12}\otimes V)\star K_{23} \simeq (K_{12}\star K_{23})\otimes V
            \simeq K_{12}\star(K_{23}\otimes V),\\
            F^*(K\star K')\simeq F^*K\star F^*K',\qquad
            F^*(K\star A)\simeq F^*K\star F^*A,
        \end{gathered}
    \]
    where $K_{i,i+1} \in \cH_{\cL_i,\cL_{i+1}}$ and $A_j\in\cM^{H^\circ}_{\cL_j}$ for $i,j \leq 3$, and $V\in \cD^{\opr{mix}}(\opr{pt})$.
\end{lemma}

\subsubsection{Standard objects}
\label{subsubsec:standard-objects}
For each $w\in W$, choose a representative $\dot w$, with $\dot1=1$.
The associated $G\times T\times T$-orbit in $\widetilde{\cB}^2$ is
\[
    \widetilde{\cO}_w
    \coloneqq\{(g_1U,g_2U)\mid g_1^{-1}g_2\in B\dot wB\}\subseteq \widetilde{\cB}^2,
\]
and $j_w:\widetilde{\cO}_w\hookrightarrow\widetilde{\cB}^2$ denotes its inclusion.
The orbit $\widetilde{\cO}_w$ has dimension $d_0+\ell_W(w)$.
By \cite[Prop.~2.2]{henderson2003character}, for any $\cL'\in\fo$ there is a unique irreducible
equivariant local system $\mathscr E_{w,\cL'}$ on $\widetilde{\cO}_w$ with monodromy $w\cL'\boxtimes(\cL')^{-1}$. Give $\mathscr E_{w,\cL'}$ a pure $\bF_{q^r}$-structure of weight $0$.

\begin{definition}[cf. {\cite[\S2.11]{lusztig2020endoscopy}}]
    \label{def:hecke-standard-objects}
    The \emph{standard object} in $\cH_\fo$ indexed by
    $(w,\cL')$ is
    \[
        \Delta_{w,\cL'}\coloneqq j_{w!}
        \bigl(\mathscr E_{w,\cL'}\langle d_0+\ell_W(w)\rangle\bigr)\in \cH_{w\cL', \cL'}\subseteq \cH_\fo.
    \]
\end{definition}

\begin{definition}
    \label{def:standard-basis-parameter}
    For a Kummer local system $\cL'$ on $T$, let $\Gamma^{H^\circ}_{\cL'}$ be the set of pairs $\xi=(\cO_\xi,\mathscr E_\xi)$,
    where
    \begin{itemize}
        \item $\cO_\xi$ is an $H^\circ$-orbit in $\cB$;
        \item $\mathscr E_\xi$ is a representative of an isomorphism class of irreducible $H^\circ$-equivariant $(T,\cL')$-monodromic local systems on $\widetilde{\cO}_\xi\coloneqq\pi^{-1}(\cO_\xi)$.
    \end{itemize}
    Put
    \[
        \Gamma^{H^\circ}_{\fo}\coloneqq\coprod_{\cL'\in\fo}
        \Gamma^{H^\circ}_{\cL'}.
    \]
    For $\xi\in\Gamma^{H^\circ}_{\cL'}$, write
    $\cL_\xi\coloneqq\cL'$ for its torus monodromy.
\end{definition}
We give each $\mathscr E_\xi$ a pure $\bF_{q^r}$-structure of weight $0$.
For $\xi\in\Gamma^{H^\circ}_{\fo}$, write
$j_\xi:\widetilde{\cO}_\xi\hookrightarrow\widetilde{\cB}$ and
$d_\xi\coloneqq\dim\widetilde{\cO}_\xi=\dim\cO_\xi+\dim T$.
\begin{definition}
    \label{def:relative-standard-objects}
    The \emph{standard object} in $\cM^{H^\circ}_\fo$ is
    \[
        \Delta_\xi\coloneqq j_{\xi!}(\mathscr E_\xi\langle d_\xi\rangle).
    \]
\end{definition}

\subsubsection{Standard bases and the generic Hecke algebra}
\label{subsubsec:the-generic-module-and-frobenius}
We use the weight-polynomial description of
\cite[pp.~27--28]{henderson2003character} for the weighted Grothendieck
groups of Definition~\ref{def:tate-twists-and-normalization}.
The standard classes and intersection-cohomology classes differ by a
triangular change of basis with invertible diagonal.
The next definition follows \cite[\S3.13]{lusztig2020endoscopy}.
\begin{definition}[{\cite[1.4]{lusztig2017conjugacy}}]
    \label{def:monodromic-hecke-algebra}
    The \emph{monodromic Hecke algebra}, denoted by $\mathbf H_{\fo}$, is the unital associative $\cA$-algebra with generators $\mathbf T_w$ for $w\in W$ and $1_{\cL'}$ $(\cL' \in \fo)$ and relations
    \[
        \renewcommand{\arraystretch}{1.6}
        \begin{array}{@{}ll@{}}
            1_{\mathcal L}1_{\mathcal L'}
            = \delta_{\mathcal L,\mathcal L'}1_{\mathcal L},
             & \text{for } \mathcal L,\mathcal L' \in \fo; \\

            \mathbf T_w\mathbf T_{w'} = \mathbf T_{ww'},
             & \text{if } w,w' \in W
            \text{ and } \ell_W(ww')=\ell_W(w)+\ell_W(w'); \\

            \mathbf T_w1_{\mathcal L} = 1_{w\mathcal L}\mathbf T_w,
             & \text{for } w\in W,\ \mathcal L\in\fo;      \\

            \displaystyle
            \mathbf T_s^2 = u\mathbf T_1+(u-1)
            \sum_{\mathcal L;\,s\in W_{\mathcal L}}
            \mathbf T_s1_{\mathcal L},
             & \text{for simple reflections } s\in W;      \\[6pt]

            \displaystyle
            \mathbf T_1 = 1=\sum_{\mathcal L\in\fo}1_{\mathcal L}.
             &
        \end{array}
    \]
    Set $\mathbf{H}_{\cL_1,\cL_2} \coloneqq 1_{\cL_1} \mathbf{H}_{\fo} 1_{\cL_2}$ ($\cL_1,\cL_2 \in \fo$), which is generated by $\{1_{\cL_1} \mathbf T_w 1_{\cL_2}\,|\, w\in W,\, w\cL_2 \cong \cL_1   \}$. Write $\mathbf{H}_{\cL'}\coloneqq \mathbf{H}_{\cL',\cL'}$.
\end{definition}
By \cite[\S3.13]{lusztig2020endoscopy},
$\{\mathbf T_w1_{\cL'};(w,\cL')\in W\times\fo\}$ is an $\cA$-basis of
$\mathbf H_\fo$, and the elements with $w\in W_{\cL'}$ form a
basis of $\mathbf H_{\cL'}$. The map $\chi_u$ induces isomorphisms
\[
    K_0^{\mathrm{wt}}(\cH_{\fo})\cong\mathbf{H}_\fo, \qquad
    K_0^{\mathrm{wt}}(\cH_{\cL_1,\cL_2})\cong\mathbf{H}_{\cL_1,\cL_2},
\]
such that $\chi_u(\Delta_{w,\cL'})=u^{-\ell_W(w)/2}\mathbf T_w1_{\cL'}$.
Note that the weight at
$(U,\dot wU)$ of $\Delta_{w,\cL'}$ is equal to $-d_0 -\ell_W(w)$.
This explains the normalization
$\langle-d_0\rangle$ in Definitions~\ref{def:hecke-convolution} and \ref{def:relative-convolution}.
\begin{definition}[cf. {\cite[\S3]{henderson2003character}}]
    \label{def:relative-monodromic-hecke-module}
    For a Kummer local system $\cL'$ on $T$, let
    $\mathbf M^{H^\circ}_{\cL'}$ be the free $\cA$-module with
    basis $m_\xi$ for $\xi\in\Gamma^{H^\circ}_{\cL'}$.
    The \emph{symmetric monodromic Hecke module} for $\fo$ is
    \[
        \mathbf M^{H^\circ}_{\fo}\coloneqq
        \bigoplus_{\cL'\in\fo}\mathbf M^{H^\circ}_{\cL'}.
    \]
\end{definition}
Similarly, \cite[p.~28]{henderson2003character} gives isomorphisms
\[
    K_0^{\mathrm{wt}}(\cM^{H^\circ}_{\fo})\cong\mathbf{M}^{H^\circ}_\fo, \qquad
    K_0^{\mathrm{wt}}(\cM^{H^\circ}_{\cL'})\cong\mathbf{M}^{H^\circ}_{\cL'},
\]
such that $\chi_u(\Delta_{\xi})=u^{-d_\xi/2}m_\xi$.
By \cite[\S3.13]{lusztig2020endoscopy} and
\cite[Theorem~3.2 and Proposition~3.4]{henderson2003character},
convolution gives a left $\mathbf H_{\fo}$-module
$\mathbf M^{H^\circ}_{\fo}$ with
\[
    \chi_u(K)\chi_u(A)=\chi_u(K\star A), \qquad
    A\in \cM^{H^\circ}_{\fo},\,
    K\in \cH_{\fo}.
\]
For each $\cL'\in\fo$, one has
$1_{\cL'}\mathbf M^{H^\circ}_{\fo}=\mathbf M^{H^\circ}_{\cL'}$.
Lemma~\ref{lem:convolution-associativity} also gives a compatible
action of Frobenius.

Set $\epsilon:\cA\longrightarrow\overline{\bQ}_\ell$ by $\epsilon(u^{1/2})=-1$. We also use $\epsilon$ to denote the map $\epsilon:N \to N\otimes_{\cA,\epsilon}\overline{\bQ}_\ell$ for any $\cA$-module $N$ defined by $x\mapsto x\otimes 1$. Then $\epsilon$ induces an isomorphism
\[
    \overline{\bQ}_\ell[W\ltimes \fo]
    \simeq
    \mathbf H_{\fo}
    \otimes_{\cA,\epsilon}
    \overline{\bQ}_\ell,\qquad
    \overline{\bQ}_\ell[W_{\cL'}]
    \simeq
    \mathbf H_{\cL'}
    \otimes_{\cA,\epsilon}
    \overline{\bQ}_\ell,
\]
where $\overline{\bQ}_\ell[W\ltimes \fo]$ is the action groupoid algebra. Convolution then gives an action of $W$ and $W_{\cL'}$ on
\[
    M^{H^\circ}_{\fo} \coloneqq\overline{\bQ}_\ell\otimes_{\cA,\epsilon}\mathbf M^{H^\circ}_{\fo},\qquad
    M^{H^\circ}_{\cL'} \coloneqq\overline{\bQ}_\ell\otimes_{\cA,\epsilon}\mathbf M^{H^\circ}_{\cL'},
\]
respectively.
We also write $m_\xi$ for the image of the generic standard basis
vector under specialization.

% !TEX root = ../main.tex
\subsection{Occurrence in the symmetric Hecke module}
\label{subsec:occurrence-in-the-symmetric-hecke-module}
Recall $Z_{\cL,F}$, $v$ and $\gamma_v$ from \S\ref{subsec:lusztig-theory},
and fix $E\in\Irr(W_\cL)^{\gamma_v}$.
Put $\cA_\ell\coloneqq\overline{\bQ}_\ell[u^{1/2},u^{-1/2}]$.
We prove the following proposition.
\begin{proposition}
    \label{prop:occurrence}
    If $\rho^\cL_{\widetilde E}$ is $(H^\circ)^F$-distinguished, then
    \[
        \Hom_{W_\cL}(E,M^{H^\circ}_{\cL})\ne0.
    \]
\end{proposition}

The proof relates convolution to $(H^\circ)^F$-averaging through
the character functors of Definition~\ref{def:frobenius-character-functors}.
The Schur transporter $\mathsf S^\cL_{\widetilde E}$ of
Definition~\ref{def:schur-transporter} produces a uniform almost
character and detects occurrence of $E$ in the specialized symmetric Hecke module.

The averaging identity in Lemma~\ref{lem:categorification-of-averaging}
contains a nonzero factor $P_{H^\circ}(u)$ that may vanish at $u^{1/2}=-1$.
We therefore complete at $u^{1/2}+1$, cancel this factor, and only
then specialize.

\subsubsection{Categorification}
\label{subsubsec:categorification}
We repeat the construction of \cite[\S3]{henderson2003character} in the $F$-twisted case.
For $\cL_1,\cL_2\in\fo$, set
\[
    \cM^{H^\circ\times H^\circ}_{\cL_1,\cL_2}
    \coloneqq
    \cD^{H^\circ\times H^\circ,\opr{mix}}_{\cL_1\boxtimes\cL_2^{-1}}
    (\widetilde{\cB}^2).
\]
Consider the correspondence
\[
    \begin{tikzcd}[column sep=large]
        \widetilde{\cB}^2
        &G\times\widetilde{\cB}
        \ar[l,"r_F"']\ar[r,"q"]
        &G\times\cB\ar[r,"p"]&G,
    \end{tikzcd}
\]
where $r_F(x,yU)\coloneqq(yU,xF(y)U)$,  $q\coloneqq\opr{id}\times\pi$, and
$p(x,yB)\coloneqq x$.
These maps are $G\times G$-equivariant for the following actions
\begin{align*}
    (g_1,g_2)\cdot(xU,yU)  & =(g_1xU,g_2yU),           &
    (g_1,g_2)\cdot_F(x,yU) & =(g_2xF(g_1)^{-1},g_1yU),   \\
    (g_1,g_2)\cdot_F(x,yB) & =(g_2xF(g_1)^{-1},g_1yB), &
    (g_1,g_2)\cdot_F x     & =g_2xF(g_1)^{-1}.
\end{align*}
Restricting the actions makes these maps $H^\circ\times H^\circ$-equivariant
and equivariant for the diagonal copy of $G$.
The induced action of $G$ on itself is
$\opr{ad}_F(g)\cdot x\coloneqq gxF(g)^{-1}$.

The map $q$ is a $T$-torsor for the action
$(x,yU)\mapsto(x,yt^{-1}U)$, and $r_F$ is $T$-equivariant for the action $(xU,yU) \mapsto (xt^{-1}U,yF(t)^{-1}U)$.
For a $(T\times T,\cL_1\boxtimes\cL_2^{-1})$-monodromic sheaf $K$,
the pullback $r_F^*K$ is
$(T,\cL_1\otimes F^*(\cL_2^{-1}))$-monodromic.
If $\cL_1\cong F^*\cL_2$, the pullback is $T$-equivariant and
descends along $q$.
Write $q_\flat$ for the inverse of $q^*$, without a shift.

\begin{definition}
    \label{def:frobenius-character-functors}
    Define
    \[
        \begin{aligned}
            \opr{Ch}^G_F:
            \cH_{F^*\cL,\cL}
             & \longrightarrow\cD^{(G,\ad_F),\opr{mix}}(G) \\
            \opr{Ch}^{H^\circ\times H^\circ}_F:
            \cM^{H^\circ\times H^\circ}_{F^*\cL,\cL}
             & \longrightarrow
            \cD^{(H^\circ\times H^\circ,\cdot_F),\opr{mix}}(G)
        \end{aligned}
    \]
    by the same formula
    \[
        \opr{Ch}_F(K)\coloneqq p_!q_\flat r_F^*K\langle-d_0\rangle.
    \]
\end{definition}

As in \cite[\S3]{henderson2003character}, set
\[
    \widetilde{\cB}^2
    \xleftarrow{\ \alpha\ }
    H^\circ\times\widetilde{\cB}^2
    \xrightarrow{\ \beta\ }
    \widetilde{\cB}^2,\qquad
    \alpha(h,xU,yU)\coloneqq(hxU,yU),\quad
    \beta(h,xU,yU)\coloneqq(xU,yU).
\]
By \cite[p. 30]{henderson2003character},
\[
    K_0^{\mathrm{wt}}(\cM^{H^\circ\times H^\circ}_{\cL_1,\cL_2})
    \simeq \mathbf{M}^{H^\circ}_{\cL_1}
    \otimes_{\cA} \mathbf{M}^{H^\circ}_{\cL_2^{-1}}.
\]
$\beta_!\alpha^*$ induces a morphism $[\beta_!\alpha^*]:\mathbf{H}_{\cL_1,\cL_2} \to \mathbf{M}^{H^\circ}_{\cL_1}
    \otimes_{\cA} \mathbf{M}^{H^\circ}_{\cL_2^{-1}}$.
The proof of \cite[Proposition~3.8]{henderson2003character} gives
\begin{equation}
    \label{eq:module-action-as-a-functor}
    [\beta_!\alpha^*](\chi_u(K))
    =
    (\chi_u(K)\otimes 1)[\beta_!\alpha^*](1_{\cL_2})
    \in \mathbf{M}^{H^\circ}_{\cL_1}
    \otimes_{\cA} \mathbf{M}^{H^\circ}_{\cL_2^{-1}}, \qquad K\in \cH_{\cL_1,\cL_2}.
\end{equation}
Set
\[
    G
    \xleftarrow{\ \alpha_F\ }
    H^\circ\times G
    \xrightarrow{\ \beta_F\ }
    G,\qquad
    \alpha_F(h,x)\coloneqq xF(h)^{-1},
    \qquad \beta_F(h,x)\coloneqq x,
\]
and consider the functor
\[
    \beta_{F!}\alpha_F^*:
    \cD^{(G,\ad_F),\opr{mix}}(G)\to\cD^{(H^\circ\times H^\circ,\cdot_F),\opr{mix}}(G).
\]
\begin{lemma}[cf. {\cite[Proposition~3.10]{henderson2003character}}]
    \label{lem:ch-alpha-beta-commute}
    The following diagram commutes up to a natural isomorphism:
    \begin{equation}
        \label{eq:ch-alpha-beta-commute}
        \begin{tikzcd}[column sep = huge]
            \cH_{F^*\cL,\cL}
            \ar[r,"\opr{Ch}^G_F"]
            \ar[d,"\beta_!\alpha^*"]
            &
            \cD^{(G,\ad_F),\opr{mix}}(G)
            \ar[d,"\beta_{F!}\alpha_F^*"]
            \\
            \cM^{H^\circ\times H^\circ}_{F^*\cL,\cL}
            \ar[r,"\opr{Ch}^{H^\circ\times H^\circ}_F"]
            &
            \cD^{(H^\circ\times H^\circ,\cdot_F),\opr{mix}}(G)
        \end{tikzcd}
    \end{equation}
\end{lemma}
\begin{proof}
    Consider the diagram
    \[
        \begin{tikzcd}[column sep = 3.5cm]
            \widetilde{\cB}^2
            & H^\circ \times \widetilde{\cB}^2
            \ar[r,"\alpha"]
            \ar[l,"\beta",swap]
            & \widetilde{\cB}^2
            \\
            G \times \widetilde{\cB}
            \ar[u,"r_F"]
            \ar[d,"q"]
            & H^\circ \times G  \times \widetilde{\cB}
            \ar[u,"\id\times r_F"]
            \ar[d,"\id\times q",swap]
            \ar[r,"{(h,x,yU) \mapsto(\alpha_F(h,x),hyU)}"]
            \ar[l, "\beta_F\times \id",swap]
            & G \times \widetilde{\cB}
            \ar[u,"r_F"]
            \ar[d,"q"]
            \\
            G\times \cB
            \ar[d,"p"]
            & H^\circ \times G\times \cB
            \ar[d,"\id\times p",swap]
            \ar[r,"{(h,x,yB) \mapsto(\alpha_F(h,x),hyB)}"]
            \ar[l, "\beta_F\times \id",swap]
            & G \times \cB
            \ar[d,"p"]
            \\
            G
            & H^\circ \times G
            \ar[r,"\alpha_F"]
            \ar[l,"\beta_F",swap]
            & G
        \end{tikzcd}
    \]
    where the top-left and bottom-right squares are Cartesian.
    The middle squares are compatible with $T$-torsor descent along $q$.
    Equivariant base change and descent give the claim.
\end{proof}
Let $i_1:\opr{pt}\to G$ be the inclusion of $1$, and put
\[
    P_{H^\circ}(u)\coloneqq
    \chi_u\!\left(i_1^*(\beta_{F!}\alpha_F^* \underline{\overline{\bQ}_\ell})\right)\in\cA.
\]
We use the representation-valued $\chi_u(i_1^*A)$ of
Definition~\ref{def:tate-twists-and-normalization} when
applying $\operatorname{Av}_{(H^\circ)^F}$.

\begin{lemma}
    \label{lem:categorification-of-averaging}
    For any $A\in\cD^{(G,\ad_F),\opr{mix}}(G)$,
    \[
        \chi_u\!\left(i_1^*(\beta_{F!}\alpha_F^*A)\right)
        =P_{H^\circ}(u)\,\opr{Av}_{(H^\circ)^F}(\chi_u(i_1^*A)) \in \cA.
    \]
    Moreover, $P_{H^\circ}(u)\ne0$. The identity extends linearly over $\cA_\ell$.
\end{lemma}
\begin{proof}
    Consider the fiber diagram
    \[
        \begin{tikzcd}
            H^\circ
            \arrow[r, "j"']
            \arrow[d]
            \arrow[dr, phantom, "\lrcorner", very near start]
            &
            H^\circ \times G
            \arrow[d, "\beta_F"]
            \arrow[r, "\alpha_F"']
            & G
            \\
            \mathrm{pt}
            \arrow[r,"i_1"]
            &
            G
        \end{tikzcd}
    \]
    where $j(h)=(h,1)$.
    The map $\alpha_F\circ j:h\mapsto F(h)^{-1}$ factors through
    a finite universal homeomorphism $H^\circ\to H^\circ$
    defined over $\bF_q$.
    By base change, we obtain
    \[
        i_1^*(\beta_{F!}\alpha_F^*A)
        \cong R\Gamma_c(H^\circ, (\alpha_F\circ j)^* A)
        \cong R\Gamma_c(H^\circ, A|_{H^\circ}).
    \]

    The Lang cover $h\mapsto hF(h)^{-1}$ of $H^\circ$ has
    deck group $(H^\circ)^F$.
    Descent along this finite cover gives
    \[
        R\Gamma_c(H^\circ, A|_{H^\circ})
        \cong
        \bigl(R\Gamma_c(H^\circ,\underline{\overline{\bQ}_\ell})
        \otimes i_1^*A\bigr)^{(H^\circ)^F}
        \cong
        i_1^*(\beta_{F!}\alpha_F^*\underline{\overline{\bQ}_\ell})
        \otimes
        \bigl(i_1^*A\bigr)^{(H^\circ)^F}.
    \]
    Since $H^\circ$ is connected, right translations act trivially
    on $H_c^i(H^\circ,\underline{\overline{\bQ}_\ell})$.
    Together with exactness of $(H^\circ)^F$-invariants,
    this proves the second isomorphism.
    Applying $\chi_u$ gives the required identity.

    As $H^\circ$ is smooth and connected,
    $H_c^{2\dim H^\circ}(H^\circ)=\underline{\overline{\bQ}_\ell}(-\dim H^\circ)$.
    For $i<2\dim H^\circ$, the weights of $H_c^i(H^\circ)$ are at most $i$.
    Thus the term $u^{\dim H^\circ}$ in $P_{H^\circ}(u)$ has coefficient one and
    cannot cancel. Additivity and compatibility with tensoring
    by complexes on a point give the final assertion.
\end{proof}
\begin{remark}
    \label{rem:connected-equivariance-required}
    Connectedness of $H^\circ$ is used both for the Lang cover and
    for triviality of translations on its compactly supported cohomology.
\end{remark}

For $a\in\mathbf{H}_{F^*\cL,\cL}$,
Lemma~\ref{lem:categorification-of-averaging} and
Lemma~\ref{lem:ch-alpha-beta-commute} give
\begin{equation}
    \label{eq:occurrence-weighted-average}
    [i_1^*\opr{Ch}^{H^\circ\times H^\circ}_F\beta_!\alpha^*](a)
    = [i_1^*\beta_{F!}\alpha_F^*\opr{Ch}^{G}_F](a)
    = P_{H^\circ}(u)\opr{Av}_{(H^\circ)^F}
    \bigl([i_1^*\opr{Ch}^G_F](a)\bigr)
    .
\end{equation}
In particular, when $a\in \mathbf{H}_{F^*\cL,\cL}$ annihilates $\mathbf{M}^{H^\circ}_{\cL}$,
$[\beta_!\alpha^*](a)=0$ and thus
\[
    P_{H^\circ}(u)\opr{Av}_{(H^\circ)^F}
    \bigl([i_1^*\opr{Ch}^G_F](a)\bigr)=0.
\]

\subsubsection{Schur transporters}
\label{subsubsec:schur-transporters}
By Lang's theorem we have
\[
    i_1^*:\cD^{(G,\ad_F)}(G)\xrightarrow{\sim}
    D^b(\Rep_{\overline{\bQ}_\ell}(G^F)).
\]
It also identifies the mixed categories with mixed complexes
of $G^F$-representations.
Recall $Z_{\cL,F}$, $v$, $\Gamma_{\widetilde E}$, $R^G_{T_z}(\psi_z)$, and
$\rho^\cL_{\widetilde E}$ from \S\ref{subsec:lusztig-theory}.
Fix $z\in Z_{\cL,F}$. The inverse image of $\widetilde{\cO}_z$
under $yU\mapsto r_F(1,yU)=(yU,F(y)U)$ is
\[
    \{yU\in G/U\mid y^{-1}F(y)\in B\dot zB\}.
\]
On this locus, the torus coordinate $a\in T$ is determined by
$y^{-1}F(y)\in Ua\dot zU$.
The restriction of $r_F^*\mathscr E_{z,\cL}$ to $x=1$ is the
pullback of $z\cL\cong F^*\cL$ along this coordinate map.
The action $yU\mapsto yt^{-1}U$ induces
\[
    a\longmapsto taF_z(t)^{-1},\qquad
    F_z\coloneqq\opr{Int}_{\dot z}\circ F\big|_T.
\]

The orbit map at $a=1$ is the Lang torsor
$t\mapsto tF_z(t)^{-1}$, with stabilizer $T^{F_z}$.
Under descent along this torsor, the equivariant local system
$F^*\cL$ corresponds to its normalized $F_z$-trace character on
$T^{F_z}$; see \cite[Lemma~3.2.4]{eteve2024free}.
Transporting this character through
$\opr{Int}_{g_z}:T^{F_z}\xrightarrow{\sim}T_z^F$ gives
$\psi_z$ by Definition~\ref{lem:kummer-to-rational-character}.
Definition~\ref{def:hecke-standard-objects} and the standard-basis
normalization give
\[
    \chi_u(j_{z!}\mathscr E_{z,\cL})
    =u^{(d_0+\ell_W(z))/2}\chi_u(\Delta_{z,\cL})
    =u^{d_0/2}\mathbf T_z1_\cL.
\]
The shift $\langle-d_0\rangle$ in
Definition~\ref{def:frobenius-character-functors} contributes the
same factor $u^{d_0/2}$. These factors cancel, so
\[
    [i_1^*\opr{Ch}_F^G](\mathbf T_z1_\cL)
    =\chi_u\!\left(i_1^*p_!q_\flat r_F^*
    j_{z!}\mathscr E_{z,\cL}\right).
\]
By \cite[Lemma~3.2.7(ii)]{eteve2024free}, the unshifted complex on the
right computes Deligne--Lusztig cohomology with character $\psi_z$.
Specialization takes its alternating character and therefore gives
\[
    \epsilon\!\left(
    [i_1^*\opr{Ch}^G_F](\mathbf T_z1_\cL)\right)
    =R^G_{T_z}(\psi_z) \in \opr{CF}(G^F), \qquad z\in Z_{\cL,F}.
\]
This identity motivates the following definition.
\begin{definition}
    \label{def:schur-transporter}
    Define the \emph{Schur transporter}
    \[
        \mathsf S^\cL_{\widetilde E}\coloneqq
        \frac{1}{|W_\cL|}\sum_{w\in W_\cL}
        \opr{Tr}(\Gamma_{\widetilde E}E(w))\,\mathbf T_{vw}1_\cL
        \in \mathbf{H}_{F^*\cL,\cL}\otimes_\cA \cA_\ell,
    \]
    and the \emph{central primitive idempotent associated with $E$} by
    \[
        e_E\coloneqq
        \frac{\dim E}{|W_\cL|}
        \sum_{w\in W_\cL}\opr{Tr}(E(w^{-1}))\,w
        \in\overline{\bQ}_\ell[W_\cL].
    \]
\end{definition}
We have
\begin{equation}
    \label{eq:categorification-of-almost-character}
    \epsilon\!\left(
    [i_1^*\opr{Ch}^G_F](\mathsf S^\cL_{\widetilde E})\right)
    =\rho^\cL_{\widetilde E}.
\end{equation}
In addition, $e_E$ acts on every $W_\cL$-module as the projection
onto its $E$-isotypic component.
Since $W_\cL$ is a Weyl group, $E$ is self-dual.
Character orthogonality therefore gives
\begin{equation}
    \label{eq:occurrence-specialized-character-support}
    \epsilon(\mathsf S^\cL_{\widetilde E})e_E
    =\epsilon(\mathsf S^\cL_{\widetilde E}).
\end{equation}

Put
\[
    \widehat{\cA}_\ell\coloneqq\overline{\bQ}_\ell[[\varpi]],
    \qquad
    \cA\longrightarrow\widehat{\cA}_\ell,
    \quad u^{1/2}\longmapsto-1+\varpi.
\]
For a finite free $\cA$-module $N$, put
$\widehat N\coloneqq N\otimes_\cA\widehat\cA_\ell$.
This is the completion of $N\otimes_\cA\cA_\ell$ at $u^{1/2}+1$.
The following diagram commutes.
\begin{equation}
    \label{eq:reduction-specialization}
    \begin{tikzcd}
        N
        \ar[rr,"n\mapsto n\otimes 1"]
        \ar[dr,swap,"\epsilon"]
        &&
        \widehat{N}
        \ar[dl,"\mod\varpi"]
        \\
        &N \otimes_{\cA,\epsilon}  \overline{\bQ}_\ell
    \end{tikzcd}
\end{equation}
All induced functor maps are extended
$\widehat{\cA}_\ell$-linearly.

\begin{lemma}
    \label{lem:non-occurrence-to-vanishing-transporter}
    Assume $\Hom_{W_\cL}(E,M^{H^\circ}_{\cL})=0$.
    There exists
    \[
        \widehat{\mathsf S}_E\in
        1_{F^*\cL}\widehat{\mathbf{H}}_{\fo}1_\cL
    \]
    such that
    \[
        [\beta_!\alpha^*](\widehat{\mathsf S}_E)=0, \qquad
        \widehat{\mathsf S}_E\bmod\varpi=\epsilon(\mathsf S^\cL_{\widetilde E}).
    \]
\end{lemma}
\begin{proof}
    Lift $e_E$ to an idempotent
    $\widehat e_E\in1_\cL\widehat{\mathbf H}_{\fo}1_\cL$.
    To construct the lift, choose any lift $b_0$ and set
    $b_{n+1}\coloneqq3b_n^2-2b_n^3$.
    Since $b_0^2-b_0\equiv0\pmod\varpi$, one has
    \[
        \begin{gathered}
            b_{n+1}^2 -b_{n+1} = (b_n^2-b_n)^2(4b_n^2 -4b_n -3) \equiv 0 \mod\varpi^{2^{n}}, \\
            b_{n+1}-b_n = (b_n^2-b_n)(-2b_n +1) \equiv 0 \mod\varpi^{2^{n}}.
        \end{gathered}
    \]
    Since $1_\cL\widehat{\mathbf H}_{\fo}1_\cL$ is complete, $\lim_{n\to \infty} b_n$ is idempotent.

    Idempotency gives the direct-sum decomposition
    \[
        \widehat{\mathbf M}^{H^\circ}_\cL
        =\widehat e_E\widehat {\mathbf{M}}^{H^\circ}_\cL
        \oplus(1-\widehat e_E)\widehat {\mathbf{M}}^{H^\circ}_\cL.
    \]
    By assumption, the reduction modulo $\varpi$ of the first summand is
    $e_EM^{H^\circ}_{\cL}=0$.
    The module $\widehat{\mathbf M}^{H^\circ}_\cL$ is finite free over $\widehat{\cA}_\ell$.
    Nakayama's lemma therefore gives
    $\widehat e_E\widehat{\mathbf M}^{H^\circ}_\cL=0$.
    Define
    \[
        \widehat{\mathsf S}_E\coloneqq
        \mathsf S^\cL_{\widetilde E}\widehat e_E.
    \]
    Then $\widehat{\mathsf S}_E$ annihilates $\widehat{\mathbf M}^{H^\circ}_\cL$. The claim now follows from
    \eqref{eq:module-action-as-a-functor}, \eqref{eq:occurrence-specialized-character-support} and \eqref{eq:reduction-specialization} extended to $\widehat{\cA}_\ell$.
\end{proof}

\begin{proof}[Proof of Proposition~\ref{prop:occurrence}]
    Suppose that $E$ does not occur in $M^{H^\circ}_{\cL}$.
    Lemma~\ref{lem:non-occurrence-to-vanishing-transporter}
    supplies a completed class $\widehat{\mathsf S}_E$ with
    $[\beta_!\alpha^*](\widehat{\mathsf S}_E)=0$.
    By \eqref{eq:occurrence-weighted-average} extended to $\widehat{\cA}_\ell$, we have
    \[
        P_{H^\circ}((-1+\varpi)^2)\opr{Av}_{(H^\circ)^F}
        \bigl([i_1^*\opr{Ch}^G_F](\widehat{\mathsf S}_E)\bigr)=0
        \quad\text{in }\widehat{\cA}_\ell.
    \]
    By Lemma~\ref{lem:categorification-of-averaging}, $P_{H^\circ}((-1+\varpi)^2)\neq 0$.
    Since $\widehat{\cA}_\ell$ is an integral domain, it follows that
    \[
        \opr{Av}_{(H^\circ)^F}
        \bigl([i_1^*\opr{Ch}^G_F](\widehat{\mathsf S}_E)\bigr)=0.
    \]

    The reduction of $\widehat{\mathsf S}_E$ modulo $\varpi$ is
    $\epsilon(\mathsf S^\cL_{\widetilde E})$.
    Reducing the vanishing identity for the average of
    $[i_1^*\opr{Ch}^G_F](\widehat{\mathsf S}_E)$ modulo $\varpi$
    and applying \eqref{eq:categorification-of-almost-character} gives
    \[
        \opr{Av}_{(H^\circ)^F}(\rho^\cL_{\widetilde E})=0,
    \]
    contradicting $(H^\circ)^F$-distinction.
\end{proof}

\begin{remark}
    \label{rem:occurrence-completion-needed}
    Using $\mathsf S^\cL_{\widetilde{E}}$ instead of $\widehat{\mathsf S}_E$ in the proof of Proposition~\ref{prop:occurrence} gives
    \[
        P_{H^\circ}(u)\opr{Av}_{(H^\circ)^F}
        \bigl([i_1^*\opr{Ch}^G_F](\mathsf S^\cL_{\widetilde{E}})\bigr)\equiv 0 \mod u^{1/2}+1.
    \]
    If $\dim H^\circ>0$, then $P_{H^\circ}(u)\equiv0\pmod{u^{1/2}+1}$, so this congruence gives no information on $\rho^\cL_{\widetilde E}$.
\end{remark}

% !TEX root = ../main.tex
\subsection{Classification of standard objects}
\label{subsec:classification-of-standard-objects}
We describe the standard-basis parameter set $\Gamma^{H^\circ}_{\fo}$
of Definition~\ref{def:standard-basis-parameter} in terms of component-group characters.
By \cite[Proposition~6.8]{helminck1993rationality},
for any $H^\circ$-orbit $\cO$ in $\cB$ there is an element $f\in G$,
unique up to left $H^\circ$-multiplication and right $T$-multiplication,
with $f^{-1}\theta(f)\in N_G(T)$ and $\cO=H^\circ fB/B$.
Let $\mathscr V$ contain one such representative $f$ for each orbit.

Assume that $\fo$ is $F$-stable.
By Remark~\ref{rem:uniform-kummer-exponent}, $[N]^*\cL'$ is
trivial for every $\cL'\in\fo$.
For any $f\in \mathscr V$, set
\begin{gather*}
    \theta_{f} \coloneqq \opr{Int}_{f^{-1}}\circ \theta \circ \opr{Int}_{f},\qquad
    H_f \coloneqq f^{-1}Hf=G^{\theta_f},\qquad
    T_{f,N} \coloneqq [N]^{-1}(T\cap H_f^\circ).
\end{gather*}
Note that $T$ is $\theta_f$-stable if and only if $f^{-1}\theta(f)\in N_G(T)$.
For a finite abelian group $A$, write
$\widehat A\coloneqq\Hom(A,\overline{\bQ}_\ell^\times)$.
We also use $\mu\in\widehat{\pi_0(T_{f,N})}$ to denote its
pullback to $T_{f,N}$ along the component map.
\begin{lemma}
    \label{lem:combi-parameter-rel-standard-obj}
    For $f\in\mathscr V$ and
    $\mu\in\widehat{\pi_0(T_{f,N})}$, put
    \[
        \cL_{f,\mu}\coloneqq
        \bigl([N]_*\overline{\bQ}_\ell\bigr)_{\mu|_{T[N]}}.
    \]
    We have a bijection
    \[
        \Gamma^{H^\circ}_{\fo}
        \longleftrightarrow
        \Big\{
        (f,\mu) \,|\,
        f\in \mathscr{V},\,
        \mu \in \widehat{\pi_0(T_{f,N})},\,
        \cL_{f,\mu} \in \fo
        \Big\}.
    \]
    Under this bijection, $(f,\mu)\in\Gamma^{H^\circ}_{\cL'}$
    if and only if $\cL_{f,\mu}\cong\cL'$.
    In particular, $\cL_\xi\cong\cL_{f,\mu}$ for $\xi=(f,\mu)$.
\end{lemma}
\begin{proof}
    Fix $f\in\mathscr V$ and put
    $\cO\coloneqq H^\circ fB/B$.
    Left multiplication by $f^{-1}$ identifies
    $\pi^{-1}(\cO)$ with $H_f^\circ B/U$ and conjugates $H^\circ$ to $H_f^\circ$.
    We classify irreducible $H_f^\circ\times T$-equivariant local systems
    on $H_f^\circ B/U$ for the action
    \[
        (a,t)\cdot xU\coloneqq ax t^{-N}U.
    \]

    Since $[N]:T\to T$ is surjective, this action is transitive.
    Its stabilizer at $U$ is
    \[
        S_f\coloneqq\opr{Stab}_{H_f^\circ\times T}(U)
        =\{(a,t)\in H_f^\circ\times T\mid a\in t^NU\}.
    \]
    The homomorphism $\psi:T_{f,N}\to S_f$ defined by
    $\psi(t)\coloneqq(t^N,t)$ induces an isomorphism on connected
    components, as we now show.

    Set $V_f\coloneqq U\cap\theta_f(U)$.
    Since $T$ is $\theta_f$-stable, we have a $\theta_f$-stable decomposition
    \[
        B\cap\theta_f(B)=T\ltimes V_f.
    \]
    Uniqueness of this decomposition
    therefore gives
    $B\cap H_f=T^{\theta_f}\ltimes V_f^{\theta_f}$.

    The group $V_f^{\theta_f}$ is connected. Indeed, choose a cocharacter
    $\nu$ of $T$ whose pairing with every root of $U$ is positive.
    If $\alpha$ is a root of $V_f$, both $\alpha$ and $\theta_f(\alpha)$
    are roots of $U$. Hence $\alpha$ has positive pairing with both
    $\nu$ and $\theta_f(\nu)$. Conjugation by the $\theta_f$-fixed cocharacter
    $\nu+\theta_f(\nu)$ contracts $V_f$ to $1$.
    This contraction preserves $V_f^{\theta_f}$, proving connectedness. In particular,
    $V_f^{\theta_f}\subseteq H_f^\circ$, and hence
    \[
        B\cap H_f^\circ
        =(T\cap H_f^\circ)\ltimes V_f^{\theta_f}.
    \]
    Projection onto the second coordinate therefore gives a split exact sequence
    \[
        1\longrightarrow V_f^{\theta_f}
        \xrightarrow{\ v\mapsto(v,1)\ }S_f
        \xrightarrow{\ (a,t)\mapsto t\ } T_{f,N}
        \longrightarrow1,
    \]
    with section $\psi$.
    Since $V_f^{\theta_f}$ is connected, $\psi$ induces
    $\pi_0(T_{f,N})\cong\pi_0(S_f)$.

    The group $H_f^\circ\times T$ is connected. Thus taking the stalk at $U$
    identifies irreducible equivariant local systems on $H_f^\circ B/U$
    with irreducible representations of $\pi_0(S_f)$.
    Since $T_{f,N}$ is a subgroup of the torus $T$, its component group is
    abelian, so these representations are characters.

    The local system associated with $(f,\mu)$ has deck character
    $\mu|_{T[N]}$, since $\psi(z)=(1,z)$ for $z\in T[N]$.
    By Definition~\ref{def:monodromic}, its torus monodromy is $\cL_{f,\mu}$.
    Definition~\ref{def:standard-basis-parameter} therefore identifies
    $\Gamma^{H^\circ}_{\cL'}$ with the pairs satisfying
    $\cL_{f,\mu}\cong\cL'$. Taking the union over $\cL'\in\fo$
    proves the bijection.
\end{proof}

\begin{lemma}
    \label{lem:anti-invariant-lift}
    Let $f\in\mathscr V$ and $\mu\in\widehat{\pi_0(T_{f,N})}$.
    There exists a character $\tilde{\mu}\in X^*(T)$ such that
    \[
        \tilde{\mu}|_{T_{f,N}} = \mu, \qquad \theta_f \tilde{\mu} = - \tilde{\mu}.
    \]
    The local system $\cL_{f,\mu}$ corresponds to the class
    $\tilde\mu/N+X^*(T)$ in Definition~\ref{def:kummer-local-systems}.
    In particular, $\theta_f^*\cL_{f,\mu}\cong\cL_{f,\mu}^{-1}$.
\end{lemma}
\begin{proof}
    The component group of the closed subgroup $T_{f,N}\subseteq T$
    has order prime to $p$. Using the fixed identification of roots of
    unity, regard $\mu$ as an algebraic character of $T_{f,N}$ trivial
    on $T_{f,N}^\circ$. Restriction of algebraic characters from a torus
    to a closed subgroup is surjective. Hence $\mu$ extends to a character
    $\tilde\mu\in X^*(T)$.

    The torus $(T^{\theta_f})^\circ$ is contained in $H_f^\circ$
    and is preserved by $[N]$, so
    $(T^{\theta_f})^\circ\subseteq T_{f,N}^\circ$.
    The homomorphism $t\mapsto t\theta_f(t)$ has connected image
    contained in $T^{\theta_f}$, hence in $(T^{\theta_f})^\circ$.
    Therefore, for every $t\in T$,
    \[
        \tilde\mu(t)\tilde\mu(\theta_f(t))
        =\tilde\mu(t\theta_f(t))=1.
    \]
    Thus $\theta_f\tilde\mu=-\tilde\mu$.

    By Lemma~\ref{lem:combi-parameter-rel-standard-obj}, the deck character
    of $\cL_{f,\mu}$ is $\mu|_{T[N]}=\tilde\mu|_{T[N]}$.
    Definition~\ref{def:kummer-local-systems} therefore identifies
    $\cL_{f,\mu}$ with $\tilde\mu/N+X^*(T)$.
    The equality $\theta_f\tilde\mu=-\tilde\mu$ gives
    $\theta_f^*\cL_{f,\mu}\cong\cL_{f,\mu}^{-1}$.
\end{proof}

\begin{lemma}
    \label{lem:untwisting-the-center}
    Assume that $Z\coloneqq Z(G)$ is connected.
    Let $\xi=(f,\mu)\in\Gamma^{H^\circ}_\fo$.
    There exists a Kummer local system $\cL_0$ on $T$ such that
    \[
        \cL_0|_Z\cong\overline{\bQ}_\ell,\qquad
        {\theta_{f}}^*\cL_0\cong\cL_0^{-1},\qquad
        W_{\cL_0}=W_{\cL_\xi}.
    \]
\end{lemma}
\begin{proof}
    By Lemma~\ref{lem:anti-invariant-lift}, choose an extension
    $\tilde\mu\in X^*(T)$ of $\mu$ with
    $\theta_f\tilde\mu=-\tilde\mu$.
    Since $2N$ is prime to $p$, the class
    $\tilde\mu|_Z/(2N)+X^*(Z)$ defines a Kummer local system
    $\mathcal K$ on $Z$. The lemma gives
    \[
        \mathcal K^{\otimes2}\cong\cL_\xi|_Z,
        \qquad (\theta_f|_Z)^*\mathcal K\cong\mathcal K^{-1}.
    \]

    Put $A\coloneqq G/G_{\opr{der}}$ and let $\pi_A:G\to A$ be the
    quotient map. By \cite[Example~19.25]{Milne_algebraic_groups},
    the restriction $\pi_A|_Z$ is an isogeny of tori.
    The classification in \cite[(2.2.1)]{lusztig1985character1}
    and divisibility of the group of prime-to-$p$ roots of unity imply
    that pullback along an isogeny of tori is surjective on isomorphism
    classes of Kummer local systems.
    Choose a Kummer local system $\cN$ on $A$ with
    $(\pi_A|_Z)^*\cN\cong\mathcal K$.

    Define
    \[
        d:T\longrightarrow A,\qquad
        d(t)\coloneqq\pi_A(t\theta_f(t)^{-1}).
    \]
    Multiplicativity gives
    \[
        (d|_Z)^*\cN
        \cong\mathcal K\otimes\bigl((\theta_f|_Z)^*\mathcal K\bigr)^{-1}
        \cong\mathcal K^{\otimes2}\cong\cL_\xi|_Z.
    \]
    Thus $\cL_0\coloneqq\cL_\xi\otimes(d^*\cN)^{-1}$ is trivial on $Z$.
    Since $d(\theta_f(t))=d(t)^{-1}$ and
    $\theta_f^*\cL_\xi\cong\cL_\xi^{-1}$, we obtain
    $\theta_f^*\cL_0\cong\cL_0^{-1}$.
    Finally, $\theta_f$ induces an automorphism of $A$, so $d$ factors
    through the $W$-invariant map $\pi_A|_T$.
    Thus $d^*\cN$ is $W$-invariant, which gives
    $W_{\cL_0}=W_{\cL_\xi}$ by Definition~\ref{def:relative-weyl-group}.
\end{proof}

% !TEX root = ../main.tex
\subsection{Cross action and the specialized operator formula}
\label{subsec:cross-action-and-specialized-operator-formula}
Retain the parameters $\xi=(f,\mu)\in\Gamma^{H^\circ}_{\fo}$ of
Lemma~\ref{lem:combi-parameter-rel-standard-obj}.
We recall the cross action and the specialized operator statement
from \cite{dagan-monodromic-hecke}.
These results provide the input for the filtration constructed in
\S\ref{subsec:the-associated-graded-module}.
We say that a root $\alpha\in\Phi(G,T)$ is
\emph{$\theta_f$-real} if $\theta_f(\alpha)=-\alpha$.

\begin{definition}[cf. {\cite[\S7.5(c)--(c$'$)]{lusztig-vogan2014quasisplit}}]
    \label{def:cross-action}
    Let $w\in W$ and $(f,\mu)\in \Gamma^{H^\circ}_\fo$.
    For a positive system $P\subseteq\Phi(G,T)$, set
    \[
        \opr{Inv}_{f,P}(w)\coloneqq
        \{\beta\in P\mid\theta_f(\beta)=-\beta,\ w\beta\in-P\}.
    \]
    Let $\dot{w}\in N_G(T)$ be a lift of $w$.
    Let $w\times f$ be the unique element in $\mathscr{V}\cap H^\circ f\dot w^{-1}B$.
    Conjugation by $\dot w$ identifies $T_{f,N}$ with
    $T_{w\times f,N}$; write $(w\mu)(t)\coloneqq\mu(w^{-1}(t))$
    for the transported character.
    For $t\in T_{w\times f,N}$, set
    \[
        \begin{aligned}
            (w\times_{f,P}\mu)(t)
             & \coloneqq \mu(w^{-1}(t))
            \prod_{\beta\in\opr{Inv}_{w\times f,P}(w^{-1})}\beta(t^N), \\
            w\times_P(f,\mu)
             & \coloneqq (w\times f,w\times_{f,P}\mu).
        \end{aligned}
    \]
    Whenever the positive system is omitted from the notation,
    we mean $P=\Phi^+$, the positive system determined by $B$.
\end{definition}
\begin{lemma}[{\cite[Lemma~4.2.2]{dagan-monodromic-hecke}}]
    \label{lem:cross-action-well-defined}
    For each positive system $P$, Definition~\ref{def:cross-action}
    defines an action $\times_P$ of $W$ on $\Gamma^{H^\circ}_\fo$.
    For $w\in W$ and $\xi=(f,\mu)$, one has
    \[
        \theta_{w\times f}=w\theta_fw^{-1}\quad\text{on }T,
        \qquad
        \cL_{w\times_P\xi}\cong w\cL_\xi.
    \]
    For $P=\Phi^+$, $s\in S$, and $t\in T_{s\times f,N}$,
    \[
        (s\times_f\mu)(t)=
        \begin{cases}
            \mu(s(t))\alpha_s(t^N), & \alpha_s\text{ is }\theta_f\text{-real}, \\
            \mu(s(t)),              & \text{otherwise}.
        \end{cases}
    \]
\end{lemma}

\begin{lemma}
    \label{lem:positive-system-comparison}
    Let $\xi=(f,\mu)\in\Gamma^{H^\circ}_\fo$ and $w\in W$
    satisfy $w\times f=f$.
    Suppose that
    \[
        (1-w)\gamma\in(1-\theta_f)X^*(T)
        \qquad
        \text{for every $\theta_f$-real root $\gamma$}.
    \]
    Then the character $w\times_{f,P}\mu$ on $T_{f,N}$
    is independent of the positive system $P$.
\end{lemma}
\begin{proof}
    Since $w\times f=f$, Lemma~\ref{lem:cross-action-well-defined}
    gives $w\theta_f=\theta_fw$. Thus $w$ permutes the real roots.
    Put
    \[
        \rho_{f,P}\coloneqq
        \frac12\sum_{\substack{\gamma\in P\\
                \theta_f\gamma=-\gamma}}\gamma.
    \]
    Since $w$ preserves the real roots,
    \[
        \sum_{\beta\in\opr{Inv}_{f,P}(w^{-1})}\beta
        =(1-w)\rho_{f,P}.
    \]
    For positive systems $P$ and $Q$, the difference
    $\rho_{f,Q}-\rho_{f,P}$ is an integral sum of real roots.
    The hypothesis therefore gives
    \[
        (1-w)(\rho_{f,Q}-\rho_{f,P})
        \in(1-\theta_f)X^*(T).
    \]
    Every character in $N(1-\theta_f)X^*(T)$ is trivial on
    $T_{f,N}$, since $t^N\in T^{\theta_f}$ for $t\in T_{f,N}$.
    Thus the correction factor in $w\times_{f,P}\mu$ is independent
    of $P$. Since $w\mu$ is also independent of $P$, the assertion follows.
\end{proof}

\begin{definition}
    \label{def:theta-degree}
    For $\xi=(f,\mu)\in\Gamma^{H^\circ}_{\fo}$, put
    \[
        T^{-\theta_f}\coloneqq\{t\in T\mid\theta_f(t)=t^{-1}\},
        \qquad r_\theta(\xi)\coloneqq\dim (T^{-\theta_f})^\circ.
    \]
\end{definition}

The conjugacy identity in Lemma~\ref{lem:cross-action-well-defined}
gives $T^{-\theta_{w\times f}}=w(T^{-\theta_f})$.
Consequently $r_\theta(w\times\xi)=r_\theta(\xi)$ for every $w\in W$.

\begin{lemma}[{\cite[Lemma~6.1.3]{dagan-monodromic-hecke}}]
    \label{lem:specialized-rank-one-leading-term}
    Let $s\in S$ and $\xi=(f,\mu)\in\Gamma^{H^\circ}_{\fo}$.
    With $\epsilon(u^{1/2})=-1$, one has
    \[
        \epsilon(\mathbf T_s)m_\xi-c_{s,\xi}m_{s\times\xi}
        \in\opr{span}_{\overline{\bQ}_\ell}
        \{m_\zeta\mid\zeta\in\Gamma^{H^\circ}_{\fo},\
        r_\theta(\zeta)>r_\theta(\xi)\},
    \]
    where
    \[
        c_{s,\xi}\coloneqq
        \begin{cases}
            -1, & \theta_f(\alpha_s)=-\alpha_s, \\
            1,  & \text{otherwise}.
        \end{cases}
    \]
\end{lemma}

% !TEX root = ../main.tex
\subsection{The associated graded module}
\label{subsec:the-associated-graded-module}
We use Lemma~\ref{lem:specialized-rank-one-leading-term} to construct
a $W_\cL$-stable filtration on $M^{H^\circ}_\cL$.
We prove that its associated graded is monomial and determine
the stabilizers and their characters.
The resulting decomposition gives a criterion for occurrence of
a simple $W_\cL$-module.

\subsubsection{Construction of the associated graded}
\label{subsubsec:construction-of-the-associated-graded}
Retain the full monodromy sum $M^{H^\circ}_{\fo}$,
where $\fo=W\cdot\cL$, and its specialization $\epsilon(u^{1/2})=-1$.
Use the degree $r_\theta$ of Definition~\ref{def:theta-degree}.

\begin{definition}
    \label{def:relative-module-filtration}
    For every integer $r$, set
    \[
        \mathrm{Fil}^{\ge r}M^{H^\circ}_{\fo}
        \coloneqq\opr{span}_{\overline{\bQ}_\ell}
        \{m_\xi\mid\xi\in\Gamma^{H^\circ}_{\fo},\ r_\theta(\xi)\ge r\}.
    \]
\end{definition}

Define
\[
    \mathrm{Fil}^{\ge r}M^{H^\circ}_{\cL}
    \coloneqq 1_\cL \mathrm{Fil}^{\ge r}M^{H^\circ}_{\fo}.
\]
For each parameter $\xi$, let $\overline m_\xi$ denote the class
of $m_\xi$ in the summand of degree $r_\theta(\xi)$.
These classes form a basis of the associated graded.
\begin{lemma}
    \label{lem:relative-module-filtration-stability}
    The filtration $\mathrm{Fil}^\bullet M^{H^\circ}_{\fo}$ is
    $\overline{\bQ}_\ell[W\ltimes\fo]$-stable.
    Hence $\mathrm{Fil}^\bullet M^{H^\circ}_{\cL}$ is $W_\cL$-stable.
    On $\opr{gr}_{\mathrm{Fil}} M^{H^\circ}_{\fo}$, $w\in W$
    takes the line of $\overline m_\xi$ to the line of
    $\overline m_{w\times\xi}$.
\end{lemma}
\begin{proof}
    Lemma~\ref{lem:specialized-rank-one-leading-term} gives
    stability under $\epsilon(\mathbf T_s)$ for every $s\in S$.
    Each idempotent $1_{\cL'}$, for $\cL'\in\fo$, projects onto
    the basis vectors indexed by $\Gamma^{H^\circ}_{\cL'}$,
    so it also preserves the filtration.
    These operators generate $\overline{\bQ}_\ell[W\ltimes\fo]$.

    Since $\fo=W\cdot\cL$, the equivalence between representations
    of $W\ltimes\fo$ and representations of $W_\cL$ is given
    by the exact functor $V\mapsto 1_\cL V$.
    Applying this functor to the filtration gives
    $W_\cL$-stability and identifies
    \[
        \opr{gr}_{\mathrm{Fil}} M^{H^\circ}_{\cL}
        \cong 1_\cL\opr{gr}_{\mathrm{Fil}} M^{H^\circ}_{\fo}.
    \]
    On the full associated graded, each $\epsilon(\mathbf T_s)$ sends
    $\overline m_\xi$ to $c_{s,\xi}\overline m_{s\times\xi}$.
    Composing these operators and applying
    Lemma~\ref{lem:cross-action-well-defined} gives the asserted action
    on basis lines for every $w\in W$.
\end{proof}

\subsubsection{Monomial decomposition}
\label{subsubsec:monomial-decomposition}
\begin{definition}
    \label{def:monomial-module}
    A representation $V$ of a group $K$ is \emph{monomial with respect to a basis
        $\{v_i\}_{i\in I}$} if $K$ permutes the basis lines $\{\overline{\bQ}_\ell v_i\mid i\in I\}$.
    The representation $V$ is \emph{monomial} if it admits such a basis.
\end{definition}

\begin{lemma}
    \label{lem:monomial-decomposition}
    Let $K$ be a finite group and $V$ a monomial $K$-representation.
    Choose a basis line $L$ in each orbit of basis lines.
    Let $K_L$ be its stabilizer and $\eta_L$ its character.
    Then
    \[
        V\simeq\bigoplus_{[L]}\Ind_{K_L}^{K}\eta_L.
    \]
\end{lemma}
\begin{proof}
    For each orbit of basis lines, the action map
    \[
        \overline{\bQ}_\ell[K]
        \otimes_{\overline{\bQ}_\ell[K_L]}L
        \longrightarrow V,\qquad g\otimes v\longmapsto gv
    \]
    identifies $\Ind_{K_L}^K\eta_L$ with the span of that orbit:
    coset representatives translate $L$ to distinct basis lines.
    These orbit spans sum directly to $V$.
\end{proof}

Lemma~\ref{lem:relative-module-filtration-stability} makes
$\opr{gr}_{\mathrm{Fil}} M^{H^\circ}_{\cL}$ monomial with respect to $\{\overline{m}_\xi\}_{\xi\in \Gamma^{H^\circ}_\cL}$.
For every $\xi=(f,\mu)\in\Gamma^{H^\circ}_\fo$, define
\begin{equation}
    \label{eq:parameter-stabilizer-and-character}
    W_\xi
    \coloneqq\{w\in W\mid w\times\xi=\xi\},
    \qquad
    w\overline m_\xi=\eta_\xi(w)\overline m_\xi
    \quad(w\in W_\xi).
\end{equation}
Definition~\ref{def:cross-action} gives
$W_\xi\subseteq W_{\cL_\xi}$.
By Lemma~\ref{lem:relative-module-filtration-stability},
$W_\xi$ is the stabilizer of the basis line
$\overline{\bQ}_\ell\overline m_\xi$, and $\eta_\xi$ is its character.
Semisimplicity gives the same irreducible constituents in
$M^{H^\circ}_\cL$ and its associated graded.
Lemma~\ref{lem:monomial-decomposition} and Frobenius reciprocity
therefore give the following corollary.

\begin{corollary}
    \label{cor:relative-module-graded-constituents}
    A simple $W_\cL$-module $E$ occurs in $M^{H^\circ}_{\cL}$
    if and only if $E$ is $(W_\xi,\eta_\xi)$-distinguished
    for some $\xi\in\Gamma^{H^\circ}_{\cL}$.
\end{corollary}

\begin{lemma}
    \label{lem:cross-action-stabilizer-criterion}
    Let $w\in W$ and $\xi=(f,\mu)\in\Gamma^{H^\circ}_{\fo}$.
    Choose $\tilde\mu\in X^*(T)$ extending $\mu$ as in
    Lemma~\ref{lem:anti-invariant-lift}. Then:
    \begin{enumerate}
        \item $w\times f = f$ if and only if $w$ admits a representative in $N_{H^\circ_f}(T)$;
        \item if $w\times f=f$ and
              \[
                  w\tilde\mu-\tilde\mu
                  +N\sum_{\beta\in\opr{Inv}_f(w^{-1})}\beta
                  \in N(1-\theta_f)X^*(T),
              \]
              then $w\times_f\mu=\mu$, and hence $w\in W_\xi$.
    \end{enumerate}
\end{lemma}
\begin{proof}
    By Definition~\ref{def:cross-action}, $w\times f=f$ means
    $f\dot w^{-1}=hft$ for some $h\in H^\circ$ and $t\in T$.
    This equality is equivalent to
    $t\dot w=f^{-1}h^{-1}f\in N_{H_f^\circ}(T)$,
    proving the first assertion.

    If $w\times f=f$, conjugation by such a representative preserves
    $T_{f,N}$. By Definition~\ref{def:cross-action}, the character
    $(w\times_f\mu)\mu^{-1}$ is the restriction to $T_{f,N}$ of
    \[
        w\tilde\mu-\tilde\mu
        +N\sum_{\beta\in\opr{Inv}_f(w^{-1})}\beta.
    \]
    Every character in $N(1-\theta_f)X^*(T)$ restricts trivially
    to $T_{f,N}$. Indeed, for $t\in T_{f,N}$ and $\lambda\in X^*(T)$,
    \[
        \bigl(N(1-\theta_f)\lambda\bigr)(t)
        =\frac{\lambda(t^N)}{\lambda(\theta_f(t^N))}=1,
    \]
    since $t^N\in T\cap H_f^\circ$. This proves the second assertion.
    The equality $w\times\xi=\xi$ then follows from
    Definition~\ref{def:cross-action}.
\end{proof}

\begin{lemma}
    \label{lem:stabilizer-character-positive-system-independence}
    \label{lem:stabilizer-character-real-inversions}
    For $\xi=(f,\mu)\in\Gamma^{H^\circ}_\fo$, $w\in W_\xi$,
    and any positive system $P$, one has
    \[
        \eta_\xi(w)=(-1)^{|\opr{Inv}_{f,P}(w)|},
    \]
    with $\opr{Inv}_{f,P}(w)$ defined in
    Definition~\ref{def:cross-action}.
\end{lemma}
\begin{proof}
    Since $w\in W_\xi$, we have $w\times f=f$.
    Lemma~\ref{lem:cross-action-well-defined} gives
    $w\theta_f=\theta_fw$, so $w$ permutes the real roots.
    Permuting the pairs $\{\gamma,-\gamma\}$ contributes even sign,
    whereas reversing a pair contributes $-1$.
    Hence the sign of the permutation of all real roots is
    $(-1)^{|\opr{Inv}_{f,P}(w)|}$, which is independent of $P$.
    Thus it is enough to prove the formula for $P=\Phi^+$.

    Choose a reduced expression $w=s_r\cdots s_1$ and put
    $v_0\coloneqq1$, $v_i\coloneqq s_i\cdots s_1$.
    The roots $\beta_i\coloneqq v_{i-1}^{-1}\alpha_{s_i}$
    are distinct and enumerate
    $\{\beta\in\Phi^+\mid w\beta\in-\Phi^+\}$,
    by the root-sequence identity for the reduced expression
    $w^{-1}=s_1\cdots s_r$.
    Lemma~\ref{lem:cross-action-well-defined} gives, on $T$,
    $\theta_{v_{i-1}\times f}=v_{i-1}\theta_fv_{i-1}^{-1}$.
    Hence $(s_i,v_{i-1}\times\xi)$ is real exactly when
    $\theta_f\beta_i=-\beta_i$.

    Lemma~\ref{lem:specialized-rank-one-leading-term} therefore gives
    \[
        w\overline m_\xi
        =\left(\prod_{i=1}^r c_{s_i,v_{i-1}\times\xi}\right)
        \overline m_{w\times\xi}
        =(-1)^{|\opr{Inv}_f(w)|}\overline m_{w\times\xi}.
    \]
    Since $w\in W_\xi$, we have $\overline m_{w\times\xi}=\overline m_\xi$.
    The definition of $\eta_\xi$ proves the formula for $P=\Phi^+$.
\end{proof}

% !TEX root = ../main.tex
\subsection{Family distinction}
\label{subsec:family-distinction}
\begin{theorem}
    \label{thm:family-distinction-criterion}
    Let $\chi\in\Irr(G^F)$ be $(H^\circ)^F$-distinguished.
    There exist
    \begin{enumerate}
        \item an element $f\in\mathscr V$;
        \item a character $\mu\in\widehat{\pi_0(T_{f,N})}$ such that
              the orbit $W\cdot\cL_{f,\mu}$
              is $F$-stable;
        \item a simple module $E\in\Irr(W_{\cL_{f,\mu}})^{\gamma_v}$;
    \end{enumerate}
    such that, for $\xi\coloneqq(f,\mu)$,
    \[
        \chi\in\opr{Fam}^{G^F}_{\cL_\xi}(E), \qquad
        E\text{ is }(W_{\xi},\eta_\xi)\text{-distinguished}.
    \]
\end{theorem}
\begin{proof}
    By Proposition~\ref{prop:distinguished-uniform-almost-characters},
    choose $\cL$ with $F$-stable $W$-orbit and
    $E\in\Irr(W_\cL)^{\gamma_v}$ such that
    \[
        \chi\in\opr{Fam}^{G^F}_\cL(E),\qquad
        \opr{Av}_{(H^\circ)^F}(\rho^\cL_{\widetilde E})\ne0.
    \]
    Proposition~\ref{prop:occurrence} implies that $E$ occurs
    in $M^{H^\circ}_\cL$. By
    Corollary~\ref{cor:relative-module-graded-constituents}, there is
    $\xi\in\Gamma^{H^\circ}_\cL$ such that $E$ is
    $(W_\xi,\eta_\xi)$-distinguished.
    Lemma~\ref{lem:combi-parameter-rel-standard-obj}, applied with
    $\fo\coloneqq W\cdot\cL$, writes
    $\xi=(f,\mu)$ with $f\in\mathscr V$,
    $\mu\in\widehat{\pi_0(T_{f,N})}$, and $\cL_{f,\mu}\cong\cL$.
    These data satisfy all the assertions.
\end{proof}

% !TEX root = ../main.tex
\section{Family distinction to twisting invariance}
\label{sec:family-distinction-to-twisting-invariance}
Throughout this section, assume that $Z(G)$ is connected.
In \S\ref{subsec:character-twisting-to-family-twisting}, we establish a formula for the family of a $\theta$-twist.
Let $\chi$ be an $(H^\circ)^F$-distinguished irreducible character.
Theorem~\ref{thm:family-distinction-criterion}
and Proposition~\ref{prop:induced-weyl-group-involution} reduce Theorem~\ref{intthm:lp-conjecture}
to showing that a simple module $E$ of the relative Weyl group $W_\cL$, indexing the family of $\chi$,
satisfies $E\cong E\circ \delta_f$ for a given involution $\delta_f$ on $W_\cL$.
Theorem~\ref{thm:family-distinction-criterion} supplies a parameter
$\xi=(f,\mu)$ with $\cL\coloneqq\cL_\xi$
such that $E$ is $(W_\xi,\eta_\xi)$-distinguished.

In \S\ref{subsec:invariance-under-root-system-involutions},
we classify the obstructions to $E\cong E\circ\delta_f$ for
involutions induced by root-system involutions.
They arise from unequal constituents on exchanged components or from
degenerate constituents on preserved components of type $D_{2m}$
with outer action.

In \S\ref{subsec:family-distinction-to-invariance-condition},
we prove diagonal distinction and exclude the degenerate
constituents needed in simple connected-center type.
In \S\ref{sec:proofs-of-the-main-theorems}
we then apply these criteria
to the simple connected-center and Galois cases, then use
Proposition~\ref{prop:connected-center-reduction} to conclude
both main theorems.

% !TEX root = ../main.tex
\subsection{Character twisting to family twisting}
\label{subsec:character-twisting-to-family-twisting}

Let $\mathcal L$ be a Kummer local system on $T$ with $F$-stable $W$-orbit.
Choose $f\in G$ such that
\[
    \theta_f(T)=T,\qquad \theta_f^*\mathcal L\cong\mathcal L^{-1}.
\]
Put $\delta_f\coloneqq(\theta_f)_W$.
Then $\delta_f$ preserves $W_{\mathcal L}$.

\begin{proposition}
    \label{prop:induced-weyl-group-involution}
    For every $E\in\opr{Irr}(W_{\mathcal L})^{\gamma_v}$, the representation $E\circ\delta_f\in\opr{Irr}(W_{\mathcal L})^{\gamma_v}$ satisfies
    \[
        \chi\in\opr{Fam}^{G^F}_{\mathcal L}(E)
        \quad\Longrightarrow\quad
        \chi^\vee\circ\theta
        \in\opr{Fam}^{G^F}_{\mathcal L}(E\circ\delta_f).
    \]
\end{proposition}
We first prove a few auxiliary lemmas.
\begin{lemma}
    \label{lem:dl-theta-inv-functoriality}
    Set $h \coloneqq F^{-1}(f)$ and $c_h\coloneqq h^{-1}\theta(h)$. Then $c_hF(c_h)^{-1}\in N_G(T)$. Denote its image in $W$ by $\kappa_{f}$.
    The map
    \[
        \alpha_f:Z_{\cL,F}\longrightarrow Z_{\cL,F},
        \qquad \alpha_f(z)\coloneqq \kappa_{f}\delta_f(z)
    \]
    is an involution, and
    \[
        \bigl(R^G_{T_z}(\psi_z)\bigr)^\vee\circ\theta
        =R^G_{T_{\alpha_f(z)}}(\psi_{\alpha_f(z)}).
    \]
\end{lemma}
\begin{proof}
    For $z\in Z_{\cL,F}$, set $g'_z\coloneqq\theta(g_z)c_h^{-1}$.
    Then $g'_zT(g'_z)^{-1}=\theta(T_z)$ and
    \[
        (g'_z)^{-1}F(g'_z)
        =c_h\theta(\dot z)F(c_h)^{-1}
        =c_hF(c_h)^{-1}\theta_f(\dot z).
    \]
    By definition, the image of this element in $W$ is $\alpha_f(z)$.

    Since $\theta\circ\opr{Int}_{g'_z}=\opr{Int}_{g_z}\circ\theta_h$, one has
    \[
        \opr{Int}_{g'_z}^*\theta^*\cL_z^{-1}
        \cong\theta_h^*F^*(\cL^{-1})
        \cong F^*\cL,
    \]
    where $\cL_z\coloneqq (\opr{Int}_{g_z^{-1}})^*(F^*\cL)$.
    Since $\theta^*\cL_z^{-1}$ is $F$-stable,
    the criterion in Definition~\ref{lem:kummer-to-rational-character} gives
    $\alpha_f(z)\in Z_{\cL,F}$.
    Independence of choices and functoriality of normalized traces identify
    $(\theta(T_z),\psi_z^{-1}\circ\theta)$ with
    $(T_{\alpha_f(z)},\psi_{\alpha_f(z)})$ up to $G^F$-conjugacy.
    Lemma~\ref{lem:dl-functoriality} now proves $\bigl(R^G_{T_z}(\psi_z)\bigr)^\vee\circ\theta
        =R^G_{T_{\alpha_f(z)}}(\psi_{\alpha_f(z)})$.

    Finally, $\theta(c_h)=c_h^{-1}$ implies
    \[
        \bigl(c_hF(c_h)^{-1}\bigr)
        \theta_f\bigl(c_hF(c_h)^{-1}\bigr)=1.
    \]
    Hence $\kappa_{f}\delta_f(\kappa_{f})=1$, and $\alpha_f^2=1$.
\end{proof}
\begin{lemma}
    \label{lem:unif-alm-char-theta-inv-functoriality}
    For $E\in\Irr(W_\cL)^{\gamma_v}$, the representation
    $E\circ\delta_f$ is $\gamma_v$-stable.
    For the chosen finite-order extensions, there is a root of unity
    $\zeta$ such that
    \[
        \bigl(\rho^\cL_{\widetilde E}\bigr)^\vee\circ\theta
        =\zeta\,\rho^\cL_{\widetilde{E\circ\delta_f}}.
    \]
    If $E\cong E\circ\delta_f$, then $\zeta=1$.
\end{lemma}
\begin{proof}
    Let $F_*$ denote the map induced by $F$ on $X_*(T)$, and put
    $\varphi\coloneqq q^{-1}F_*\circ v$ on
    $\Span_\bR(\Phi_\cL^\vee)$.
    The relation $v\cL\cong F^*\cL$ implies that $\varphi$ preserves
    $\Phi_\cL^\vee$.
    Minimality of $v$ and preservation of $(\Phi^+)^\vee$ by
    $q^{-1}F_*$ show that $\varphi$ preserves $(\Phi_\cL^+)^\vee$.
    Conjugation by $\varphi$ induces $\gamma_v$ on $W_\cL$.
    On $\Span_\bR(\Phi_\cL^\vee)$, write $\theta_f=x\tau$, where $x\in W_\cL$
    and $\tau$ preserves $(\Phi_\cL^+)^\vee$, by
    \cite[Proposition~C.53]{Milne_algebraic_groups}.
    We also write $\tau$ for the induced automorphism of $W_\cL$;
    thus $\delta_f=\opr{Int}_x\circ\tau$.

    The identities $\theta_f\circ F=F\circ\theta_h$ and
    $\theta_h=\opr{Int}_{c_hF(c_h)^{-1}}\circ\theta_f$ on $T$ give
    \[
        \theta_f\varphi\theta_f^{-1}
        =q^{-1}F_*\kappa_f\delta_f(v)
        \in\varphi W_\cL.
    \]
    Projecting onto the automorphism group of the Dynkin diagram
    gives $\tau\varphi=\varphi\tau$.
    Hence $\theta_f\varphi\theta_f^{-1}=x\varphi x^{-1}$, and
    \[
        v^{-1}\alpha_f(vw)
        =\gamma_v^{-1}(x)\tau(w)x^{-1},\qquad w\in W_\cL.
    \]

    Since $\Gamma_{\widetilde E}E(\gamma_v^{-1}(x))
        =E(x)\Gamma_{\widetilde E}$, cyclicity of trace gives
    \begin{align*}
        \opr{Tr}(\Gamma_{\widetilde E}E(v^{-1}\alpha_f(vw)))
         & =\opr{Tr}(E(x)\Gamma_{\widetilde E}E(\tau(w))E(x)^{-1}) \\
         & =\opr{Tr}(\Gamma_{\widetilde E}E(\tau(w))).
    \end{align*}
    By Lemma~\ref{lem:dl-theta-inv-functoriality} and $\alpha_f^2=1$,
    reindexing the defining sum yields
    \[
        \bigl(\rho^\cL_{\widetilde E}\bigr)^\vee\circ\theta
        =\frac{1}{|W_\cL|}\sum_{w\in W_\cL}
        \opr{Tr}(\Gamma_{\widetilde E}E(\tau(w)))R^G_{T_{vw}}(\psi_{vw}).
    \]
    Since $\tau$ commutes with $\gamma_v$, the operator
    $\Gamma_{\widetilde E}$ defines a finite-order extension of
    $E\circ\tau\cong E\circ\delta_f$.
    Thus $E\circ\delta_f$ is $\gamma_v$-stable.
    Comparing extensions as after Definition~\ref{def:uniform-almost-character}
    proves the equality up to a root of unity.

    If $E\cong E\circ\delta_f$, then $E$ is stable under the abelian
    group of diagram automorphisms $\langle\varphi,\tau\rangle$.
    By \cite[Lemma~0.2]{digneMichel1990}, $E$ extends to
    $W_\cL\rtimes\langle\varphi,\tau\rangle$.
    In this extension, the operator representing $\varphi$ is a scalar
    multiple of $\Gamma_{\widetilde E}$ by Schur's lemma.
    Conjugating by the operator representing $\tau$ therefore gives
    \[
        \opr{Tr}(\Gamma_{\widetilde E}E(\tau(w)))
        =\opr{Tr}(\Gamma_{\widetilde E}E(w)).
    \]
    The reindexed sum equals
    $\rho^\cL_{\widetilde E}$, proving $\zeta=1$.
\end{proof}

\begin{proof}[Proof of Proposition~\ref{prop:induced-weyl-group-involution}]
    The map $\chi\mapsto\chi^\vee\circ\theta$ permutes $\Irr(G^F)$ and
    preserves the character pairing. For
    $E_1,E_2\in\Irr(W_\cL)^{\gamma_v}$,
    Lemma~\ref{lem:unif-alm-char-theta-inv-functoriality} shows that
    twisting sends a common nonzero constituent of
    $\rho^\cL_{\widetilde E_1}$ and $\rho^\cL_{\widetilde E_2}$
    to a common nonzero constituent
    of $\rho^\cL_{\widetilde{E_1\circ\delta_f}}$ and
    $\rho^\cL_{\widetilde{E_2\circ\delta_f}}$.
    Since $\delta_f^2=1$, the converse holds as well.
    Thus $E\mapsto E\circ\delta_f$ preserves the relation in
    \eqref{eq:family-defining-relation} in both directions.
    Passing to its equivalence closure gives
    \[
        \opr{Fam}^{W_\cL}_{\gamma_v}(E\circ\delta_f)
        =\{E'\circ\delta_f\mid E'\in\opr{Fam}^{W_\cL}_{\gamma_v}(E)\}.
    \]

    If $\chi\in\opr{Fam}^{G^F}_\cL(E)$, choose
    $E'\in\opr{Fam}^{W_\cL}_{\gamma_v}(E)$ with
    $\langle\chi,\rho^\cL_{\widetilde E'}\rangle_{G^F}\ne0$.
    The same twisting identity gives
    \[
        \bigl\langle\chi^\vee\circ\theta,
        \rho^\cL_{\widetilde{E'\circ\delta_f}}
        \bigr\rangle_{G^F}\ne0.
    \]
    Since $E'\circ\delta_f\in\opr{Fam}^{W_\cL}_{\gamma_v}(E\circ\delta_f)$,
    Definition~\ref{def:lusztig-families-of-irreducible-characters}
    gives the conclusion.
\end{proof}

% !TEX root = ../main.tex
\subsection{Invariance under root-system involutions}
\label{subsec:invariance-under-root-system-involutions}
In this subsection, let $\Psi\subseteq X$ be a finite reduced
crystallographic root system spanning a real vector space $X$.
In the application, $\Psi=\Phi_\cL$ and $W(\Psi)=W_\cL$,
where $\cL$ is a Kummer local system on $T$ and $Z(G)$ is connected.

Write
\[
    \Psi=\bigsqcup_{i=1}^n\Psi_i,
    \qquad W(\Psi)=\prod_{i=1}^n W_i,
    \qquad W_i\coloneqq W(\Psi_i),
\]
where the $\Psi_i$ are the irreducible components of $\Psi$.
Let $\delta$ be an involution on $\Psi$.
Denote the induced involution on $W(\Psi)$ by $\delta$ as well.
% The orbits of $\delta$ on the factors $W_i$ are singletons
% and exchanged pairs.
For $E\in\opr{Irr}(W(\Psi))$, write
\[
    E\cong E_1\boxtimes\cdots\boxtimes E_n,
    \qquad E_i\in\opr{Irr}(W_i).
\]

Set
\[
    I_\delta\coloneqq
    \big\{(i,j)\mid 1\leq i<j\leq n,\;
    \delta(W_i)=W_j\big\}.
\]
For $(i,j)\in I_\delta$ define
the \emph{twisted diagonal subgroup}
\[
    \Delta_{\delta,i,j}\coloneqq
    \big\{(w,\delta(w))\mid w\in W_i\big\}
    \subseteq W_i\times W_j.
\]
Set
\[
    \Delta_\delta^{\mathrm{sw}}\coloneqq
    \prod_{(i,j)\in I_\delta}\Delta_{\delta,i,j}\subseteq W(\Psi).
\]

Next we describe the action of $\delta$ on a preserved component.
For $r\geq4$, use the signed-permutation realizations
\[
    C_2\coloneqq\{\pm1\},\qquad
    W(B_r)\coloneqq C_2^r\rtimes S_r,
    \qquad
    W(D_r)\coloneqq\left\{(\mathbf a,\pi)\in W(B_r)\,\Big|\,
    \prod_{i=1}^r a_i=1\right\}.
\]
Here $S_r$ permutes the coordinates of $C_2^r$.

\begin{lemma}
    \label{lem:involution-on-component-is-outer-only-in-d2m}
    Suppose $\delta(\Psi_i)=\Psi_i$.
    Then $\delta|_{W_i}$ is inner unless $\Psi_i$ has type
    $D_{2m}$ for some $m\geq2$.
    In the outer case, there exists a root-system isomorphism
    \(\iota_i:\Psi_i\to\Phi(D_{2m})\) satisfying
    \[
        \iota_i\circ\delta\circ\iota_i^{-1}=\tau,
        \qquad
        \tau\in W(B_{2m})\setminus W(D_{2m}),\qquad \tau^2=1.
    \]
\end{lemma}
\begin{proof}
    Choose a base $S_i$ of $\Psi_i$.
    By \cite[Proposition~C.53]{Milne_algebraic_groups},
    \[
        \operatorname{Aut}(\Psi_i)
        =W_i\rtimes\operatorname{Aut}(\opr{Dyn}_{\Psi_i,S_i}).
    \]
    Write $\delta|_{\operatorname{span}_{\bR}(\Psi_i)}=w\eta$
    with $w\in W_i$ and $\eta$ a diagram involution.
    On $W_i$, we obtain $\delta|_{W_i}=\operatorname{Int}_w\circ\eta_{W_i}$, where $\eta_{W_i}$
    is induced by $\eta$.

    The classification of Dynkin diagrams
    (\cite[Theorem~C.56]{Milne_algebraic_groups}) shows that
    nontrivial diagram automorphisms occur only in types
    $A_r$ ($r\geq2$), $D_r$ ($r\geq4$), and $E_6$.
    Let $w_0$ be the longest element of $W_i$.
    In types $A_r$, $D_{2m+1}$, and $E_6$, the nontrivial diagram
    automorphism is $-w_0$ on the root space.
    Its action on the Weyl group is therefore $\operatorname{Int}_{w_0}$.
    Hence an outer action can occur only in type $D_{2m}$.
    In type $D_{2m}$, a nontrivial $\eta$ is the fork exchange for a suitable choice of coordinates.
    In type $D_4$, choose the identification so that $\eta$ becomes the chosen transposition of the three outer nodes.
    In these coordinates, $\delta$ has the required form.
\end{proof}
For any element $\tau=(\mathbf a,\pi)\in W(B_{2m})\setminus W(D_{2m})$ with $\tau^2=1$, we have
\[
    \pi^2=1, \qquad
    a_i = a_{\pi(i)}, \qquad
    \pi(i_0) = i_0 \text{ for some } i_0\leq 2m.
\]

We recall the irreducible representations of $W(D_r)$ from
\cite[\S\S2.2--2.3]{taylor2015induced}.
Let $(\lambda,\mu)$ be a bipartition of $r$,
and let $S^\lambda$ and $S^\mu$ be the Specht representations of
$S_{|\lambda|}$ and $S_{|\mu|}$ indexed by $\lambda$ and $\mu$, respectively.
Define
\begin{equation}
    \label{eq:type-B-irreps}
    \chi_{\lambda,\mu}: C_2^r \to \{\pm 1\},\qquad
    \chi_{\lambda,\mu}(\mathbf a) \coloneqq \prod_{i=|\lambda|+1}^{r}a_i,\qquad
    V_{\lambda,\mu} \coloneqq \ind_{C_2^r \rtimes S_{|\lambda|}\times S_{|\mu|}}^{W(B_{r})} \chi_{\lambda,\mu} \otimes (S^\lambda\boxtimes S^\mu).
\end{equation}
As $(\lambda,\mu)$ ranges over the bipartitions of $r$, the representations
$V_{\lambda,\mu}$ exhaust $\opr{Irr}(W(B_r))$ without repetition.

For $\lambda\ne\mu$, the restriction
\[
    E_{\{\lambda,\mu\}}\coloneqq
    V_{\lambda,\mu}|_{W(D_r)}
    \cong V_{\mu,\lambda}|_{W(D_r)}
\]
is irreducible and depends only on the unordered pair $\{\lambda,\mu\}$.
For $\lambda=\mu$, the rank $r$ is even and
\[
    V_{\lambda,\lambda}|_{W(D_r)}
    \cong E_\lambda^+\oplus E_\lambda^-,
\]
where the two summands are irreducible and non-isomorphic.
These representations exhaust $\opr{Irr}(W(D_r))$ with no repetition.

\begin{definition}
    \label{def:degenerate-d2m-modules}
    For $m\geq2$, an irreducible $W(D_{2m})$-module is
    \emph{degenerate} if it is isomorphic to $E_\lambda^+$ or
    $E_\lambda^-$ for a partition $\lambda$ of $m$.
\end{definition}

On each $\delta$-stable component $W_i$ on which $\delta$ acts outerly,
fix $\iota_i:W_i\to W(D_{2m})$ as in
Lemma~\ref{lem:involution-on-component-is-outer-only-in-d2m}.
We say that an irreducible module $E_i$ of $W_i$ is degenerate if $E_i\circ\iota_i^{-1}$ is.

\begin{proposition}
    \label{prop:non-fixed-module-criteria}
    Let $E\in\opr{Irr}(W(\Psi))$ and write $E\cong\boxtimes_{i=1}^n E_i$.
    Then $E\not\cong E\circ\delta$ if and only if at least one
    of the following conditions holds:
    \begin{enumerate}
        \item $E$ is not $\Delta_\delta^{\mathrm{sw}}$-distinguished.
        \item For some $\delta$-stable component $W_i$ on which
              $\delta$ acts outerly, $E_i$ is degenerate.
    \end{enumerate}
\end{proposition}
\begin{proof}
    $E\cong E\circ\delta$ if and only if
    $E_i\cong E_i\circ\delta|_{W_i}$ on each preserved component
    and, for every $(i,j)\in I_\delta$,
    \[
        E_i\cong E_j\circ\delta|_{W_i}.
    \]
    Every irreducible representation of a finite Weyl group in
    characteristic zero is self-dual. Hence
    \[
        (E_i\boxtimes E_j)^{\Delta_{\delta,i,j}}
        \cong\operatorname{Hom}_{W_i}
        (E_i^*,E_j\circ\delta|_{W_i})
    \]
    is nonzero precisely when $E_i\cong E_j\circ\delta|_{W_i}$.
    $E_i\cong E_j\circ\delta|_{W_i}$ for all $(i,j)\in I_\delta$
    if and only if $E$ is $\Delta_\delta^{\mathrm{sw}}$-distinguished.

    On a preserved component, inner twisting fixes the isomorphism
    class of $E_i$.
    In the outer case,
    Lemma~\ref{lem:involution-on-component-is-outer-only-in-d2m}
    identifies the action of $\delta$ on $W_i$ with conjugation by
    $\tau\in W(B_{2m})\setminus W(D_{2m})$, where $\tau^2=1$.
    Every non-degenerate representation extends to $W(B_{2m})$,
    so conjugation by $\tau$ preserves its isomorphism class.
    The degenerate representations satisfy
    $E_\lambda^+\circ\opr{Int}_{\tau}\cong E_\lambda^-$.
    Consequently $E_i\not\cong E_i\circ\delta|_{W_i}$ if and only if
    $E_i\circ\iota_i^{-1}$ is degenerate.

    The exchanged pairs give condition~(1), and the preserved
    components give condition~(2).
\end{proof}

\begin{lemma}
    \label{lem:fixed-label-implies-fixed-family}
    Assume that every irreducible component of $\Psi$ exchanged
    with another component by $\delta$ has type $A$.
    Let $\cF_W$ be a Lusztig family of $W(\Psi)$ containing
    $E_0\cong E_0\circ\delta$. Then
    $E\cong E\circ\delta$ for every $E\in\cF_W$.
\end{lemma}
\begin{proof}
    Write $\cF_W=\boxtimes_i\cF_{W_i}$ and
    $E_0=\boxtimes_i E_{0,i}$ using the product decomposition
    of Weyl-group families \cite[\S4.2]{lusztig1984characters}.
    On exchanged components, the families $\cF_{W_i}$ are
    singletons because the components have type $A$.
    Thus the matching identities for $E_{0,i}$ hold for the factors
    of every $E\in\cF_W$.

    On a preserved component with outer action,
    Proposition~\ref{prop:non-fixed-module-criteria} shows that
    $E_{0,i}$ is non-degenerate.
    Each degenerate representation forms a family by itself
    \cite[(4.6.10)]{lusztig1984characters}.
    Hence $\cF_{W_i}$ contains no degenerate representation.
    The same proposition now gives $E\cong E\circ\delta$
    for every $E\in\cF_W$.
\end{proof}

The following distinction condition excludes degenerate $W(D_{2m})$-modules.
\begin{lemma}
    \label{lem:no-degenerate-in-induced-d2m-module}
    Let $\pi \in S_{2m}$ satisfy $\pi^2=1$ and admit a fixed point.
    Set
    \[
        A_{\pi}\coloneqq
        \left\{
        \mathbf b \in C_2^{2m}\mid
        b_i = b_{\pi(i)},\, \prod_{i\leq 2m}b_i=1
        \right\}\subseteq C_2^{2m}\cap W(D_{2m}),\qquad
        \eta_\pi(\mathbf b)\coloneqq
        \prod_{i < \pi(i)}b_i.
    \]
    No degenerate $W(D_{2m})$-module is $(A_\pi,\eta_\pi)$-distinguished.
\end{lemma}
\begin{proof}
    Suppose a degenerate module $E_\lambda^\pm$ is
    $(A_\pi,\eta_\pi)$-distinguished.
    Then $V_{\lambda,\lambda}$ is also
    $(A_\pi,\eta_\pi)$-distinguished.
    Define
    \[
        \eta_I(\mathbf b)\coloneqq\prod_{i\in I}b_i,\qquad
        \mathbf b\in C_2^{2m}, \
        I\subseteq\{1,\ldots,2m\}.
    \]
    By \eqref{eq:type-B-irreps} and Mackey theory, we obtain
    \[
        \eta_\pi = \eta_{I_\pi}|_{A_\pi}
        \quad\text{for some }I_\pi\subseteq\{1,\ldots,2m\},\quad |I_\pi|=m.
    \]
    Set $L_\pi\coloneqq\{i\mid i<\pi(i)\}$ and
    $S_\pi\coloneqq\{i\mid i=\pi(i)\}$. We have $2|L_\pi|+|S_\pi|=2m$, $|S_\pi|>0$.

    For each $i\in L_\pi$, take $\mathbf b\in A_\pi$ such that $b_i = b_{\pi(i)}=-1$ and $b_{j}=1$ for $j\notin \{i,\pi(i)\}$. Then $\eta_\pi(\mathbf b)=-1$, forcing $I_\pi$ to contain exactly one of $\{i,\pi(i)\}$.
    For distinct $i,j\in S_\pi$, take $\mathbf b \in A_\pi$ such that $b_i = b_j = -1$ and $b_{k}=1$ for $k\notin \{i,j\}$.
    Then $\eta_\pi(\mathbf b)=1$, forcing $I_\pi$ to contain both or none of $\{i,j\}$.
    Thus $|I_\pi|\in\{|L_\pi|,|L_\pi|+|S_\pi|\}$, contradicting $|I_\pi|=m$.
\end{proof}

% !TEX root = ../main.tex
\subsection{Family distinction to invariance condition}
\label{subsec:family-distinction-to-invariance-condition}
Choose $\xi=(f,\mu)$ and $E$ as in Theorem~\ref{thm:family-distinction-criterion},
and write $\cL\coloneqq\cL_\xi$.
Retain the stabilizer $W_\xi$ and its character $\eta_\xi$ from
\S\ref{subsec:the-associated-graded-module}.
Choose
$\tilde{\mu}\in X^*(T)$ such that
\[
    \tilde{\mu}|_{T_{f,N}}=\mu,\qquad
    \theta_f\tilde{\mu}=-\tilde{\mu},
\]
as in Lemma~\ref{lem:anti-invariant-lift}.
The lemma also gives $\theta_f^*\cL\cong\cL^{-1}$.
Set $\delta_f\coloneqq(\theta_f)_W$.
Since $Z\coloneqq Z(G)$ is connected,
Lemma~\ref{lem:connected-center-character-weyl-group} gives
$W_\cL=W(\Phi_\cL)$.
Apply \S\ref{subsec:invariance-under-root-system-involutions} with
$\Psi=\Phi_\cL$, $W(\Psi)=W_\cL$, and $\delta=\delta_f$.
Write
\[
    \Phi_\cL=\bigsqcup_{i=1}^n\Phi_{\cL,i},
    \qquad
    W_\cL=\prod_{i=1}^n W_{\cL,i},
    \qquad W_{\cL,i}\coloneqq W(\Phi_{\cL,i}),
\]
where the $\Phi_{\cL,i}$ are the irreducible components.
The involution $\delta_f$ of Proposition~\ref{prop:induced-weyl-group-involution}
preserves $\Phi_\cL$.

\subsubsection{Simple derived groups}
\label{subsubsec:simple-derived-groups}

\begin{lemma}
    \label{lem:involution-away-from-a-d-e6-are-inner}
    Assume that $G_{\opr{der}}$ is simple not of type $A$, $D$, or $E_6$. Then
    $\delta_f|_{W_\cL}$ is inner.
\end{lemma}
\begin{proof}
    The Dynkin diagram of the remaining types
    has no nontrivial automorphisms; see the proof of
    Lemma~\ref{lem:involution-on-component-is-outer-only-in-d2m}.
    Hence
    the action of $\theta_f$ on roots is given by some $w\in W$
    (\cite[Appendix~C]{Milne_algebraic_groups}). Since roots
    generate $X^*(T/Z)$, the actions of $w$ and ${\theta_{f}}$ on
    $T/Z$ agree. Also $\delta_f=\opr{Int}_w$ and $w^2=1$.

    By Lemma~\ref{lem:untwisting-the-center}, there exists a Kummer local system $\cL_0$ on $T$ such that
    \[
        \cL_0|_Z\cong\overline{\bQ}_\ell,\qquad
        {\theta_{f}}^*\cL_0\cong\cL_0^{-1},\qquad
        W_{\cL_0}=W_\cL.
    \]
    The local system $\cL_0$ descends to $T/Z$, where $-1$ acts
    by inversion. Hence
    \[
        w^*\cL_0\cong\cL_0^{-1},\qquad
        (-1)^*\cL_0\cong\cL_0^{-1}.
    \]
    It follows that $a\coloneqq(-1)w$ belongs to
    $W_{\cL_0}=W_\cL$. Since $-1$ is central in $W$, we have
    $a^2=1$ and $\delta_f=\opr{Int}_a$, as required.
\end{proof}

\begin{lemma}
    \label{lem:exchanged-components-have-type-a}
    Assume that $G_{\opr{der}}$ is simple of type $A$, $D$, or $E_6$.
    Every irreducible component of $\Phi_\cL$ exchanged with
    another component by $\theta_f$ has type $A$.
\end{lemma}
\begin{proof}
    In type $A$, write the root-space projection of
    $\tilde\mu/N$ as $\sum_j t_j e_j$.
    A root $e_j-e_l$ belongs to $\Phi_\cL$ precisely when
    $t_j-t_l\in\bZ$. The blocks of equal coordinate residues
    therefore give components of type $A$.

    In type $D_r$, choose coordinates on the span of the roots in which the induced action of $\theta_f$ is a signed permutation. Such coordinates exist also for $D_4$:
    the image of an involution in its diagram group $S_3$ is
    trivial or a transposition, which can be identified with
    the fork exchange, as in
    Lemma~\ref{lem:involution-on-component-is-outer-only-in-d2m}.
    Write the root-space projection of $\tilde\mu/N$ as
    $\sum_j t_j e_j$. Then $e_j\pm e_l\in\Phi_\cL$ precisely
    when $t_j\pm t_l\in\bZ$.
    Residue pairs $\{\epsilon,-\epsilon\}$ with
    $\epsilon\ne-\epsilon$ give type $A$ components.
    The remaining roots lie on the two blocks
    \[
        \{j\mid t_j\equiv0\pmod{\bZ}\},\qquad
        \{j\mid t_j\equiv\tfrac12\pmod{\bZ}\}.
    \]
    The identity $\theta_f\tilde\mu=-\tilde\mu$ implies that
    the signed permutation preserves each of these blocks.
    On a block of size at least four, the roots form one
    preserved component of type $D$.
    Smaller blocks give only type $A$ components.

    In type $E_6$, every component is simply laced.
    A component not of type $A$ has rank at least four.
    Two such components cannot occur in a root system of rank
    at most six, so none can be exchanged with another component.
\end{proof}

\subsubsection{Diagonal distinction}
\label{subsubsec:diagonal-distinction}
\begin{lemma}
    \label{lem:diagonal-subgroup-in-parameter-stabilizer}
    If $E$ is $(W_\xi,\eta_\xi)$-distinguished, then $E$ is
    $\Delta_{\delta_f}^{\mathrm{sw}}$-distinguished.
\end{lemma}
\begin{proof}
    It suffices to prove
    $\Delta_{\delta_f}^{\mathrm{sw}}\subseteq W_\xi$ and
    $\eta_\xi|_{\Delta_{\delta_f}^{\mathrm{sw}}}=1$.
    Fix $(i,j)\in I_{\delta_f}$ and a root $\alpha$ in $\Phi_{\cL,i}$.
    The roots $\alpha$ and $\theta_f(\alpha)$ lie in distinct
    components of $\Phi_\cL$, so they are orthogonal.
    Every $\theta_f$-fixed root $\gamma$ satisfies
    \[
        \langle\tilde\mu,\gamma^\vee\rangle
        =\langle\theta_f\tilde\mu,\theta_f(\gamma^\vee)\rangle
        =-\langle\tilde\mu,\gamma^\vee\rangle=0,
    \]
    and hence belongs to $\Phi_\cL$.
    A root proportional to $\alpha+\theta_f(\alpha)$ would therefore
    belong to $\Phi_\cL$ and be nonorthogonal to both $\alpha$ and
    $\theta_f(\alpha)$, contradicting their distinct components.
    Since $s_\alpha$ sends $\alpha-\theta_f(\alpha)$ to
    $-\alpha-\theta_f(\alpha)$, no root is proportional to
    $\alpha-\theta_f(\alpha)$ either.
    In particular, $\alpha\pm\theta_f(\alpha)\notin\Phi(G,T)$,
    so the groups $\langle U_{\pm\alpha}\rangle$ and
    $\langle U_{\pm{\theta_f(\alpha)}}\rangle$ commute.

    Choose a representative
    $n_\alpha\in N_G(T)\cap\langle U_{\pm\alpha}\rangle$
    of $s_\alpha$. The map $g\mapsto g\theta_f(g)$ from the connected
    group $\langle U_{\pm\alpha}\rangle$ to $H_f$ is a homomorphism,
    so its image lies in $H_f^\circ$. Thus
    $n_\alpha\theta_f(n_\alpha)\in N_{H_f^\circ}(T)$ represents
    $w_{\alpha}\coloneqq s_\alpha s_{\theta_f(\alpha)}$. Both reflections preserve
    $\cL$, so $w_{\alpha}\in W_\cL$.
    The representative in $N_{H_f^\circ}(T)$ gives
    $w_\alpha\times f=f$ under Definition~\ref{def:cross-action}.
    Moreover, $w_{\alpha}^2=1$ and $w_{\alpha}\theta_f=\theta_f w_{\alpha}$.

    Let $\beta\in X^*(T)$ such that $\theta_f\beta =- \beta$.
    Put $m_\beta\coloneqq\langle\beta,\alpha^\vee\rangle\in\bZ$.
    Since the character--cocharacter pairing is $\theta_f$-invariant, we have
    \[
        \beta - w_{\alpha}\beta
        = m_\beta(1-\theta_f)\alpha.
    \]
    Use $\opr{Inv}_f(w_\alpha)$ from Definition~\ref{def:cross-action}.
    As $w_{\alpha}^2=1$, we have
    $\opr{Inv}_f(w_{\alpha})=\opr{Inv}_f(w_{\alpha}^{-1})$.
    The map $\beta\mapsto-w_\alpha\beta$ is an involution of
    $\opr{Inv}_f(w_\alpha)$, since $w_\alpha^2=1$ and
    $w_\alpha\theta_f=\theta_fw_\alpha$.
    A fixed point would satisfy $2\beta=m_\beta(1-\theta_f)\alpha$,
    contradicting the absence of roots proportional to
    $\alpha-\theta_f(\alpha)$.
    Thus $|\opr{Inv}_f(w_\alpha)|$ is even.
    Each two-element orbit has sum
    $\beta-w_{\alpha}\beta=m_\beta(1-\theta_f)\alpha$. Therefore
    \[
        N\sum_{\beta \in \opr{Inv}_f(w_{\alpha}^{-1})} \beta \in N(1-\theta_f)X^*(T).
    \]
    Consequently,
    \[
        \prod_{\beta\in\opr{Inv}_f(w_\alpha^{-1})}\beta(t^N)=1
        \qquad(t\in T_{f,N}).
    \]

    Since $\alpha\in\Phi_\cL$, the description of $\cL$ by
    $\tilde\mu/N+X^*(T)$ gives
    \[
        m_{\tilde{\mu}}\coloneqq\langle\tilde{\mu},\alpha^\vee\rangle\in N\mathbb Z.
    \]
    Therefore
    \[
        \tilde{\mu}-w_\alpha\tilde{\mu}
        =m_{\tilde{\mu}}(1-\theta_f)\alpha
        =N(1-\theta_f)
        \left(\frac{m_{\tilde{\mu}}}{N}\alpha\right).
    \]
    Thus $w_\alpha\mu=\mu$. The real-root correction in
    Definition~\ref{def:cross-action} is trivial, so
    $w_\alpha\times_f\mu=\mu$ and $w_\alpha\times\xi=\xi$.
    Lemma~\ref{lem:stabilizer-character-positive-system-independence}
    gives
    \[
        \eta_\xi(w_\alpha)=(-1)^{|\opr{Inv}_f(w_\alpha)|}=1.
    \]
    As $(i,j)$ ranges over $I_{\delta_f}$ and $\alpha$ ranges over $\Phi_{\cL,i}$, the elements $w_\alpha$ generate $\Delta_{\delta_f}^{\mathrm{sw}}$. Hence $\Delta_{\delta_f}^{\mathrm{sw}}\subseteq W_\xi$ and $\eta_\xi|_{\Delta_{\delta_f}^{\mathrm{sw}}}=1$.
\end{proof}

\subsubsection{Exclusion of degenerate \texorpdfstring{$D_{2m}$}{D2m}-factors}
\label{subsubsec:exclusion-of-degenerate-d2m-factors}
Let $W_{\cL,k}$ be a $\delta_f$-stable component of $W_\cL$
on which $\delta_f$ acts outerly.
Write $E\cong\boxtimes_{i=1}^n E_i$ with
$E_i\in\Irr(W_{\cL,i})$.
Let $\tau=(\mathbf a,\pi)\in W(B_{2m})\setminus W(D_{2m})$
and $\iota_k:W_{\cL,k}\to W(D_{2m})$ be as in Lemma~\ref{lem:involution-on-component-is-outer-only-in-d2m},
and let $A_\pi$ and $\eta_\pi$ be as defined in
Lemma~\ref{lem:no-degenerate-in-induced-d2m-module}.

\begin{lemma}
    \label{lem:even-sign-subgroup-in-parameter-stabilizer}
    Assume that $G_{\opr{der}}$ is simple and simply connected
    of type $D$ or $E_6$.
    If $E$ is $(W_\xi,\eta_\xi)$-distinguished,
    then $E_k$ is non-degenerate.
\end{lemma}

By Lemma~\ref{lem:no-degenerate-in-induced-d2m-module},
it is enough to show
\[
    \iota_k^{-1}(A_\pi)\subseteq W_\xi,
    \qquad
    \eta_\xi(\iota_k^{-1}(x))=\eta_\pi(x)
    \quad(x\in A_\pi).
\]
\begin{lemma}
    \label{lem:even-d-component-is-levi}
    Assume that $G_{\opr{der}}$ is simple and simply connected
    of type $D$ or $E_6$.
    There exists a $\theta_f$-stable Levi subgroup $M_k\subseteq G$
    containing $T$ such that $\Phi(M_k,T)=\Phi_{\cL,k}$.
    Its derived group $M_{k,\opr{der}}$ is $\theta_f$-stable and
    isomorphic to $\opr{Spin}_{4m}$.
\end{lemma}
\begin{proof}
    By Lemma~\ref{lem:involution-on-component-is-outer-only-in-d2m},
    $\Phi_{\cL,k}$ has type $D_{2m}$ with $m\geq2$. We claim that
    \[
        \Phi(G,T)\cap\Span_{\bR}(\Phi_{\cL,k})=\Phi_{\cL,k}.
    \]
    In type $D_r$, the proof of
    Lemma~\ref{lem:exchanged-components-have-type-a} shows that
    $\Phi_{\cL,k}$ consists of all roots $\pm e_j\pm e_l$
    supported on one coordinate block.
    Its real span is the coordinate subspace of that block.
    Hence
    $\Phi(G,T)\cap\Span_{\bR}(\Phi_{\cL,k})=\Phi_{\cL,k}$.

    In type $E_6$, the rank bound gives $2m\in\{4,6\}$.
    Normalize roots to have squared length two.
    Any ambient root in $\Span_{\bR}(\Phi_{\cL,k})$ lies in the
    $D_{2m}$ weight lattice $\bZ^{2m}\cup(\bZ+\tfrac12)^{2m}$.
    Half-integral vectors have squared length in $m/2+2\bZ$, which excludes two.
    The integral vectors of squared length two are precisely $\pm e_j\pm e_l$,
    proving the claim.

    Set
    \[
        S_k\coloneqq\left(\bigcap_{\alpha\in\Phi_{\cL,k}}\ker\alpha\right)^\circ_{\opr{red}},
        \qquad M_k\coloneqq C_G(S_k).
    \]
    By \cite[\S25.25]{Milne_algebraic_groups}, $M_k$ is a Levi subgroup;
    the claim gives $\Phi(M_k,T)=\Phi_{\cL,k}$.
    Since $T$ and $\Phi_{\cL,k}$ are $\theta_f$-stable, so are $S_k$, $M_k$,
    and $M_{k,\opr{der}}$.
    Apply \cite[Corollary~9.5.11]{conrad2020reductive} to the Levi subgroup
    $M_k\cap G_{\opr{der}}$ of the simply connected group $G_{\opr{der}}$.
    It follows that $M_{k,\opr{der}}$ is simply connected of type $D_{2m}$,
    hence isomorphic to $\opr{Spin}_{4m}$.
\end{proof}

For the fixed $\tau=(\mathbf a,\pi)\in W(B_{2m})\setminus W(D_{2m})$, set
\[
    \begin{gathered}
        P_\tau\coloneqq\{i\mid\pi(i)=i,\ a_i=1\},\qquad
        C_\tau\coloneqq\{i\mid\pi(i)=i,\ a_i=-1\},\\
        L_\pi\coloneqq\{i\mid i<\pi(i)\}.
    \end{gathered}
\]
As $\tau^2=1$, the sets $P_\tau$ and $C_\tau$ have odd cardinality.

\begin{lemma}
    \label{lem:even-sign-generators-have-connected-fixed-representatives}
    Assume that $G_{\opr{der}}$ is simple and simply connected
    of type $D$ or $E_6$.
    Every $x\in \iota_k^{-1}A_\pi$ has a representative in
    $N_{H_f^\circ}(T)\cap M_{k,\opr{der}}$.
\end{lemma}
\begin{proof}
    By Lemma~\ref{lem:even-d-component-is-levi},
    $M_{k,\opr{der}}\cong\opr{Spin}_{4m}$ is $\theta_f$-stable.
    By \cite[Theorem~6.1.17]{Con14}, choose the isomorphism with
    $\opr{Spin}_{4m}$ to induce $\iota_k$, also in type $D_4$.
    Composing with the standard double covering $\opr{Spin}(V)\to\SO(V)$ gives
    \[
        q:M_{k,\opr{der}}\cong\opr{Spin}(V)\longrightarrow\opr{SO}(V),
        \qquad V=\bigoplus_{i=1}^{2m}(ku_i\oplus kv_i),
        \qquad (u_i,v_j)=\delta_{ij},
    \]
    with $(u_i,u_j)=(v_i,v_j)=0$, so that
    $q(T\cap M_{k,\opr{der}})$ is the diagonal torus.
    In the $\iota_k$-coordinates, the characters of $T\cap M_{k,\opr{der}}$
    trivial on $\ker(q)$ form $\bZ^{2m}$.
    Since $\tau$ preserves this lattice, $\theta_f$ preserves $\ker(q)$.

    The involution $\underline{\theta}_f$ of $\SO(V)$ induced by $\theta_f$
    differs from the standard fork exchange
    by an inner automorphism (\cite[Proposition~7.1.6]{Con14}).
    The fork exchange is induced by an orthogonal transformation
    of determinant $-1$. Thus $\underline{\theta}_f = \opr{Int}_r$
    for some $r\in O(V)$ with $\det(r)=-1$.
    Since $\theta_f^2=1$, we have $r^2=\pm1$.
    If $r^2=-1$, the two eigenspaces of $r$ are isotropic and paired
    non-degenerately, forcing $\det(r)=1$. Thus $r^2=1$.
    Put $V_\pm\coloneqq\ker(r\mp1)$. Then
    \[
        (\opr{SO}(V)^{\opr{Int}_r})^\circ
        =\opr{SO}(V_+)\times\opr{SO}(V_-).
    \]

    As $r$ normalizes the torus and preserves the bilinear form, the chosen labeling gives $r(\{ku_i,kv_i\}) = \{ku_{\pi(i)},kv_{\pi(i)}\}$.
    Therefore, for any $i=\pi(i)$ one can choose hyperbolic bases such that
    \[
        r|_{ku_p\oplus kv_p}=\delta_p\id\quad(p\in P_\tau,\ \delta_p\in\{\pm1\}),
        \qquad ru_c=v_c,\quad rv_c=u_c\quad(c\in C_\tau).
    \]
    For each pair $i<\pi(i)$, choose
    $\{r(u_i),r(v_i)\}$ as the basis of
    $\Span\{u_{\pi(i)},v_{\pi(i)}\}$.
    These choices may interchange $u_i$ and $v_i$,
    but preserve all coordinate sign changes,
    hence the identification of $A_\pi$. In this basis
    $r$ has the following form
    \[
        \bigoplus_{p\in P_\tau}\delta_p I_2
        \;\oplus\;
        \bigoplus_{c\in C_\tau}
        \begin{pmatrix}0&1\\1&0\end{pmatrix}
        \;\oplus\;
        \bigoplus_{i\in L_\pi}
        \begin{pmatrix}
            0 & 0 & 1 & 0 \\
            0 & 0 & 0 & 1 \\
            1 & 0 & 0 & 0 \\
            0 & 1 & 0 & 0
        \end{pmatrix}.
    \]

    Let $J_i\in O(V)$ interchange $u_i,v_i$ and fix all other planes.
    For $x=(\mathbf b,1)\in A_\pi$, set
    \[
        J_x\coloneqq\prod_{b_i=-1}J_i.
    \]
    Since $b_i=b_{\pi(i)}$ and
    $\prod_i b_i=1$, the element $J_x$ commutes with $r$
    and belongs to $\opr{SO}(V)$.
    It normalizes the diagonal torus and induces $x$.

    Fix $c_0\in C_\tau$, and let $Z_{c_0}$ act by $-1$ on
    $ku_{c_0}\oplus kv_{c_0}$ and trivially on the other planes.
    The torus element $Z_{c_0}$ commutes with $r$ and has
    determinant $-1$ on each of $V_+$ and $V_-$.
    Since $\det(J_x)=1$, the determinants of $J_x$ on $V_+$ and
    $V_-$ agree. Thus either $J_x$ or $J_xZ_{c_0}$ lies in
    $\opr{SO}(V_+)\times\opr{SO}(V_-)$ and still induces $x$.

    The central isogeny $q$ induces a surjection
    \[
        (M_{k,\opr{der}}^{\theta_f})^\circ
        \longrightarrow\opr{SO}(V_+)\times\opr{SO}(V_-).
    \]
    Indeed, $\ker(q)$ has order two and the characteristic is odd,
    so the differential identifies the fixed Lie algebras.
    The image is therefore a closed connected subgroup of full dimension.
    Lift the chosen representative along this surjection.
    The lift lies in $H_f^\circ\cap M_{k,\opr{der}}$ and normalizes
    $T\cap M_{k,\opr{der}}$, the full preimage of the diagonal torus.
    It centralizes $Z(M_k)^\circ$, so it normalizes
    $T=Z(M_k)^\circ(T\cap M_{k,\opr{der}})$.
    Its Weyl-group image is $\iota_k^{-1}(x)$, as required.
\end{proof}

\begin{lemma}
    \label{lem:even-sign-parameter-character-lattice-inclusion}
    Assume that $G_{\opr{der}}$ is simple and simply connected
    of type $D$ or $E_6$.
    Then $\iota_k^{-1}(A_\pi)\subseteq W_\xi$.
\end{lemma}
\begin{proof}
    Fix $x=(\mathbf b,1)\in A_\pi$ and put $w\coloneqq\iota_k^{-1}(x)$.
    Lemma~\ref{lem:even-sign-generators-have-connected-fixed-representatives}
    gives a representative of $w$ in $N_{H_f^\circ}(T)$.
    Thus $w\times f=f$ by Lemma~\ref{lem:cross-action-stabilizer-criterion}.
    For a $W$- and $\theta_f$-invariant inner product, let $\kappa$ be the
    projection of $\tilde{\mu}/N$ onto $\Span_{\bR}(\Phi_{\cL,k})$,
    in the coordinates defining $\iota_k$.
    The projection commutes with $\theta_f$.
    Since $\theta_f\tilde\mu=-\tilde\mu$ and
    $\iota_k\theta_f\iota_k^{-1}=\tau$, we have $\tau\kappa=-\kappa$.
    The coroot pairings of $\tilde\mu/N$ on $\Phi_{\cL,k}$ are integral,
    so $\kappa_i\pm \kappa_j\in\bZ$ for $i\ne j$.
    Fix $p_0\in P_\tau$. The equality $\tau \kappa=-\kappa$ gives $\kappa_{p_0}=0$,
    hence $\kappa\in\bZ^{2m}$. Since $x\tau=\tau x$, we obtain
    \[
        w\tilde{\mu}-\tilde{\mu}=N(x\kappa-\kappa)
        =-N(1-\theta_f)
        \sum_{\substack{i\ne p_0\\b_i=-1}}
        \kappa_i(e_i-e_{p_0})
        \in N(1-\theta_f)X^*(T).
    \]
    Here $w$ fixes the orthogonal complement of $\Span_{\bR}(\Phi_{\cL,k})$,
    and each $e_i-e_{p_0}$ is a root. Thus $w\mu=\mu$.

    Choose a positive system $P$ making $\Phi_{\cL,k}$ a standard Levi
    subsystem with positive roots $e_i\pm e_j$ for $i<j$.
    All inversions of $w$ relative to $P$ lie in $\Phi_{\cL,k}$.
    The real inversion set relative to $P$ is
    \begin{equation}
        \label{eq:even-sign-real-inversions}
        \begin{aligned}
            \opr{Inv}_{f,P}(w)
            ={} & \{e_i-a_i e_{\pi(i)}\mid i\in L_\pi,\ b_i=-1\}         \\
                & \sqcup\{e_c\pm e_d\mid c,d\in C_\tau,\ c<d,\ b_c=-1\}.
        \end{aligned}
    \end{equation}
    Each root in the first set contributes $e_i-a_i e_{\pi(i)}$
    to the inversion sum, and each pair $e_c\pm e_d$ contributes $2e_c$.
    These contributions lie in $(1-\theta_f)X^*(T)$, since
    \[
        e_i-a_i e_{\pi(i)}=(1-\tau)(e_i-e_{p_0}),
        \qquad 2e_c=(1-\tau)(e_c-e_{p_0}).
    \]
    Let $\gamma$ be a $\theta_f$-real root.
    Its projection onto $\Span_{\bR}(\Phi_{\cL,k})$ has integral
    coroot pairings and zero $p_0$-coordinate, so its coordinates
    are integral.
    Applying the calculation of $w\tilde\mu-\tilde\mu$ with
    $\tilde\mu$ replaced by $N\gamma$ gives
    \[
        (1-w)\gamma\in(1-\theta_f)X^*(T).
    \]
    Lemma~\ref{lem:positive-system-comparison} therefore gives
    $w\times_{f,P}\mu=w\times_f\mu$.
    The identities for $e_i-a_i e_{\pi(i)}$ and $2e_c$ show
    that the correction factor in $w\times_{f,P}\mu$ is trivial
    on $T_{f,N}$, since $w^{-1}=w$.
    Together with $w\mu=\mu$, Definition~\ref{def:cross-action}
    gives $w\times_f\mu=w\times_{f,P}\mu=\mu$.
    Hence $w\times\xi=\xi$, so $w\in W_\xi$.
\end{proof}

\begin{lemma}
    \label{lem:even-sign-parameter-line-character}
    Assume that $G_{\opr{der}}$ is simple and simply connected
    of type $D$ or $E_6$.
    Then
    \[
        \eta_\xi(\iota_k^{-1}(x))=\eta_\pi(x)
        \quad(x\in A_\pi).
    \]
\end{lemma}
\begin{proof}
    Write $x=(\mathbf b,1)\in A_\pi$ and $w\coloneqq\iota_k^{-1}(x)$.
    Choose the Levi positive system $P$ used in
    \eqref{eq:even-sign-real-inversions}.
    Lemma~\ref{lem:stabilizer-character-positive-system-independence} gives
    $\eta_\xi(w)=(-1)^{|\opr{Inv}_{f,P}(w)|}$.
    The pairs $e_c\pm e_d$ contribute an even number of inversions.
    Equation~\eqref{eq:even-sign-real-inversions} therefore gives
    \[
        \eta_\xi(w)=(-1)^{|\{i\in L_\pi\mid b_i=-1\}|}
        =\prod_{i\in L_\pi}b_i=\eta_\pi(x).
    \]
\end{proof}

\begin{proof}[Proof of Lemma~\ref{lem:even-sign-subgroup-in-parameter-stabilizer}]
    By Lemmas~\ref{lem:even-sign-parameter-character-lattice-inclusion} and~\ref{lem:even-sign-parameter-line-character}
    we have
    \[
        \iota_k^{-1}(A_\pi)\subseteq W_\xi,
        \qquad
        \eta_\xi(\iota_k^{-1}(x))=\eta_\pi(x)
        \quad(x\in A_\pi).
    \]
    The subgroup $\iota_k^{-1}(A_\pi)$ acts only on the $k$th factor of $E$.
    A nonzero $\eta_\xi$-eigenvector in $E$ therefore gives
    $(A_\pi,\eta_\pi)$-distinction of $E_k\circ\iota_k^{-1}$.
    The permutation $\pi$ has a fixed point, as noted after
    Lemma~\ref{lem:involution-on-component-is-outer-only-in-d2m}.
    Lemma~\ref{lem:no-degenerate-in-induced-d2m-module} therefore
    implies that $E_k$ is non-degenerate.
\end{proof}

% !TEX root = ../main.tex
\section{Proofs of the main theorems}
\label{sec:proofs-of-the-main-theorems}
For triples of simple connected-center type, we prove that every
uniform almost character in a family containing a distinguished
character is fixed by twisted duality.
We then treat Galois triples separately. Proposition~\ref{prop:connected-center-reduction}
then gives both main theorems.

% !TEX root = ../main.tex
\subsection{Invariance of uniform almost characters}
\label{subsec:uniform-almost-characters-invariance}
\begin{theorem}
    \label{thm:distinguished-family-uniform-invariance}
    Assume that $(G,F,\theta)$ is of simple connected-center type.
    Let $\cF_G$ be a Lusztig family containing an
    $(H^\circ)^F$-distinguished irreducible character.
    Every uniform almost character $\rho^\cL_{\widetilde E}$
    with $\opr{Fam}^{G^F}_\cL(E)=\cF_G$ satisfies
    \[
        \bigl(\rho^\cL_{\widetilde E}\bigr)^\vee\circ\theta
        =\rho^\cL_{\widetilde E}.
    \]
\end{theorem}
\begin{proof}
    Choose an $(H^\circ)^F$-distinguished $\chi\in\cF_G$.
    Theorem~\ref{thm:family-distinction-criterion}
    supplies a parameter $\xi=(f,\mu)$ and a simple module $E_0$
    such that
    \[
        \cF_G=\opr{Fam}^{G^F}_{\cL_\xi}(E_0),\qquad
        \Hom_{W_\xi}(\eta_\xi,E_0)\ne0.
    \]
    Here $\theta_f(T)=T$, the orbit $W\cdot{\cL_\xi}$ is $F$-stable,
    and $E_0\in\Irr(W_{\cL_\xi})^{\gamma_v}$.
    Let $\cF_W$ be the Weyl-group family containing $E_0$.

    By Lemma~\ref{lem:anti-invariant-lift}, choose
    $\tilde{\mu}\in X^*(T)$ with
    $\tilde{\mu}|_{T_{f,N}}=\mu$ and
    $\theta_f\tilde{\mu}=-\tilde{\mu}$.
    The lemma also gives
    $\theta_f^*{\cL_\xi}\cong{\cL_\xi}^{-1}$.
    Set $\delta_f\coloneqq(\theta_f)_W$.
    By Lemma~\ref{lem:connected-center-character-weyl-group},
    $W_{\cL_\xi}=W(\Phi_{\cL_\xi})$, and $\theta_f$ preserves $\Phi_{\cL_\xi}$.
    Outside types $A$, $D$, and $E_6$,
    Lemma~\ref{lem:involution-away-from-a-d-e6-are-inner}
    gives $E\cong E\circ\delta_f$ for every $E\in\cF_W$.

    Suppose now that $G$ has type $A$, $D$, or $E_6$.
    Lemma~\ref{lem:diagonal-subgroup-in-parameter-stabilizer}
    gives $E_0^{\Delta_{\delta_f}^{\mathrm{sw}}}\ne0$.
    Write $E_0\cong\boxtimes_i E_{0,i}$ for
    $W_{\cL_\xi}=\prod_i W_{{\cL_\xi},i}$.
    By Proposition~\ref{prop:non-fixed-module-criteria}, it suffices
    for $E_0\cong E_0\circ\delta_f$ to exclude degenerate
    $E_{0,k}$ on preserved components with outer action.
    In type $A$, every component of $\Phi_{\cL_\xi}$ has type $A$,
    so no preserved component admits an outer Weyl-group action.
    In types $D$ and $E_6$,
    Lemma~\ref{lem:even-sign-subgroup-in-parameter-stabilizer}
    excludes the degenerate constituents.
    Thus $E_0\cong E_0\circ\delta_f$.
    Lemma~\ref{lem:exchanged-components-have-type-a} and
    Lemma~\ref{lem:fixed-label-implies-fixed-family} now give
    $E\cong E\circ\delta_f$ for every $E\in\cF_W$.

    The uniform almost characters belonging to $\cF_G$
    are indexed by
    $\cF_W\cap\Irr(W_{\cL_\xi})^{\gamma_v}$.
    Apply Lemma~\ref{lem:unif-alm-char-theta-inv-functoriality}
    to each of these labels.
\end{proof}

% !TEX root = ../main.tex
\subsection{Galois type}
\label{subsec:galois-type}
Lemma~\ref{lem:lp-for-galois-type} and
Proposition~\ref{prop:classical-galois-twisted-self-duality} are
consequences of \cite[Theorem~1.4(a)--(c)]{lusztig2000g} and
compatibility of Jordan decomposition with twisted duality.
We give a proof of family invariance using the occurrence criterion
of Theorem~\ref{thm:family-distinction-criterion} and the diagonal-distinction
argument of Lemma~\ref{lem:diagonal-subgroup-in-parameter-stabilizer}.

For every Galois-type triple, Definition~\ref{def:galois-type-triple}
identifies $H$ with the diagonal subgroup
$\{(g,g):g\in G'\}\subseteq G'\times G'$.
Since $G'$ is connected, $H=H^\circ$ and hence
$H^F=(H^\circ)^F$.

\begin{lemma}
    \label{lem:lp-for-galois-type}
    Let $(G,F,\theta)$ be a Galois-type triple arising from a group
    $G'$ of simple connected-center type.
    Every $H^F$-distinguished irreducible character $\chi$ of $G^F$
    lies in the same Lusztig family as $\chi^\vee\circ\theta$.
\end{lemma}
\begin{proof}
    Write $(G,F,\theta)$ as in Definition~\ref{def:galois-type-triple},
    so $H=\Delta G'$ is connected.
    Theorem~\ref{thm:family-distinction-criterion} supplies $\xi=(f,\mu)$ and $E$
    with $\chi\in\opr{Fam}^{G^F}_{\cL_\xi}(E)$ and
    $\Hom_{W_\xi}(\eta_\xi,E)\ne0$.
    Set $\delta_f\coloneqq(\theta_f)_W$.
    The involution $\theta_f$ exchanges the two factors of $G$,
    hence all irreducible components of $\Phi(G,T)$ are exchanged.
    Lemma~\ref{lem:diagonal-subgroup-in-parameter-stabilizer} gives
    $E^{\Delta_{\delta_f}^{\mathrm{sw}}}\ne0$.
    Every component of $\Phi_{\cL_\xi}$ is also exchanged, so
    Proposition~\ref{prop:non-fixed-module-criteria} gives
    $E\cong E\circ\delta_f$.
    By Lemma~\ref{lem:anti-invariant-lift},
    $\theta_f^*\cL_\xi\cong\cL_\xi^{-1}$.
    Proposition~\ref{prop:induced-weyl-group-involution} now proves
    the family assertion.
\end{proof}

In the application of the following lemma,
$\Gamma$ is Lusztig's finite group attached to the
$F'$-stable family $\cF_W$ of the Weyl group $W_s$.
The automorphism $\tau$ is induced by $F'$.

\begin{lemma}[cf. {\cite[\S1.3 and Theorem~1.4(c)]{lusztig2000g}}]
    \label{lem:square-root-parameter-invariance}
    Let $\Gamma \cong C_2^r$  where $C_2$ is the group
    of two elements and $r\in \bZ_{\ge0}$.
    Let $\tau\in\operatorname{Aut}(\Gamma)$.
    Put $\dot\Gamma\coloneqq\Gamma\rtimes\langle\gamma\rangle$,
    where $\langle\gamma\rangle$ is infinite cyclic and
    $\gamma a\gamma^{-1}=\tau(a)$ for $a\in\Gamma$.
    For $z\in\dot\Gamma$, denote its centralizer in $\Gamma$
    by $Z_\Gamma(z)$.
    Fix $x'\in\Gamma\gamma$ and put $x\coloneqq (x')^2$.
    Suppose that $\xi\in\Irr(Z_\Gamma(x))$ occurs in
    $\Ind_{Z_\Gamma(x')}^{Z_\Gamma(x)}1$.
    Then the pairs
    \[
        (x,\xi)
        \quad\text{and}\quad
        (\gamma x\gamma^{-1},\xi\circ\tau^{-1})
    \]
    are $\Gamma$-conjugate.
\end{lemma}
\begin{proof}
    Since $\Gamma$ is abelian, $Z_\Gamma(x')=\Gamma^\tau$
    and $Z_\Gamma(x)=\Gamma^{\tau^2}$.
    The character $\xi$ is one-dimensional, so Frobenius reciprocity gives
    $\xi|_{Z_\Gamma(x')}=1$.
    For $h\in Z_\Gamma(x)$, we have $h\tau(h)\in Z_\Gamma(x')$
    and $h^2=1$, so
    \[
        \xi(\tau(h))=\xi(h)^{-1}=\xi(h).
    \]
    Thus $\xi\circ\tau^{-1}=\xi$ on $Z_\Gamma(x)$.
    Write $x'=a\gamma$ with $a\in\Gamma$.
    Since $(x')^2=x$, we have
    \[
        a(\gamma x\gamma^{-1})a^{-1}=x.
    \]
    Conjugation by $a$ acts trivially on $\Gamma$.
    Hence it sends $(\gamma x\gamma^{-1},\xi\circ\tau^{-1})$
    to $(x,\xi)$.
\end{proof}

\begin{lemma}
    \label{lem:classical-galois-unipotent-parameter}
    Let $G'$ be of simple connected-center classical type, with
    Frobenius $F'$ and dual group $\cG'$.
    Write $F'$ also for the dual Frobenius.
    Let $\chi\in\Irr(G'^{F'^2})$ be $G'^{F'}$-distinguished.
    There is a semisimple element $s\in\cG'^{F'^2}$ for $\chi$
    with $F'(s)=s^{-1}$.
    Let $J(\chi)$ be the unipotent irreducible character of
    $C_{\cG'}(s)^{F'^2}$ associated to $\chi$ by Jordan
    decomposition.
    Then
    \[
        J(\chi)^\vee=J(\chi),
        \qquad
        J(\chi)\circ F'=J(\chi).
    \]
\end{lemma}
\begin{proof}
    Since $G'$ has connected center and is simple modulo its
    center, \cite[Theorem~1.4(a)]{lusztig2000g} supplies
    $s\in\cG'^{F'^2}$ with $F'(s)=s^{-1}$.
    The group $C_{\cG'}(s)$ is connected and $F'$-stable, and its root
    system has only classical components.
    Choose an $F'$-stable Borel--torus pair in $C_{\cG'}(s)$, and let
    $W_s$ be its Weyl group.
    The family $\cF_W$ associated to $\chi$ is $F'$-stable by
    \cite[Theorem~1.4(b)]{lusztig2000g}.

    By Lemma~\ref{lem:classical-unipotent-uniform-separation}, the unipotent
    characters of $C_{\cG'}(s)^{F'^2}$ are determined by their
    uniform coefficients.
    Since every $R_T^{C_{\cG'}(s)}(1)$ is fixed by contragredience,
    $J(\chi)$ and $J(\chi)^\vee$ have the same uniform
    coefficients. Thus $J(\chi)^\vee=J(\chi)$.

    Let $\Gamma$ be the group attached to $\cF_W$, and let
    $\tau$ be the automorphism induced by $F'$.
    In classical type, $\Gamma\cong C_2^r$ for some $r\geq0$,
    and $|\cF_W|\ne2$.
    Use the semidirect product $\dot\Gamma$ of
    Lemma~\ref{lem:square-root-parameter-invariance}.
    In the parametrization of
    \cite[\S\S1.2--1.3]{lusztig2000g}, write $(x,\xi)$ for
    the parameter of $\chi$, where
    \[
        x\in\Gamma\gamma^2,
        \qquad
        \xi\in\Irr(Z_\Gamma(x)).
    \]
    Jordan decomposition identifies this parameter with
    the parameter of $J(\chi)$
    \cite[Theorem~7.1(i)]{digneMichel1990}.
    Set $S_x\coloneqq\{x'\in\Gamma\gamma:(x')^2=x\}$.
    The $Z_\Gamma(x)$-module with basis the conjugation orbit
    of $x'\in S_x$ is isomorphic to
    $\Ind_{Z_\Gamma(x')}^{Z_\Gamma(x)}1$.
    Thus \cite[Theorem~1.4(c)]{lusztig2000g} gives
    \[
        0<
        \langle\chi,1\rangle_{G'^{F'}}
        =
        \sum_{[x']\in Z_\Gamma(x)\backslash S_x}
        \langle \xi,\Ind_{Z_\Gamma(x')}^{Z_\Gamma(x)}1 \rangle_{Z_\Gamma(x)}.
    \]
    Hence $\xi$ occurs in $\Ind_{Z_\Gamma(x')}^{Z_\Gamma(x)}1$
    for some $x'\in S_x$.
    Lemma~\ref{lem:square-root-parameter-invariance} therefore
    shows that the parameter of $J(\chi)$ is $F'$-fixed.

    The unipotent parametrization is compatible with
    $F'$-pullback.
    Indeed, Deligne--Lusztig induction commutes with
    $F'$-pullback \cite[Corollary~9.2]{digneMichel1990},
    and the pairing defining uniform coefficients is
    invariant under the induced action on parameters
    \cite[Theorem~4.23]{lusztig1984characters}.
    Since uniform coefficients determine classical
    unipotent characters, the fixed parameter gives
    $J(\chi)\circ F'=J(\chi)$.
\end{proof}

\begin{proposition}
    \label{prop:classical-galois-twisted-self-duality}
    Let $(G,F,\theta)$ be a Galois-type triple arising from
    a group $G'$ of simple connected-center classical type.
    Every $H^F$-distinguished irreducible character $\chi$
    of $G^F$ satisfies
    \[
        \chi^\vee\circ\theta=\chi.
    \]
\end{proposition}
\begin{proof}
    Lemma~\ref{lem:cyclic-projection-preserves-families}
    identifies $(G^F,H^F,\theta)$ with
    \[
        (G'^{F'^2},G'^{F'},F'|_{G'^{F'^2}}).
    \]
    Write $\chi$ also for the corresponding character of
    $G'^{F'^2}$.
    By Lemma~\ref{lem:classical-galois-unipotent-parameter},
    choose a semisimple parameter $s$ with $F'(s)=s^{-1}$.
    Its Jordan correspondent $J(\chi)$ is self-dual and
    $F'$-invariant.

    Jordan decomposition commutes with Frobenius pullback
    by \cite[Theorem~7.1]{digneMichel1990}, and with
    contragredience by
    \cite[proof of Theorem~4.1]{srinivasanVinroot2015real}.
    Thus $\chi^\vee\circ F'$ has semisimple parameter
    $F'(s)^{-1}=s$, and
    \[
        J(\chi^\vee\circ F')
        =J(\chi)^\vee\circ F'
        =J(\chi).
    \]
    Injectivity of Jordan decomposition gives
    $\chi^\vee\circ F'=\chi$.
\end{proof}

% !TEX root = ../main.tex
\subsection{Family invariance and equality in classical type}
\label{subsec:family-invariance-and-equality-in-classical-type}
\begin{theorem}\label{thm:lp-conjecture}
    Let $G$ be a connected reductive group over the algebraic closure
    of a finite field of odd characteristic.
    Let $F$ be a Frobenius endomorphism and $\theta$ an involution
    commuting with $F$. Put $H\coloneqq G^\theta$.
    If $\chi\in\Irr(G^F)$ is $H^F$-distinguished, then
    $\chi$ and $\chi^\vee\circ\theta$ lie in the same Lusztig family.
\end{theorem}
\begin{proof}
    By Proposition~\ref{prop:connected-center-reduction}(1), it suffices
    to treat simple connected-center triples and Galois triples
    arising from groups of simple connected-center type.
    Lemma~\ref{lem:lp-for-galois-type} settles the Galois cases.
    Assume henceforth that $(G,F,\theta)$ is of simple connected-center type.
    Let $\cF_G$ be the family of $\chi$, and choose a uniform
    almost character $\rho^\cL_{\widetilde E}$ with
    $\langle\chi,\rho^\cL_{\widetilde E}\rangle_{G^F}\ne0$.
    Twisted duality preserves the character pairing.
    Theorem~\ref{thm:distinguished-family-uniform-invariance} gives
    \[
        \bigl\langle\chi^\vee\circ\theta,
        \rho^\cL_{\widetilde E}\bigr\rangle_{G^F}
        =\langle\chi,\rho^\cL_{\widetilde E}\rangle_{G^F}\ne0.
    \]
    Thus $\chi^\vee\circ\theta$ also belongs to $\cF_G$.
\end{proof}

\begin{theorem}
    \label{thm:classical-twisted-self-duality}
    Suppose that $Z(G)$ is connected and every simple
    factor of $G_{\opr{der}}$ is of classical type.
    Every $(H^\circ)^F$-distinguished $\chi\in\Irr(G^F)$ satisfies
    \[
        \chi^\vee\circ\theta=\chi.
    \]
    If $G$ is of simple connected-center type, every irreducible
    character in a Lusztig family containing an $(H^\circ)^F$-distinguished
    character satisfies the same equality.
\end{theorem}
\begin{proof}
    Apply Proposition~\ref{prop:connected-center-reduction}(2) with
    the simple derived types restricted to classical types.
    Proposition~\ref{prop:classical-galois-twisted-self-duality}
    gives equality for Galois triples
    arising from groups of simple connected-center classical type.
    It remains to treat simple connected-center classical triples.
    We prove the stronger second assertion.
    Let $\cF_G$ be such a family, and let $\psi\in\cF_G$.
    Choose a uniform almost character belonging to $\cF_G$
    with nonzero coefficient at $\psi$.
    Theorem~\ref{thm:distinguished-family-uniform-invariance}
    fixes that uniform almost character, so $\psi^\vee\circ\theta$
    also belongs to $\cF_G$.
    The same theorem gives
    \[
        \bigl\langle\psi^\vee\circ\theta,
        \rho^\cL_{\widetilde E}\bigr\rangle_{G^F}
        =\langle\psi,\rho^\cL_{\widetilde E}\rangle_{G^F}
    \]
    for every uniform almost character belonging to $\cF_G$.
    Both coefficients vanish for uniform almost characters outside
    $\cF_G$, since $\psi,\psi^\vee\circ\theta\in\cF_G$.
    Lemma~\ref{lem:classical-uniform-separation} yields
    $\psi^\vee\circ\theta=\psi$.
\end{proof}

% !TEX root = ../main.tex
\section{Failure of equality in disconnected center}
\label{sec:failure-of-equality-in-disconnected-center}
We give examples with disconnected center where
$\chi\ne\chi^\vee\circ\theta$ under the hypotheses of
Theorem~\ref{intthm:lp-conjecture}.
In Lemma~\ref{lem:lp-counterexample-sl2} and
Proposition~\ref{prop:lp-counterexample-disconnected-center},
$\theta$ is conjugation by an element of $G^F$,
so $\chi^\vee\circ\theta=\chi^\vee$ for every character $\chi$ of $G^F$.

\begin{lemma}
    \label{lem:lp-counterexample-sl2}
    Let $G=\SL_2$ over $\bF_q$ with $q\equiv3\pmod4$,
    and let $F$ be its Frobenius endomorphism.
    Let $\theta$ be conjugation by $s\coloneqq\left(\begin{smallmatrix}
                0 & -1 \\ 1 & 0
            \end{smallmatrix}\right)$,
    and set $H\coloneqq G^\theta$.
    For every nontrivial additive character
    $\psi:\bF_q\to\overline{\bQ}_\ell^\times$, the odd Weil representation
    $\omega_\psi^-$ of $G^F$ is irreducible, cuspidal,
    $H^F$-distinguished, and not self-dual.
    As $\theta$ is inner,
    \[
        \omega_\psi^-
        \not\cong
        (\omega_\psi^-)^\vee\circ \theta
        \cong
        (\omega_\psi^-)^\vee.
    \]

\end{lemma}
\begin{proof}
    Let $\omega_\psi$ be the Weil representation on functions on $\bF_q$,
    and let $\omega_\psi^-$ be its restriction to odd functions.
    Let $U$ be the upper unitriangular subgroup of $G$.
    For $b\in\bF_q$, set $u(b)\coloneqq\left(\begin{smallmatrix}
                1 & b \\ 0 & 1
            \end{smallmatrix}\right)$.
    For $a\in\bF_q$, define $\psi_a(u(b))\coloneqq\psi(ab)$.
    Choose representatives $R\subset\bF_q^\times$ for
    $\bF_q^\times/\{\pm1\}$.
    We use the Schr\"odinger model of
    \cite[\S2.2, (2.1)]{liuZhang2018fourierJacobi}.
    In particular,
    \[
        \omega_\psi^-|_{U^F}
        =\bigoplus_{x\in R}
        \left(\Span\{\delta_x-\delta_{-x}\},\psi_{x^2}\right).
    \]

    Character orthogonality implies that every nonzero
    $G^F$-invariant subspace contains a $U^F$-weight line.
    The diagonal torus permutes these lines transitively, so
    $\omega_\psi^-$ is irreducible.
    Its $U^F$-coinvariants vanish, hence it is cuspidal.

    The group $T_{\opr{ns}}\coloneqq H$ is a nonsplit maximal torus,
    and $H^F$ is cyclic of order $q+1$.
    A generator of $H^F$ is $G^F$-conjugate to its inverse.
    Thus inverse characters of $H^F$ occur with equal multiplicities
    in $\omega_\psi^-|_{H^F}$.
    Since $\dim\omega_\psi^-=(q-1)/2$ is odd, either the trivial
    or the quadratic character of $H^F$ has odd multiplicity.

    If $q>3$, the group $G^F$ is perfect, so $\det\omega_\psi^-=1$.
    The inverse pairs contribute trivially to the determinant on $H^F$.
    Hence the quadratic character has even multiplicity, and the trivial
    character has odd multiplicity.
    If $q=3$, the one-dimensional representation $\omega_\psi^-$ factors
    through $\SL_2(3)^{\opr{ab}}\cong C_3$.
    Its restriction to the group $H^F$ of order four is therefore trivial.
    In either case, $\omega_\psi^-$ is $H^F$-distinguished.

    The $U^F$-weights of $(\omega_\psi^-)^\vee$ are
    $\{\psi_{-x^2}\mid x\in R\}$.
    Since $-1$ is not a square in $\bF_q$, these weights are disjoint
    from those of $\omega_\psi^-$.
    Thus $\omega_\psi^-\not\simeq(\omega_\psi^-)^\vee$.
\end{proof}

We use parabolic induction to extend
Lemma~\ref{lem:lp-counterexample-sl2} to higher rank.

\begin{proposition}
    \label{prop:lp-counterexample-disconnected-center}
    Let $G=\Sp_{4n+2}$ over $\bF_q$ with $n\geq1$ and
    $q\equiv3\pmod4$, and let $F$ be its Frobenius endomorphism.
    Fix a symplectic basis
    $(x_1,x_1^*,\ldots,x_{2n+1},x_{2n+1}^*)$ over $\bF_q$.
    Let $\theta$ be conjugation by
    $\opr{diag}(\overbrace{s,\ldots,s}^{2n+1})$ in this basis,
    where $s\coloneqq\left(\begin{smallmatrix}
                0 & -1 \\ 1 & 0
            \end{smallmatrix}\right)$.
    Set $H\coloneqq G^\theta$.
    There exists an $H^F$-distinguished non-self-dual irreducible
    character of $G^F$.
\end{proposition}
\begin{proof}
    Let $P$ stabilize the isotropic flag
    $\langle x_1\rangle\subset\cdots\subset\langle x_1,\ldots,x_{2n}\rangle$.
    Its standard Levi subgroup is
    $L\coloneqq\GL_1^{2n}\times\SL_2$.
    The involution $\theta$ sends $P$ to its opposite parabolic $P^-$
    with common Levi $L$.
    Consequently,
    \[
        H\cap P=L^\theta=\mu_2^{2n}\times T_{\opr{ns}},
    \]
    where $T_{\opr{ns}}$ lies in the last $\SL_2$ factor.

    Fix a nontrivial additive character $\psi$ of $\bF_q$ and set
    $\tau_0\coloneqq1^{\boxtimes2n}\boxtimes\omega_\psi^-$ on $L$.
    Let $\rho$ be the parabolic induction of $\tau_0$ along $P$.
    The identity double coset in Mackey's formula and
    Lemma~\ref{lem:lp-counterexample-sl2} give $\rho^{H^F}\ne0$.
    Thus $\rho$ has an $H^F$-distinguished irreducible constituent $\rho_0$.

    The representation $\tau_0$ is cuspidal, so $\rho_0$ has cuspidal
    support $(L,\tau_0)$, up to $G^F$-conjugacy.
    Every $g\in N_G(L)^F$ preserves the symplectic plane supporting
    the $\SL_2$ factor and restricts to an element of $\Sp_2(q)=\SL_2(q)$.
    Thus $\tau_0\circ\opr{Int}_g\simeq\tau_0$ for every such $g$.
    Since $\tau_0\not\simeq\tau_0^\vee$, the pairs $(L,\tau_0)$ and
    $(L,\tau_0^\vee)$ are not $G^F$-conjugate.
    Since duality commutes with parabolic induction, uniqueness of
    cuspidal support gives
    $\rho_0\not\simeq\rho_0^\vee$.
\end{proof}
\begin{remark}
    \label{rem:lp-counterexample-sl2-q-one-mod-four}
    Let $G=\SL_2$ over $\bF_q$ with $q\equiv1\pmod4$,
    and let $F$ be its Frobenius endomorphism.
    Fix a nontrivial additive character $\psi$ of $\bF_q$.
    Choose $\kappa\in\bF_q^\times\setminus(\bF_q^\times)^2$ and put
    $s\coloneqq\left(\begin{smallmatrix}0&\kappa\\1&0\end{smallmatrix}\right)$,
    $\theta\coloneqq\opr{Int}_s$, and $H\coloneqq G^\theta$.
    The analogue of Lemma~\ref{lem:lp-counterexample-sl2} uses
    the even Weil representation $\omega_\psi^+$ of
    \cite[\S2.2]{liuZhang2018fourierJacobi} in place of $\omega_\psi^-$.
    A similar argument shows that $\omega_\psi^+$ is irreducible and $H^F$-distinguished, and
    $\omega_\psi^+\not\cong(\omega_\psi^+)^\vee\circ\theta$.

    Here $\theta$ is outer on $G^F$, and $\omega_\psi^+$ is self-dual
    (\cite[Lemma~2.5]{liuZhang2018fourierJacobi}) and non-cuspidal.
    The cuspidal-support argument of
    Proposition~\ref{prop:lp-counterexample-disconnected-center}
    therefore does not extend to this case.
\end{remark}

\end{document}